\documentclass[notitlepage,leqno,10pt]{article}

\usepackage{latexsym}
\usepackage{amsmath}
\usepackage{amsfonts}
\usepackage{amssymb}
\usepackage{graphicx}
\usepackage{hyperref}
\renewcommand{\a }{\alpha }

\renewcommand{\d}{\delta }
\newcommand{\D }{\Delta }

\newcommand{\e }{\varepsilon }
\newcommand{\g }{\gamma}

\newcommand{\G }{\Gamma }

\newcommand{\n }{\nabla }

\newcommand{\s }{\sigma }

\renewcommand{\th }{\theta }
\renewcommand{\o }{\omega }
\renewcommand{\O }{\Omega }

\newcommand{\z }{\zeta }
\newcommand{\ov}{\overline}
\newcommand{\pa}{\partial}
\newcommand{\tPhi}{\tilde{\Phi}}

\newcommand{\be}{\begin{equation}}
\newcommand{\ee}{\end{equation}}
\newenvironment{pf}{\noindent{\sc Proof}.\enspace}{\rule{2mm}{2mm}\medskip}

\newtheorem{remark}{Remark}[section]
\newcommand{\R}{\mathbb{R}}

\newcommand{\N}{\mathbb{N}}

\newcommand{\no}{\noindent}

\author{Fethi MAHMOUDI$^{\rm a}$ \and Andrea MALCHIODI$^{\rm a}$ \and Marcelo MONTENEGRO$^{\rm b}$}

\date{}

\title{Solutions to the nonlinear Schr\"odinger equation \\
carrying momentum along a curve. \\ Part I: study of the limit set
and approximate solutions}

\begin{document}

\newtheorem{lem}{Lemma}[section]
\newtheorem{pro}[lem]{Proposition}
\newtheorem{thm}[lem]{Theorem}
\newtheorem{rem}[lem]{Remark}
\newtheorem{cor}[lem]{Corollary}
\newtheorem{df}[lem]{Definition}

\maketitle

\begin{center}

$^{\rm a}${\small SISSA, Sector of Mathematical Analysis \\ Via
Beirut 2-4, 34014 Trieste, Italy}

$^{\rm b}${Universidade Estadual de Campinas, IMECC, Departamento de
Matem\'atica, \\ Caixa Postal 6065, CEP 13083-970, Campinas, SP,
Brasil}

\end{center}

\footnotetext[1]{E-mail addresses: mahmoudi@ssissa.it
(F.Mahmoudi), malchiod@sissa.it (A. Malchiodi), msm@ime.unicamp.br
(M. Montenegro)}

\noindent {\sc abstract}. We prove existence of a special class of
solutions to the (elliptic) Nonlinear Schr\"odinger Equation $- \e^2
\D \psi + V(x) \psi = |\psi|^{p-1} \psi$, on a manifold or in the
Euclidean space. Here $V$ represents the potential, $p$ an exponent
greater than $1$ and $\e$ a small parameter corresponding to the
Planck constant. As $\e$ tends to zero (namely in the semiclassical
limit) we prove existence of complex-valued solutions which
concentrate along closed curves, and whose phase is highly
oscillatory. Physically, these solutions carry quantum-mechanical
momentum along the limit curves. In this first part we provide the
characterization of the limit set, with natural stationarity and
non-degeneracy conditions. We then construct an approximate solution
up to order $\e^2$, showing that these conditions appear naturally
in a Taylor expansion of the equation in powers of $\e$. Based on
these, an existence result will be proved in the second part
\cite{mmm2}.

\begin{center}

\bigskip\bigskip


\bigskip

\centerline{\bf AMS subject classification: 34B18, 35B25, 35B34,
35J20, 35J60}

\end{center}

\section{Introduction}\label{s:i}

In this paper we are concerned with concentration phenomena for
solutions of the singularly-perturbed elliptic equation
\begin{equation}\label{eq:pe}\tag{$NLS_\e$}
    - \e^2 \D_g \psi + V(x) \psi = |\psi|^{p-1} \psi \qquad
    \hbox{ on } M,
\end{equation}
where $M$ is an $n$-dimensional compact manifold (or the flat
Euclidean space $\R^n$), $V$ a smooth positive function on $M$
satisfying the properties
\begin{equation}\label{eq:VV}
    0 < V_1 \leq V \leq V_2; \qquad \qquad \|V\|_{C^3} \leq V_3,
\end{equation}
(for some fixed constants $V_1, V_2, V_3$) $\psi$ a complex-valued
function, $\e
> 0$ a small parameter and $p$ an exponent greater than $1$. Here
$\D_g$ stands for the Laplace-Beltrami operator on $(M,g)$.

\eqref{eq:pe} arises from the study of the Nonlinear Schr\"odinger
Equation
\begin{equation}\label{eq:nlse}
    i \hbar \frac{\pa \tilde{\psi}}{\pa t} = - \hbar^2 \D \tilde{\psi}
    + V(x) \tilde{\psi} - |\tilde{\psi}|^{p-1} \tilde{\psi} \qquad \quad
    \hbox{ on } M \times [0, + \infty),
\end{equation}
where $\tilde{\psi} = \tilde{\psi}(x,t)$ is the {\em wave
function}, $V(x)$ a potential, and $\hbar$ the {\em Planck
constant}. A special class of solutions to \eqref{eq:nlse} is
constituted by the functions whose dependence on the variables $x$
and $t$ is of the form $\tilde{\psi}(x,t) = e^{- i \frac{\o
t}{\hbar}} \psi(x)$. Such solutions are called {\em standing
waves} and up to substituting $V(x)$ with $V(x) - \o$, they give
rise to solutions of \eqref{eq:pe}, for $\e = \hbar$.

An interesting case is the {\em semiclassical limit} $\e \to 0$,
where one should expect to recover the Newton law of classical
mechanics. In particular, near stationary points of the potential,
one is lead to search highly concentrated solutions, which could
mimic point-particles at rest.

In recent years, a lot of attention has been devoted to the study of
the above problem: one of the first results in this direction is due
to Floer and Weinstein in \cite{fw}, where the case of $M = \R$ and
$p = 3$ is considered, and where existence of solutions highly
concentrated near critical points of $V$ has been proved. This
result has then been extended by Oh, \cite{o}, to the case of $\R^n$
for arbitrary $n$, provided $1 < p < \frac{n+2}{n-2}$. The profile
of these solutions is given by the {\em ground state} $U_{x_0}$
(namely the solution with minimal energy, which is real-valued,
everywhere positive and can be assumed radial) of the following {\em
limit equation}
\begin{equation}\label{eq:vxo}
    - \D u + V(x_0) u = u^p \qquad \quad \hbox{ in } \R^n,
\end{equation}
where $x_0$ is the concentration point. The solutions $u_\e$
obtained in the aforementioned papers behave qualitatively like
$u_\e(x) \simeq U_{x_0} \left( \frac{x-x_0}{\e} \right)$ as $\e$
tends to zero, and since $U_{x_0}$ decays exponentially to zero at
infinity, $u_\e$ vanishes rapidly away from $x_0$.

Two comments are in order: first of all the criticality of $V$ at
$x_0$ is a {\em necessary condition} for such a behavior, as shown
in \cite{wx}. Secondly, as pointed out in \cite{cp}, also the
upper bound $p < \frac{n+2}{n-2}$ is required for having solutions
concentrating at points: indeed the well-known {\em Pohozaev's
identity} imposes this restriction for having existence of
solutions to \eqref{eq:vxo} tending to zero at infinity.

The above existence results have been extended in several
directions, including the construction of solutions with multiple
peaks, the case of degenerate potentials, potentials tending to zero
at infinity and more general nonlinearities. We refer the interested
reader for example to the (incomplete) list of works \cite{abc},
\cite{afm}, \cite{am}, \cite{amr},  \cite{ams}, \cite{as},
\cite{bw}, \cite{df}, \cite{gr}, \cite{jt} and to the bibliographies
therein.

We also mention the mathematical similarities between \eqref{eq:pe}
and the following problem
\begin{equation}\tag{$P_\e$}\label{eq:pne}
  \begin{cases}
    -\e^2 \D u + u = u^p & \text{ in } \O, \\
    \frac{\partial u}{\partial \nu} = 0 & \text{ on } \partial \O,
    \\ u > 0 & \text{ in } \O,
  \end{cases}
\end{equation}
where $\O$ is a smooth bounded domain of $\R^N$, $p > 1$, and where
$\nu$ denotes the exterior unit normal vector to $\partial \O$.
Problem \eqref{eq:pne} arises in the study of some biological
models, see for example \cite{ni} and references therein, and as
\eqref{eq:pe} it exhibits concentration of solutions at some points
of $\ov{\O}$. Since the last equation is homogeneous, the location
of the concentration points is determined by the geometry of the
domain: if it occurs at the boundary, these are critical points of
the mean curvature while if it occurs at the interior these points
are (roughly) singular points for the distance function from the
boundary. About this topic, we refer the reader to \cite{dy},
\cite{dfw}, \cite{gr2}, \cite{gui}, \cite{guiw2}, \cite{guiww},
\cite{liyy}, \cite{ln}, \cite{lnt}, \cite{nt91}, \cite{nt93},
\cite{nty}, \cite{w97}.

\

More recently, new types of solutions to \eqref{eq:pe} have been
found, since when $\e$ tends to zero they do not concentrate at
points, but instead at sets of higher dimension. Before stating our
main result, it is convenient to recall the progress on this topic
and to illustrate the new phenomena involved. Some first results in
this direction were given in \cite{bd}, \cite{bed} in the case of
radial symmetry, and later improved in \cite{amn1} (see also
\cite{amn2} for the problem in bounded domains), where necessary and
sufficient conditions for the location of the concentration set have
been given. Differently from the previous case, the limit set is not
anymore stationary for the potential $V$: indeed, from heuristic
considerations, the {\em energy} of a solution concentrated near a
sphere of radius $r$ depends both on $V$ and on its {\em volume},
which is proportional to $\e r^{n-1}$. In \cite{amn1} it was shown
that the candidate radii of concentration are the critical points of
the function $r^{n-1} V(r)^{\frac{p+1}{p-1} - \frac 12}$ (the power
of $V$ in this formula arises from some scaling argument, related to
the dependence in $V(x_0)$ of the solutions to \eqref{eq:vxo}, see
also Section \ref{s:not}). Furthermore, no upper bound on the
exponent $p$ is required: in fact the profile of these solutions is
given by the solution of \eqref{eq:vxo} in $\R^1$, and in one
dimension there is no restriction for the existence of entire
solutions.

Based on the above energy considerations, in \cite{amn1} it is
also stated a conjecture concerning concentration on
$k$-dimensional manifolds, for $k = 1, \dots, n-1$: it is indeed
expected that, under suitable non-degeneracy assumptions, the
limit set should satisfy the equation
\begin{equation}\label{eq:limitsig}
    \th_k \n^N V = V {\bf H}, \qquad \qquad \hbox{ with } \quad
    \th_k = \frac{p+1}{p-1} - \frac 12 (n-k),
\end{equation}
where $\n^N$ stands for the normal gradient, ${\bf H}$ the curvature
vector, and the profile of the solutions at a point $x_0$ in the
limit set should be asymptotic, in the normal directions, to the
ground state of
\begin{equation}\label{eq:vxorn-k}
    - \D u + V(x_0) u = u^p \qquad \qquad \hbox{ in } \R^{n-k}.
\end{equation}
Since the Pohozaev identity implies $p < \frac{n-k+2}{n-k-2}$ for
the existence of non trivial solutions, the latter condition is
expected to be a natural one for dealing with this phenomenon.

Actually, concerning \eqref{eq:pne} another conjecture has been
previously stated, asserting existence of solutions concentrating at
sets of positive dimension. Concerning the latter problem, starting
from the paper \cite{malm}, there has been some progress in the
general setting (without symmetry assumptions), and after the works
\cite{mm}, \cite{mal}, \cite{malm2}, existence in now known for
arbitrary dimension and codimension. About problem \eqref{eq:pe},
the conjecture in \cite{amn1} has been verified in \cite{dkw} for $n
= 2$ and $k = 1$. Some other (and related) results, under some
reduced symmetry assumptions (as cylindrical or similar) have been
given in \cite{bape}, \cite{dy2}, \cite{mopa}, \cite{shi}.
Specifically we mention the note \cite{wein}: instead of
\eqref{eq:nlse}, it is considered there the nonlinear wave equation
\begin{equation}\label{eq:nlwe}
    - \frac{\pa^2 \tilde{\psi}}{\pa t^2} = - \D \tilde{\psi}
    + V(x) \tilde{\psi} - |\tilde{\psi}|^{p-1} \tilde{\psi} \qquad \quad
    \hbox{ on } M \times [0, + \infty)
\end{equation}
for $p=3$. It has been proved in \cite{a1} and \cite{a2} that
$$
    - \frac{\pa^2 \tilde{\psi}}{\pa t^2} = - \D \tilde{\psi}
    + V(x) \tilde{\psi}
$$
has solutions which remain concentrated near  elliptic closed
geodesics in $M$ for long periods of time, but which eventually
drift away when $t \to \infty$. A. Weinstein in \cite{wein} proved
that there do exist periodic solutions of \eqref{eq:nlwe} remaining
concentrated for all times, whenever $p=3$, $M = S^2$ and  $V$ is
odd.

It is worth pointing out a major difference between the symmetric
and the non-symmetric situation. In fact, since the ground states
of \eqref{eq:vxo} or \eqref{eq:vxorn-k} are of mountain-pass type
(namely critical points of some Euler functional with Morse index
equal to $1$), equation \eqref{eq:pe} becomes highly resonant. To
explain the reason, we consider for example a real-valued function
$\psi$ in $\R^2$ with a radial potential. One can begin by finding
approximate (radial) solutions of the form $u_{\ov{r},\e}(r)
\simeq U_{\ov{r}}\left( \frac{r-\ov{r}}{\e}\right)$, where
$U_{\ov{r}}$ is the solution of \eqref{eq:vxo} for $n = 1$
corresponding to $V(\ov{r})$. Then, with a good choice of
$\ov{r}$, one can try to linearize the equation and find true
solutions via the implicit function theorem. The linearized
equation, taking $\psi$ real for simplicity, becomes
\begin{equation*}
    - \e^2 \D \psi + V(r) \psi - p u_{\ov{r},\e}(r)^{p-1} \psi \qquad \qquad
  \hbox{ in } \R^2.
\end{equation*}
Using polar coordinates $(r,\vartheta)$ and a Fourier decomposition
of $\psi$ with respect to $\vartheta$, $\psi(r,\vartheta) = \sum_j
e^{i j \vartheta} \psi_j(r)$, we see that on each component $\psi_j$
acts the operator
\begin{equation}\label{eq:psik}
    \underbrace{- \e^2 \psi''_j - \e^2 \frac 1r \psi'_j + V(r)
    \psi_j - p u_{\ov{r},\e}(r)^{p-1} \psi_j}_{L_{1,\e} \psi_j}
    + \frac{1}{r^2} \e^ 2 j^2 \psi_j \qquad \quad \hbox{ on }
    [0,+\infty),
\end{equation}
where $L_{1,\e}$ (apart from the term $\e^2 \frac 1r \psi'_j$ which
is not relevant to the next discussion) represents the linearized
equation of \eqref{eq:pe} in one dimension near a soliton. Since one
expects to deal with functions which are highly concentrated near $r
= \ov{r}$, the last term in the above formula naively {\em
increases} the eigenvalues by a quantity of order
$\frac{1}{\ov{r}^2} \e ^2 j^2$ compared to those of $L_{1,\e}$.

The operator $L_{1,\e}$ possesses a negative eigenvalue $\eta_\e$
lying between two negative constants independent of $\e$ (since
$U_{\ov{r}}$ is of mountain-pass type, as explained before) and a
(nearly) zero eigenvalue $\s_\e$, by the translation invariance of
\eqref{eq:vxo} in $\R^1$. As a consequence, the operator in
\eqref{eq:psik} will possess two sequences of eigenvalues
qualitatively of the form $\eta_{j,\e} \simeq \eta_\e + \e^2 j^2$
and $\s_{j,\e} \simeq \s_\e + \e^2 j^2$. This might generate two
kinds of resonances: for small values of $j$, when $\s_{j,\e} \simeq
0$, and for $j$ of order $\frac 1 \e$, when $\eta_{j,\e}$ could be
close to zero. A comment is in order on the corresponding
eigenfunctions, which can be roughly studied with a separation of
variables as before. The ones relative to $\s_{j,\e}$ (for $j$
small) are slowly oscillating along the limit set, while the ones
relative to the resonant $\eta_{j,\e}$'s are fast oscillating with a
number of oscillations proportional to $k \simeq \frac 1 \e$.

The invertibility of the linearized operator will then be equivalent
to having all the $\s_{j,\e}$'s and all the $\eta_{j,\e}$'s
different from zero. A control on the resonant $\s_{j,\e}$'s can be
obtained (via some careful expansions) from a suitable
non-degeneracy condition involving the limit set and the potential
$V$. In \cite{amn1} for example, this can be achieved from the fact
of having a non-degenerate critical point of the function $r^{n-1}
V(r)^{\frac{p+1}{p-1} - \frac 12}$. On the other hand, the possible
vanishing of some $\eta_{j,\e}$ is peculiar of this concentration
behavior and more intrinsic, so invertibility can only be achieved
by choosing suitable values of $\e$. It is interesting to compare
this phenomenon (which is also present in \eqref{eq:pne}) to a
result in \cite{dan2}, asserting that if the Morse index of a family
of solutions to \eqref{eq:pne} stays bounded as $\e \to 0$, these
must concentrate at a finite number of points.

These formal considerations can also apply to the case of
concentration near a general manifold (without symmetry) in higher
dimension or codimension. Instead of expanding in polar coordinates
one can use (naively) a Fourier decomposition with respect to the
eigenfunctions of the Laplace-Beltrami operator and the normal
Laplacian on the limit manifold, see \cite{mm}. If the latter has
dimension $k$ then the term $\e^2 j^2$, by the Weyl's asymptotic
formula, has to be substituted with a quantity behaving like $\e^2
j^{\frac 2k}$. We notice that in this way the average distance
between two consecutive $\eta_{j,\e}$'s (when they are close to
zero) is of order $\e^k$ so, even if we have invertibility, the
distance of the spectrum to zero is (in the best cases) of order
$\min \left\{ \min_j |\eta_{j,\e}|, \min_j |\s_{j,\e}| \right\}
\simeq \min\{\e^2, \e^k\}$. Therefore the inverse operator is always
large in norm. By this reason, to apply the implicit function
theorem we need first to find good approximate solutions, with a
precision depending on $k$, and then prove indeed that the
linearized operator is invertible for suitable values of $\e$. This
is indeed a rather delicate issue: for reasons of brevity we do not
discuss it here but we refer directly to \cite{malm} and
\cite{malm2}. Related phenomena appear in some geometric problems as
well, dealing with the construction of surfaces with constant mean
curvature, see \cite{mmp}, \cite{mp}.

When $\O$ is a radially symmetric domain and the potential $V$ is
radially symmetric, the problem is simpler, since working in spaces
of invariant functions avoids most of the above resonances. In this
case only finitely-many eigenvalues (depending on the dimension and
the codimension of the concentration set) can approach zero, and
localization can be determined with a finite-dimensional reduction
of the problem.

\

In this paper and in \cite{mmm2} we construct a new type of
solutions, which concentrate along some curve $\g$, and which
physically carry momentum along the limit set. Differently from the
solutions discussed before, these are complex-valued and their
profile near any point $x_0$ in the image of $\g$ is asymptotic to a
solution to \eqref{eq:vxo} which decays exponentially to zero away
from the $x_n$ axis of $\R^n$ and is periodic in $x_n$. More
precisely, we consider solutions of the form
$$
  \phi(x',x_n) = e^{-i \hat{f} x_n} \hat{U}(x'), \qquad \quad
  x' = (x_1, \dots, x_{n-1}),
$$
where $\hat{f}$ is some constant and $\hat{U}(x')$ a real function.
With this choice of $\phi$, the function $\hat{U}$ satisfies
\begin{equation}\label{eq:ovU}
    - \D \hat{U} + \left( \hat{f}^2 + V(x_0) \right) \hat{U}
    = |\hat{U}|^{p-1} \hat{U} \qquad \quad \hbox{ in } \R^{n-1},
\end{equation}
and decays to zero at infinity. Solutions to \eqref{eq:ovU} can be
found by considering the (real) function $U$ satisfying $- \D U +
U = U^p$ in $\R^{n-1}$ (decaying to zero at infinity), and by
using the scaling
\begin{equation}\label{eq:ovhovk}
    \hat{U}(x') = \hat{h} U(\hat{k} x'), \qquad \qquad \hat{h} = \left(
  \hat{f}^2 + V(x_0) \right)^{\frac{1}{p-1}}, \quad \hat{k} = \left(
  \hat{f}^2 + V(x_0) \right)^{\frac 12}.
\end{equation}
In the above formulas $\hat{f}$ can be taken arbitrarily, and
$\hat{h}, \hat{k}$ have to be chosen accordingly, depending on
$V(x_0)$. The constant $\hat{f}$ represents the speed of the phase
oscillation, and is physically related to the velocity of the
quantum-mechanical particle represented by the wave function.

We are aware of one result only in this direction, given in
\cite{dap}, where the case of an axially-symmetric potential is
considered, and our goal here is to treat this phenomenon in a
generic situation, without any symmetry restriction. Some of the
difficulties of such an extension were naively summarized in the
above discussion but some new ones arise, due to the fact that the
standing waves are complex-valued, and due to their highly
oscillatory phase. In this first part we determine the concentration
set and show that its geometric characterization appears when we
construct approximate solutions to \eqref{eq:pe}, while a full
existence result will be given in \cite{mmm2}.

Before stating our main result we discuss how to determine the limit
set: if we look for a solution $\psi$ to \eqref{eq:pe} with the
above profile, then it should qualitatively behave as
\begin{equation}\label{eq:prof}
    \psi(\ov{s},\z) \simeq e^{- i \frac{f(\ov{s})}{\e}}
    h(\ov{s}) U\left(\frac{k(\ov{s}) \z}{\e}\right),
\end{equation}
where $\ov{s}$ stands for the arc-length parameter of $\g$, and $\z$
for a system of geodesic coordinates normal to $\g$. For having more
flexibility, we chose the phase oscillation to depend on the point
$\g(\ov{s})$, while $h(\ov{s}), k(\ov{s})$ should satisfy
\begin{equation}\label{eq:ovhovk3}
    h(\ov{s}) = \left(
  (f'(\ov{s}))^2 + V(\g(\ov{s})) \right)^{\frac{1}{p-1}}, \qquad
   \qquad k(\ov{s}) = \left(
  (f'(\ov{s}))^2 + V(\g(\ov{s})) \right)^{\frac 12},
\end{equation}
which is the counterpart of \eqref{eq:ovhovk} for a variable
potential.

The function $f(\ov{s})$ can be (heuristically) determined using
an expansion of \eqref{eq:pe} at order $\e$: a computation
performed in Subsection \ref{ss:cand} (see in particular formula
\eqref{eq:f'C}) shows that
\begin{equation}\label{eq:f'Cintr}
    f'(\ov{s}) \simeq \mathcal{A} h^{\s}(\ov{s}) \qquad
    \quad \hbox{ with } \quad \s = \frac{(n-1)(p-1)}{2}-2,
\end{equation}
where $\mathcal{A}$ is an arbitrary constant. At this point, only
the curve $\g$ should be determined. First of all, we notice that
the phase should be a periodic function in the length of the curve,
and therefore by \eqref{eq:f'Cintr} it is natural to work in the
class of loops
\begin{equation}\label{eq:Gamma}
    \G := \left\{ \g : \R \to M \hbox{ periodic} \; : \;
  \mathcal{A} \int_\g h(\ov{s})^\s d \ov{s} = \hbox{constant}
  \right\},
\end{equation}
where $\ov{s}$ stands for the arc-length parameter on $\g$.
Problem \eqref{eq:pe} has a variational structure, with
Euler-Lagrange functional given by
$$
  E_\e(\psi) = \frac 12 \int_M \left( \e^2 |\n_g \psi|^2 + V(x)
  |\psi|^2 \right) - \frac{1}{p+1} \int_M |\psi|^{p+1}.
$$
For a function of the form \eqref{eq:prof}, by a scaling argument
(see \eqref{eq:appen}) one has
\begin{equation}\label{eq:reden}
    E_\e(\psi) \simeq \e^{n-1} \int_\g h(\ov{s})^\th d \ov{s}, \qquad
\quad \hbox{ with } \quad \th = p+1 - \frac 12 (p-1)(n-1),
\end{equation}
therefore a limit curve $\g$ should be a critical point of the
functional $\int_\g h(\ov{s})^\th d \ov{s}$ in the class $\G$.
With a direct computation, see Subsection \ref{ss:euler}, one can
check that the extremality condition is the following
\begin{equation}\label{eq:eulerintr}
  \n^N V = {\bf H} \left( \frac{p-1}{\th}
  h^{p-1} - 2 \mathcal{A}^2 h^{2\s} \right)
\end{equation}
where, as before, $\n^N V$ represents the normal gradient of $V$ and
${\bf H}$ the curvature vector of $\g$. Similarly, via some long but
straightforward calculation, one can find a natural non-degeneracy
condition for stationary points, which is expressed by the
invertibility of the operator in \eqref{eq:2ndvarfin4} acting on the
normal sections to $\g$ (we refer the reader to Section \ref{s:not}
for the notation used in the formula). We notice that, since formula
\eqref{eq:f'Cintr} determines only the derivative of the phase,  to
obtain periodicity we need to introduce some nonlocal terms, see
 \eqref{eq:C'}. After these preliminaries, we are in position
to state our main result.

\begin{thm}\label{t:main}
Let $M$ be a compact $n$-dimensional manifold, let $V : M \to \R$ be
a smooth positive function and let $1 < p < \frac{n+1}{n-3}$. Let
$\g$ be a simple closed curve in $M$: then there exists a positive
constant $\mathcal{A}_0$, depending on $V|_\g$ and $p$ for which the
following holds. If $0 \leq \mathcal{A} < \mathcal{A}_0$, if $\g$
satisfies \eqref{eq:eulerintr} and the operator in
\eqref{eq:2ndvarfin4} is invertible on the normal sections of $\g$,
there is a sequence $\e_k \to 0$ such that problem $(NLS_{\e_k})$
possesses solutions $\psi_{\e_k}$ having the asymptotics in
\eqref{eq:prof}, with $f$ satisfying \eqref{eq:f'Cintr}.
\end{thm}

\begin{rem}\label{r:main} (a) The statement of Theorem \ref{t:main}
remains unchanged if we replace $M$ by $\R^n$ (or with an open
manifold asymptotically Euclidean at infinity) and we assume $V$
to be bounded between two positive constants and for which $\|\n^l
V\| \leq C_l$, $l = 1, 2, 3$ for some positive constants $C_l$.

(b) The restriction on the exponent $p$ is natural, since
\eqref{eq:ovU} admits solitons if and only if $p$ is subcritical
with respect to the dimension $n-1$.

(c) The smallness requirement on $\mathcal{A}$ is technical  and we
believe this condition could be relaxed. Anyway, for
$\frac{n+2}{n-2} \leq p < \frac{n+1}{n-3}$, $\mathcal{A}$ should
have an upper bound depending on $V$, to have solvability of both
\eqref{eq:ovhovk} and \eqref{eq:f'Cintr}. About this condition see
Remark \ref{r:solv} and Remark 2.7 in \cite{mmm2}.

(d) Apart from the assumption on $\mathcal{A}$, Theorem \ref{t:main}
improves the result in \cite{dap}. In fact, in addition to removing
the symmetry condition (which is the main issue), the
characterization of the limit set is explicit, the assumptions on
$V$ are purely local, and the upper bound on $p$ is sharp.

(e) The existence of solutions to \eqref{eq:pe} only for a suitable
sequence $\e_k \to 0$ is related to the resonance phenomenon
described above. The result can be extended to a sequence of
intervals in the parameter $\e$ approaching zero but, at least with
our proof, we do not expect to find existence for all the epsilon's.
\end{rem}

\

\noindent Taking $\mathcal{A} = 0$ (which implies $f' \equiv 0$),
from \eqref{eq:ovhovk3} it follows that $V = h^{p-1}$ and that
condition \eqref{eq:eulerintr} is equivalent to \eqref{eq:limitsig},
so as a consequence of our result we can prove the conjecture in
\cite{amn1} for $k = 1$, extending the result in \cite{dkw}.

\begin{cor}\label{c:intr}
Let $M$ be a compact Riemannian $n$-dimensional manifold with metric
$g$, let $V : M \to \R$ be a  function satisfying \eqref{eq:VV} and
let $1 < p < \frac{n+1}{n-3}$. Let $\g$ be a simple closed curve
which is a non-degenerate geodesic with respect to the weighted
metric $V^{\frac{p+1}{p-1}-\frac{n-1}{2}} g$. Then  there is a
sequence $\e_k \to 0$ such that problem $(NLS_{\e_k})$ possesses
real-valued solutions $\psi_{\e_k}$ concentrating near $\g$ as $j
\to + \infty$ and having the asymptotic behavior
$$
  \psi_{\e_k}(\ov{s},\z) \simeq V(\g(\ov{s}))^{\frac{1}{p-1}}
  U\left( \frac{V(\g(\ov{s}))^{\frac{1}{2}}}{\e_k} \z \right),
$$
where $\ov{s}$ stands for the arc-length parameter of $\g$, and
$\z$ for a system geodesic coordinates normal to $\g$.
\end{cor}

\noindent Corollary \ref{c:intr} gives also some criterion for the
applicability of Theorem \ref{t:main}: in fact, starting from a
non-degenerate geodesic in the weighted metric, via the implicit
function theorem for $\mathcal{A}$ sufficiently small one obtains a
curve for which \eqref{eq:eulerintr} and the invertibility of
\eqref{eq:2ndvarfin4} hold. In particular, when $V$ is constant, one
can start with non-degenerate close geodesics on $M$ in the ordinary
sense.

\

\noindent The full proof of Theorem \ref{t:main} will be given in
the second part of the paper, \cite{mmm2}. In this first part we
derive some formal expansions of equation \eqref{eq:pe} for $\psi$
of the form \eqref{eq:prof}, and check that the assumptions of the
theorem (the stationarity and the non-degeneracy conditions)
guarantee to find approximate solutions up to any powers of $\e$.
{\em At a formal level}, we consider expansions with coefficients
depending smoothly on the variable $\ov{s}$, and the only
obstructions to an iterative solvability of \eqref{eq:pe} are given
by the presence of a kernel in the linearization of \eqref{eq:ovU},
which is generated by the functions $(\partial_{x_j} U)$, $j = 1,
\dots, n-1$ and by $i U$: this kernel arises naturally from the
invariance of \eqref{eq:ovU} by translation in $\R^n$ and by complex
rotation. However, we can guarantee approximate solvability up to
any order provided $\g$ is stationary and non-degenerate: precisely,
we prove here the following weaker version of the above theorem.

\begin{thm}\label{t:approx} Suppose the assumptions of Theorem \ref{t:main} hold.
Then for any $m \in \N$ and for $\e$ small, there exists a function
$\psi_{\e,m}$ with the profile \eqref{eq:prof} and which solves
\eqref{eq:nlse} up to order $o(\e^m)$.
\end{thm}

\no As discussed before, also some fast-oscillating functions (along
$\g$) contribute to generate some resonance, but we do not discuss
this aspect here. Postponing the description of the rigorous proof
to the introduction of \cite{mmm2}, here we limit ourselves to
mention the main new difficulty compared to the results in
\cite{dkw}, \cite{mm} \cite{mal} and \cite{malm}. In our case the
solutions are complex-valued, and this causes an extra degeneracy in
the equation, due to its invariance under multiplication by a phase
factor. As a consequence, we have a further (infinite-dimensional)
approximate kernel, corresponding roughly to factor of $\psi_\e$ in
the form $e^{- i f_1(\ov{s})}$, for $f_1$ arbitrary. The correction
in the phase can also be determined by a formal expansion in powers
of $\e$ and, as for $f'$, we still obtain nonlocal terms. Also, when
expanding formally the solutions in powers of $\e$, the highly
oscillatory behavior of solutions generates an increasing number of
derivatives in $\ov{s}$: anyway in the first part, where formal
expansions are carried out, this aspect is not very relevant.

\

\noindent The plan of the paper is the following. In Section
\ref{s:not} we study the functional in \eqref{eq:reden} constrained
to the class of curves $\G$, and we determine the Euler-Lagrange
equation together with the non-degeneracy condition. In Section
\ref{s:as} approximate solutions to \eqref{eq:pe} up to order $\e$
are found, and the error terms of order $\e^2$ are displayed. The
functions we obtain are allowed to depend on a section $\Phi$ of the
normal bundle to $\g$ and on a scalar function $f_1$. These
correspond to some tilting of the approximate solution
perpendicularly to the limit set and to a variation of the phase, in
order to have more flexibility. In Section \ref{s:proj} we consider
the terms of order $\e^2$ and we choose $f_1$ and $\Phi$ so that
even the terms of order $\e^2$ in the expansion vanish: we then
arrive to the proof of Theorem \ref{t:approx}. Finally in Section
\ref{s:app} we collect some technical material, namely some integral
identities and the verification of \eqref{eq:2ndvarfin4}.

\

\noindent The results in this paper and in \cite{mmm2}, together
with the main ideas of the proofs, are briefly summarized in the
note \cite{mmmnote}.

\

\begin{center}
{\bf Notation and conventions}
\end{center}

\no  Dealing with coordinates, capital letters like $A, B, \dots$
will vary between $1$ and $n$ while indices like $j, l, \dots$ will
run between $2$ and $n$. The symbol $i$ will always stand for the
imaginary unit.

\

\no For summations, we  use the standard convention of summing
terms where repeated indices appear.

\

\no We will choose coordinates $(x_1,\cdots,x_n)$ near a curve $\g$
and we will parameterize $\g$ by arc-length letting $x_1=\ov{s}$.
Its dilation $\g_\e:=\frac1\e\g$ will be parameterized by
$s=\frac1\e\ov{s}$. The length of $\g$ is denoted by $L$.

\

\no For simplicity, a constant $C$ is allowed to vary from one
formula to another, also within the same line.

\

\no For a real positive variable $r$ and an integer $m$, $O(r^m)$
(resp. $o(r^m)$) will denote a complex-valued quantity for which
$\left| \frac{O(r^m)}{r^m} \right|$ remains bounded (resp. $\left|
\frac{o(r^m)}{r^m} \right|$ tends to zero) when $r$ tends to zero.
We might also write $o_\e(1)$ for a quantity which tends to zero as
$\e$ tends to zero.

\

\section{Study of the reduced functional}\label{s:not}

In this section we consider the functional in the right-hand side of
\eqref{eq:reden} defined on the set $\G$, representing the
approximate energy $E_\e$ of a function concentrated near $\g$ with
the profile \eqref{eq:prof}. We first introduce a convenient set of
coordinates near an arbitrary (smooth) closed curve in $M$. Then,
using these coordinates we write the Euler equation and the second
variation formula at a stationary point.

\subsection{Geometric preliminaries}\label{ss:coord}

In this Subsection we discuss some preliminary geometric facts,
referring for example to \cite{do2}, \cite{spi}. Given an arbitrary
simple closed curve $\g$ in $M$, we choose coordinates $x_1 \dots,
x_n$ near $\g$, called {\em Fermi coordinates} in the following way.
We let $x_1$ parameterize the curve $\g$ by arc-length. At some
point $q$ in the image of $\g$ we consider an orthonormal
$(n-1)$-tuple $(Y_2, \dots, Y_n)$ which form a basis for $N_q \g$,
the normal bundle of $\g$ at $q$. We extend the $Y_l$'s as vector
fields along $\g$ via parallel transport along the curve with
respect to the normal connection $\n^N$, namely by the condition
$\n^N_{\dot{\g}} Y_l = 0$ for $l = 2, \dots, n$.

Next we parameterize a point near $\g$ using the following
coordinates $(\ov{s},y) \in \R \times \R^{n-1}$
$$
  (\ov{s}, y_2, \dots, y_n) \mapsto  \exp_{\g(\ov{s})}(y_2 Y_2 + \dots + y_n
  Y_n),
$$
where $\exp_q$ is the exponential map in $M$ through the point
$q$. In this way, fixing $\ov{s}$, each curve $t \mapsto t y$, for
$y \in \R^{n-1} \setminus \{0\}$ and $t$ close to zero, is mapped
into a geodesic in $M$ passing through $\g(\ov{s})$.

Let us now define the vector fields $E_1 = \frac{\pa}{\pa \ov{s}}$
and $E_l = \frac{\pa}{\pa y_l}$ for $l = 2, \dots, n$. We notice
that on $\g$ each $E_l$ coincides with $Y_l$, while $E_1$ on $\g$ is
nothing but $\dot{\g}$. By our choice of coordinates it follows that
$\n_E E = 0$ on $\g$ for any vector field $E$ which is a linear
combination (with coefficients depending only on $\ov{s}$) of the
$E_j$'s, $j = 2, \dots, n$. In particular, for any $l,j = 2, \dots,
n$, and for any $\a \in \R$ we have $\n_{E_l + \a E_j}\left(E_l + \a
E_j\right) = 0$ on $\g$, which implies $\n_{E_l} E_j+\n_{E_j} E_l =
0$ for every $l,j = 2, \dots, n$. Using the fact that  $E_A$'s are
coordinate vectors for $A = 1, \dots, n$ and in particular $\n_{E_A}
E_B = \n_{E_B} E_A$ for all $A, B = 1, \dots, n$, we obtain that
$\n_{E_l} E_j=0$ for every $l,j = 2, \dots, n$. This immediately
yields
$$
  \partial_m g_{lj} = E_m \langle E_l, E_j \rangle =
  \langle \n_{E_m}E_l,E_j \rangle+\langle E_l, \n_{E_m}E_j
  \rangle=0 \qquad \hbox{ on} \,\,\g, \qquad l,j,m = 2, \dots, n.
$$
Moreover, still since the $E_A$'s are coordinate vectors for $A =
1, \dots, n$, we obtain
  \begin{eqnarray*}
    \partial_m g_{1j} & = & E_m \langle E_1, E_j \rangle =
  \langle \n_{E_m}E_1,E_j \rangle+\langle E_1, \n_{E_m}E_j \rangle \\
    & = & \langle \n_{E_1} E_m,E_j \rangle + \langle E_1, \n_{E_m}E_j
    \rangle = 0 \qquad \qquad \quad \hbox{ on} \,\,\g, \qquad m,j = 2, \dots, n.
\end{eqnarray*}
Here we used the fact that $\n_{E_1}^N E_m = 0$ on $\g$, namely
that $\n_{E_1} E_m$ has zero normal components.

If ${\bf H} = H^m E_m$ is the curvature vector of $\g$ (which is
normal to the curve), then one has $\langle \n_{E_1} E_m, E_1
\rangle = - H^m$ on $\g$, so we easily deduce that
\begin{equation}\label{eq:dmg11}
    \pa_m g_{11} = E_m \langle E_1, E_1 \rangle = 2 \langle
  \n_{E_1} E_m, E_1 \rangle = - 2 H^m \qquad \hbox{ on } \g.
\end{equation}
One can also prove that the components $R_{1m1j}$ of the curvature
tensor are given by
\begin{equation}\label{eq:ctc}
R_{1m1j}=-\frac12\partial^2_{jm}g_{11}+ H^m H^j.
\end{equation}
Indeed, we have~
\begin{eqnarray*}
-R_{1m1j}&=&\langle R(E_1,E_j)E_1,E_m  \rangle=\langle
\n_{E_1}\n_{E_j}E_1,E_m \rangle-\langle \n_{E_j}\n_{E_1}E_1,E_m
\rangle\\
&=&\langle \n_{E_1}\n_{E_j}E_1,E_m \rangle-E_j\langle
\n_{E_1}E_1,E_m \rangle -\langle \n_{E_1}E_1,\n_{E_j}E_m \rangle
\\
&=&\langle \n_{E_1}\n_{E_j}E_1,E_m \rangle-E_j\langle
\n_{E_1}E_1,E_m \rangle \\
&=&\langle \n_{E_1}\n_{E_j}E_1,E_m \rangle-E_j E_1\langle E_1,E_m
\rangle+E_j \langle E_1,\n_{E_1}E_m \rangle\\
&=&\langle \n_{E_1}\n_{E_j}E_1,E_m \rangle+E_j \langle E_1,\n_{E_m}E_1 \rangle\\
&=&E_1\langle \n_{E_j}E_1,E_m \rangle-\langle
\n_{E_j}E_1,\n_{E_1}E_m \rangle+\frac12\,E_jE_m\langle E_1,E_1\rangle\\
&=&\frac12\partial^2_{jm}g_{11}-\frac14 \partial_m
g_{11}\,\partial_j g_{11},
\end{eqnarray*}
where here we have used the above properties and the fact that
$$
\n_{E_j}E_1=\n_{E_1}E_j=\frac12 \partial_j g_{11}E_1.
$$
Using \eqref{eq:dmg11} and \eqref{eq:ctc}, the above discussion can
be summarized in the following result.

\begin{lem}\label{l:me} In the coordinates $(\ov{s},y)$, for $y$
close to zero the metric coefficients satisfy
\begin{eqnarray*}
    g_{11}(y) & = & 1 - 2 \sum_{m=2}^n H^m y_m +
    \frac 12  \sum_{m,l=2}^n (H^m H^j - R_{1m1j}|_{\g})
    y_m y_l +
O(|y|^3);\\
    g_{1j}(y) & = &
    \frac 12  \sum_{m,l=2}^n \partial^2_{ml}
    g_{1j}|_{\g} y_m y_l +  O(|y|^3);
\\
    g_{kj}(y) & = & \d_{kj} + \frac 12  \sum_{m,l=2}^n
    \partial^2_{ml} g_{kj}|_{\g} y_m y_l+  O(|y|^3).
\end{eqnarray*}
\end{lem}
The second derivatives $\pa^2_{ml} g_{1j}$ and $\pa^2_{ml} g_{kj}$
could be expressed in terms of the curvature tensor and the
curvature of $\g$ reasoning as for \eqref{eq:ctc}. However for our
purposes it is not necessary to have such a formula, so we leave the
expansion of these coefficients in a generic form.

\subsection{First and second variations of the length
functional}\label{ss:1st2ndvarlen}

We recall next the formulas for the variations of the length of a
curve with respect to normal displacements. We start with a
regular closed curve $\g$ in $M$ of length $L$, which we
parameterize by arc-length, using a parameter $\ov{s} \in [0,L]$.
Then we consider a two-parameter family of closed curves
$\g_{t_1,t_2} : [0,L]\rightarrow M$, for $t_1, t_2$ in a
neighborhood of $0$ in $\R$, such that $\g_{0,0} \equiv \g$. The
length $L(t_1, t_2)$ of $\g_{t_1,t_2}$ is given by
\begin{equation*}\label{eq:Lt1t2}
    L(t_1,t_2) = \int_{\g_{t_1,t_2}} dl = \int_0^L \langle\dot{\g}_{t_1,t_2},
    \dot{\g}_{t_1,t_2}\rangle^{\frac 12} d\ov{s},
\end{equation*}
where $d l$ is the arc-length parameter and $\dot{\g}_{t_1,t_2}$
stands for $\frac{d \g_{t_1,t_2}}{d\ov{s}}$. We also define the
vector fields $\mathcal{V}, \mathcal{W}$ along $\g_{t_1,t_2}$ as
$\mathcal{V} = \frac{\pa \g_{t_1,t_2}}{\pa t_1}$ and $\mathcal{W}
= \frac{\pa \g_{t_1,t_2}}{\pa t_2}$. In the above coordinates, the
vector fields $\mathcal{V}$ and $\mathcal{W}$ along $\g$ can be
written as
$$
  \mathcal{V} = \sum_{j=2}^n \mathcal{V}^j(\ov{s}) E_j; \qquad
  \qquad \mathcal{W} = \sum_{m=2}^n \mathcal{W}^m(\ov{s}) E_m.
$$

Differentiating $L(t_1,t_2)$ with respect to $t_1$ we find
\begin{equation}\label{eq:1stvarlgen}
    \frac{\partial L(t_1,t_2)}{\partial t_1} = - \int_0^L
    \frac{\langle\n_{\mathcal{V}} \dot{\g}_{t_1,t_2},
    \dot{\g}_{t_1,t_2}\rangle}{\langle \dot{\g}_{t_1,t_2},
  \dot{\g}_{t_1,t_2} \rangle^{\frac 12}} d\ov{s}.
\end{equation}
Using \eqref{eq:dmg11}, at $(t_1, t_2) = (0,0)$ we have
\begin{equation*}
    \langle\n_{\mathcal{V}}
  \dot{\g}_{t_1,t_2}, \dot{\g}_{t_1,t_2} \rangle = - \mathcal{V}^m H^m,
\end{equation*}
therefore we can write the variation of the length at $\g$ in the
following way
\begin{equation}\label{eq:varl}
    \frac{\partial L(t_1,t_2)}{\partial t_1}|_{(t_1,t_2) = (0,0)}
    = - \int_0^L \mathcal{V}^m H^m d\ov{s} = - \int_0^L \langle \mathcal{V},
    {\bf H} \rangle d\ov{s}.
\end{equation}
Using \eqref{eq:1stvarlgen} we can evaluate the second variation
of the length as
$$
  \frac{\partial^2 L(t_1,t_2)}{\partial t_1 \pa t_2}
  = \int_0^L \left[ \frac{\langle\n_\mathcal{W} \dot{\g}_{t_1,t_2}, \n_\mathcal{V}
  \dot{\g}_{t_1,t_2}\rangle + \langle\dot{\g}_{t_1,t_2},
  \n_\mathcal{W} \n_\mathcal{V} \dot{\g}_{t_1,t_2}\rangle}{\langle\dot{\g}_{t_1,t_2},
  \dot{\g}_{t_1,t_2}\rangle^\frac12} -
  \frac{\langle\dot{\g}_{t_1,t_2}, \n_\mathcal{V}
  \dot{\g}_{t_1,t_2}\rangle \langle\dot{\g}_{t_1,t_2}, \n_\mathcal{W}
  \dot{\g}_{t_1,t_2}\rangle}{\langle\dot{\g}_{t_1,t_2},\dot{\g}_{t_1,t_2}\rangle^\frac32} \right]
  d\ov{s},
$$
so at $(t_1, t_2) = (0,0)$ we find
\begin{eqnarray*}
  \frac{\partial^2 L(t_1,t_2)}{\partial t_1 \pa t_2}|_{(t_1,t_2)=(0,0)} & = &
  \int_0^L \left[ \langle\n_\mathcal{W} \dot{\g}, \n_\mathcal{V} \dot{\g}\rangle + \langle\dot{\g},
  \n_\mathcal{W} \n_\mathcal{V} \dot{\g}\rangle - \langle\dot{\g}, \n_\mathcal{V}
  \dot{\g}\rangle \langle\dot{\g}, \n_\mathcal{W}
  \dot{\g}\rangle \right] d\ov{s}.
\end{eqnarray*}
Using the definition of the Riemann tensor and the fact that
$\mathcal{V}$ and $\mathcal{W}$ are coordinate vector fields (so
that $[\mathcal{V}, \mathcal{W}] = 0$) the last formula yields
\begin{eqnarray*}
    \frac{\partial^2 L(t_1,t_2)}{\partial t_1 \pa t_2}|_{(t_1,t_2)=(0,0)}
     & = & \int_0^L \left[ \langle\n_{\dot{\g}} \mathcal{W}, \n_{\dot{\g}}
     \mathcal{V}\rangle + \langle\dot{\g}, \n_\mathcal{W} \n_{\dot{\g}} \mathcal{V}\rangle -
  \langle\dot{\g}, \n_{\dot{\g}} \mathcal{V}\rangle
  \langle\dot{\g}, \n_{\dot{\g}} \mathcal{W}\rangle \right] d\ov{s} \\ & = & \int_0^L
  \left[ \langle\n_{\dot{\g}} \mathcal{W},
  \n_{\dot{\g}} \mathcal{V}\rangle + \langle R(\mathcal{W},\dot{\g})\mathcal{V},\dot{\g}\rangle
  - \langle\dot{\g}, \n_{\dot{\g}}
  \mathcal{V} \rangle\langle\dot{\g}, \n_{\dot{\g}} \mathcal{W}\rangle
  \right] d\ov{s} \\ & - & \int_0^L
  \langle\n_{\dot{\g}} \dot{\g}, \n_\mathcal{W} \mathcal{V}\rangle d\ov{s}.
\end{eqnarray*}
Here, we have used the fact that $g(\dot{\g}, \n_\mathcal{W}
\n_{\dot{\g}} \mathcal{V})=\langle
R(\mathcal{W},\dot{\g})\mathcal{V},\dot{\g}\rangle+\dot{\g}
\langle\n_\mathcal{W} \mathcal{V},\dot{\g}\rangle
-\langle\n_{\dot{\g}} \dot{\g},
  \n_\mathcal{W} \mathcal{V}\rangle$ and $\int_0^L \dot{\g}\langle\n_\mathcal{W}
\mathcal{V},\dot{\g}\rangle d\ov{s} = 0$. Since $\n_{E_l} E_j = 0$
on $\g$ for $l, j = 2, \dots, n$, we have
$$
  \n_\mathcal{W} \mathcal{V} =
  \sum_{j,m = 2, \dots, n}
  \mathcal{W}^m \mathcal{V}^j \n_{E_m}E_j \qquad \qquad \Longrightarrow
  \qquad \qquad
  \int_0^L \langle\n_{\dot{\g}} \dot{\g}, \n_\mathcal{W}
  \mathcal{V}\rangle d \ov{s} = 0.
$$
Moreover, recalling \eqref{eq:dmg11} we obtain
$$
  \n_{\dot{\g}} \mathcal{V} = \sum_{j=2}^n \dot{\mathcal{V}}^j E_j +
  \sum_{j = 2}^n  \mathcal{V}^j \n_{E_1}E_j =
  \sum_{j=2}^n \dot{\mathcal{V}}^j E_j -
  \sum_{j=2}^n H^j \mathcal{V}^j E_1.
$$
This implies, at $\g$~
\begin{eqnarray*}
\langle\n_{\dot{\g}} \mathcal{W},
  \n_{\dot{\g}} \mathcal{V}\rangle + \langle R(\mathcal{W},\dot{\g})
  \mathcal{V},\dot{\g}\rangle - \langle\dot{\g},
  \n_{\dot{\g}}
  V\rangle\langle\dot{\g}, \n_{\dot{\g}}
  \mathcal{W}\rangle&=&\sum_{j=2}^n\dot{\mathcal{V}}^j\dot{\mathcal{W}}^j-
  \sum_{j,l=2}^n
  R_{1j1l} \mathcal{V}^j \mathcal{W}^l.
\end{eqnarray*}
In this way the second variation of the length at $\g$ becomes
\begin{equation}\label{eq:2ndvarl}
    \frac{\partial^2 L(t_1,t_2)}{\partial t_1 \pa t_2}|_{(t_1,t_2)=(0,0)} =
    \int_0^L \left(\sum_j^n \dot{\mathcal{V}}^j \dot{\mathcal{W}}^j -
    \sum_{j,l=2}^n R_{1j1l} \mathcal{V}^j \mathcal{W}^l\right)d\ov{s}.
\end{equation}

\subsection{Determining the phase factor}\label{ss:cand}

In this section we derive formally the asymptotic profile of the
solutions to \eqref{eq:pe} which concentrate near some curve $\g$,
and we determine some necessary conditions satisfied by the limit
curve. For doing this, using the coordinates $(\ov{s},y)$ introduced
in Subsection \ref{ss:coord}, we look for approximate solutions
$\psi(\ov{s},y)$ of \eqref{eq:pe} making the following {\em ansatz}
$$
  \psi(\ov{s},y) = e^{- i \frac{f(\ov{s})}{\e}} h(\ov{s}) U \left( \frac{k(\ov{s})
  y}{\e} \right), \qquad \qquad \ov{s} \in [0, L], \quad y \in \R^{n-1},
$$
where the function $U$ is the unique radial solution (see
\cite{bl}, \cite{gnn2}, \cite{kwo}, \cite{str}) of the problem
\begin{equation}\label{eq:ovv}
  \left\{
    \begin{array}{ll}
      - \D U + U = U^p & \hbox{ in } \R^{n-1}; \\
      U(y) \to 0 & \hbox{ as } |y| \to + \infty; \\
      U > 0 & \hbox{ in } \R^{n-1},
    \end{array}
  \right.
\end{equation}
and where the functions $f$, $h$ and $k$ are periodic on $[0,L]$ and
have to be determined. With some easy computations we obtain
$$
  \frac{\partial \psi}{\partial \ov{s}} = - \frac{i f'(\ov{s})}{\e} h(\ov{s}) U
  \left( \frac{k(\ov{s}) y}{\e} \right) e^{- i \frac{f(\ov{s})}{\e}} + e^{- i
  \frac{f(\ov{s})}{\e}} h'(\ov{s}) U \left( \frac{k(\ov{s}) y}{\e} \right) +
  e^{- i \frac{f(\ov{s})}{\e}} h(\ov{s}) k'(\ov{s}) \n_y U \left( \frac{k(\ov{s}) y}{\e}
  \right) \cdot \frac{y}{\e};
$$
\begin{eqnarray*}
    \frac{\partial^2 \psi}{\partial \ov{s}^2} & = & \left[ - i
    \frac{f''(\ov{s})}{\e} h(\ov{s}) U \left( \frac{k(\ov{s}) y}{\e} \right)
    - 2 i \frac{f'(\ov{s})}{\e} h'(\ov{s}) U \left( \frac{k(\ov{s}) y}{\e}
    \right) - 2 i \frac{f'(\ov{s})}{\e} h(\ov{s}) k'(\ov{s}) \n_y U \left(
    \frac{k(\ov{s}) y}{\e} \right) \cdot \frac{y}{\e} \right. \\ & - &
    \left. \frac{(f'(\ov{s}))^2}{\e^2} h(\ov{s}) U \left( \frac{k(\ov{s}) y}{\e}
    \right) + 2 h'(\ov{s}) k'(\ov{s}) \n U \left( \frac{k(\ov{s}) y}{\e}
    \right) \cdot \frac{y}{\e} + h(\ov{s}) k''(\ov{s}) \n_y U \left(
    \frac{k(\ov{s}) y}{\e} \right) \cdot \frac{y}{\e} \right. \\ & + &
    \left. h(\ov{s}) (k'(\ov{s}))^2 \n^2_y U \left( \frac{k(\ov{s}) y}{\e} \right)
    \left[ \frac y\e, \frac y\e \right]  + h''(\ov{s}) U \left( \frac{k(\ov{s})
    y}{\e} \right) \right] e^{- i \frac{f(\ov{s})}{\e}},
\end{eqnarray*}
and also
$$
  \D_y \psi (\ov{s},y) = \frac{(k(\ov{s}))^2}{\e^2} \D_y U \left( \frac{k(\ov{s}) y}{\e}
  \right)e^{- i
  \frac{f(\ov{s})}{\e}}h(\ov{s}).
$$
Since $U$ decays to zero at infinity (exponentially indeed, by the
results in \cite{gnn2}), and since the function $\psi$ is scaled of
order $\e$ near the curve $\g$, in a first approximation we can
assume the metric $g$ of $M$ to be flat in the coordinates $(\ov{s},
y)$, see the expansions in Lemma \ref{l:me}. We look now at the
leading terms in \eqref{eq:pe}, which are of order $1$. Since $-
\D_g \psi$ is multiplied by $\e^2$, we have to focus on the terms of
order $\frac{1}{\e^2}$ in the Laplacian of $\psi$. In the above
expressions of $\frac{\pa \psi}{\pa \ov{s}}$, $\frac{\pa^2 \psi}{\pa
\ov{s}^2}$ and $\D_y \psi$, we have that the function $U$ and its
derivatives are of order $1$ when $|y| = O(\e)$, therefore when the
variables $y$ appear as factors in these expressions, we consider
them to be of order $\e$. For example, $\n^2 U \left(
\frac{k(\ov{s}) y}{\e} \right) \left[ \frac y\e, \frac y\e \right]$
will be regarded as a term of order $1$.

With these criteria, using the above computations and assumptions,
imposing the leading terms in \eqref{eq:pe} to vanish we obtain
$$
  - k^2(\ov{s}) h(\ov{s}) \D_y U \left( \frac{k(\ov{s}) y}{\e} \right) +
  h(\ov{s}) \left[ V(\ov{s}) + (f'(\ov{s}))^2 \right] U \left( \frac{k(\ov{s}) y}{\e}
  \right) = h(\ov{s})^p U \left( \frac{k(\ov{s}) y}{\e} \right)^p.
$$
{From} \eqref{eq:ovv}, we have the two relations
\begin{equation}\label{eq:hk}
    k^2(\ov{s}) = h(\ov{s})^{p-1}; \qquad \qquad \left[ V(\ov{s}) + (f'(\ov{s}))^2
  \right] = k(\ov{s})^2 = h(\ov{s})^{p-1}.
\end{equation}
We next obtain an equation for $f$, which is derived looking at
the next-order expansion of \eqref{eq:pe}. The next coefficient
arises from the terms of order $\frac 1 \e$ in $- \D_g \psi$,
which are given by
$$
    i \left[
    \frac{f''(\ov{s})}{\e} h(\ov{s}) U \left( \frac{k(\ov{s}) y}{\e} \right)
    + 2 \frac{f'(\ov{s})}{\e} h'(\ov{s}) U \left( \frac{k(\ov{s}) y}{\e}
    \right) + 2 \frac{f'(\ov{s})}{\e} h(\ov{s}) k'(\ov{s}) \n_y U \left(
    \frac{k(\ov{s}) y}{\e} \right) \cdot \frac{y}{\e} \right] e^{- i \frac{f(\ov{s})}{\e}}.
$$
Multiplying this expression by $U\left( \frac{k(\ov{s})y}{{\e}}
\right)$ and integrating in $y \in \R^{n-1}$, imposing vanishing
of this integral as well gives
\begin{eqnarray*}
 0&=&  f''(\ov{s}) h(\ov{s}) \int_{\R^{n-1}} U^{2} \left( \frac{k(\ov{s}) y}{\e}
  \right) dy + 2 h'(\ov{s}) f'(\ov{s}) \int_{\R^{n-1}} U^{2} \left( \frac{k(\ov{s})
  y}{\e} \right) dy \\ &+& 2 f'(\ov{s}) h(\ov{s}) k'(\ov{s}) \int_{\R^{n-1}} U
  \left( \frac{k(\ov{s}) y}{\e} \right) \n_y U \left( \frac{k(\ov{s}) y}{\e}
  \right) \cdot \frac y \e\; dy.
\end{eqnarray*}
Integrating by parts and reasoning as for the usual Pohozaev's
identity we obtain that $f$ must satisfy
$$
  f''(\ov{s}) h(\ov{s}) + 2 f'(\ov{s}) h'(\ov{s}) - (n-1) f'(\ov{s}) h(\ov{s})
  \frac{k'(\ov{s})}{k(\ov{s})} = 0.
$$
This is solvable in $f'(\ov{s})$ and gives, for an arbitrary
constant $\mathcal{A}$
\begin{equation}\label{eq:f'C}
    f'(\ov{s}) = \mathcal{A} k(\ov{s})^{n-1} h(\ov{s})^{-2} = \mathcal{A}
  h(\ov{s})^{\frac{(n-1)(p-1)}{2} - 2},
\end{equation}
where we have used the above equation \eqref{eq:hk} for $k$. Now we
can solve the equation for $h(\ov{s})$ depending on the potential
$V(\ov{s})$ and the above constant $\mathcal{A}$. In fact, we get
that $h(\ov{s})$ should solve
\begin{equation}\label{eq:hV}
    V(\ov{s}) + \mathcal{A}^2 h(\ov{s})^{2\s} := V(\ov{s}) + \mathcal{A}^2
    h(\ov{s})^{(n-1)(p-1) - 4} =
    h(\ov{s})^{p-1},
\end{equation}
where we have set
\begin{equation}\label{eq:s}
    \s = \frac{(n-1)(p-1)}{2} - 2.
\end{equation}

\begin{rem}\label{r:solv}
We notice that, assuming $\mathcal{A}$ to be small enough (depending
on $V$ and $p$), the above equation is always solvable in
$h(\ov{s})$. More precisely, when $p < \frac{n+2}{n-2}$ (and hence
when $2 \s < p-1$), the solution is also unique. For $p \geq
\frac{n+2}{n-2}$ there might be a second solution. In this case, we
just consider the smallest one, which stays uniformly bounded (both
from above and below) when $\mathcal{A}$ is small enough, see
Figures \ref{fig1} and \ref{fig2} below.
\begin{figure}[h] \begin{center}
 \includegraphics[angle=270,width=5.5cm]{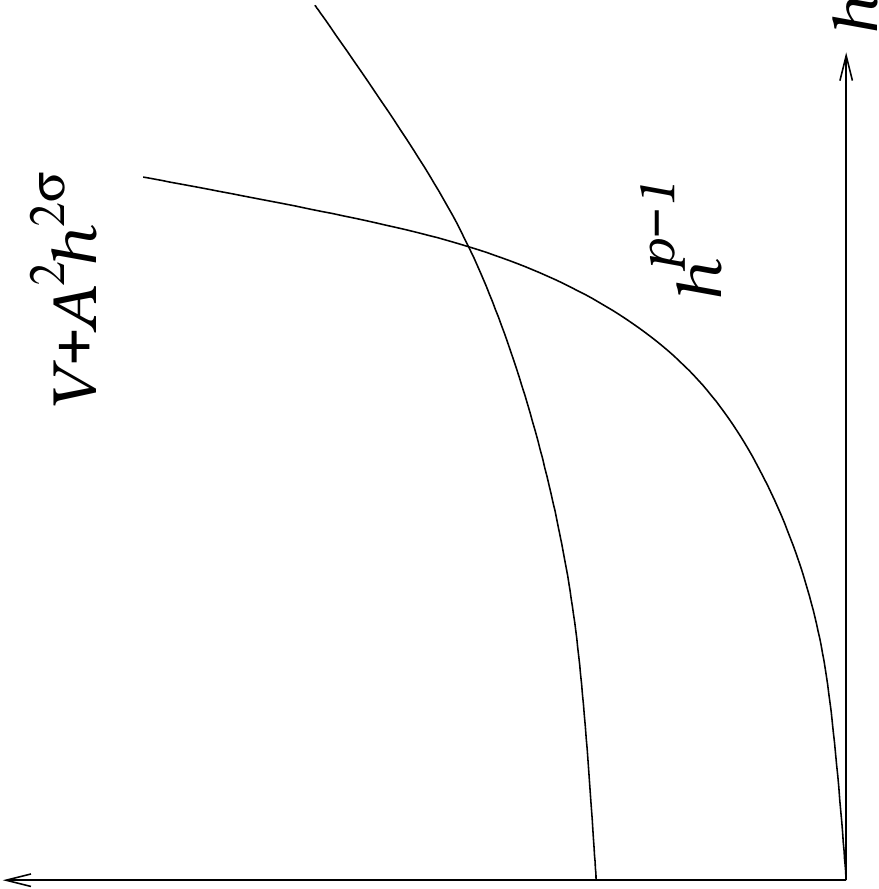} \qquad \qquad \qquad
\includegraphics[angle=270,width=5.5cm]{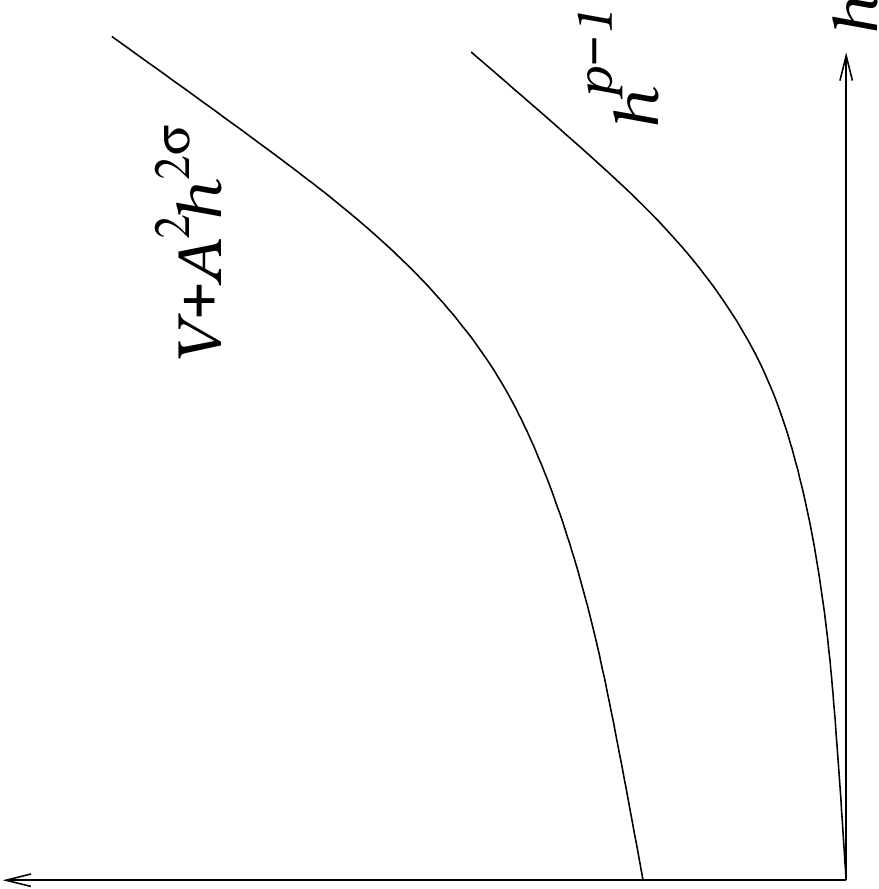}
\caption{the graphs of $V + \mathcal{A}^2 h^{2\s}$ and $h^{p-1}$
for $p < \frac{n+2}{n-2}$ and for $p = \frac{n+2}{n-2}$ with
$\mathcal{A} < 1$} \label{fig1}
\end{center}
\end{figure}

\begin{figure}[h] \begin{center}
 \includegraphics[angle=270,width=5.5cm]{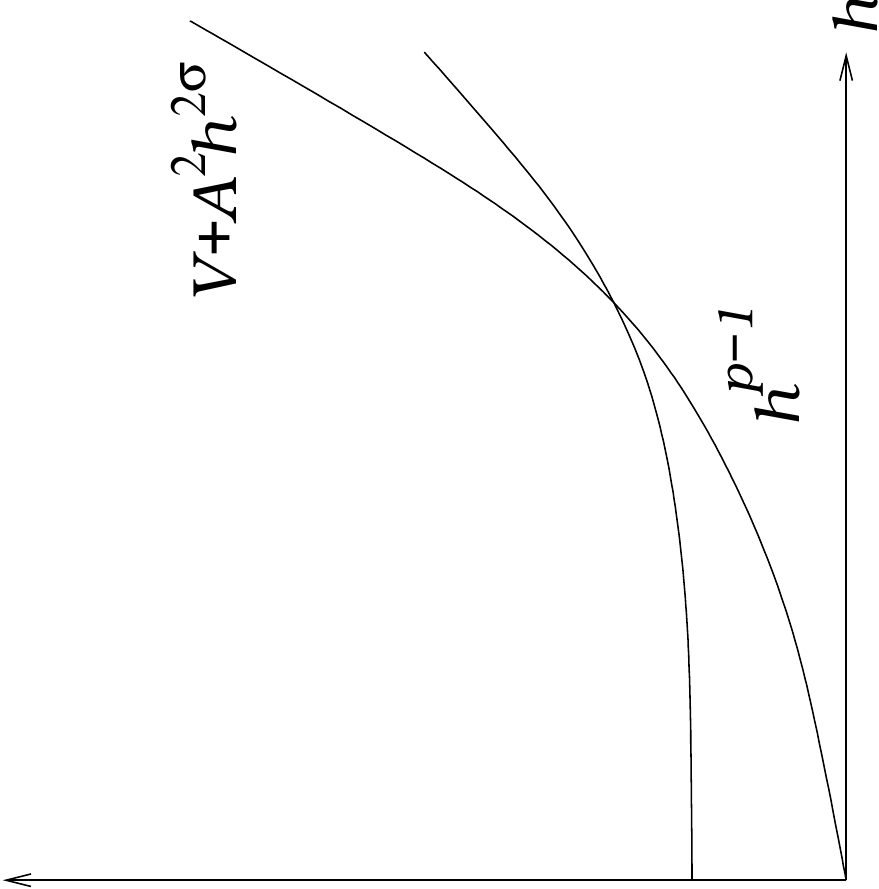} \qquad \qquad \qquad
\includegraphics[angle=270,width=5.5cm]{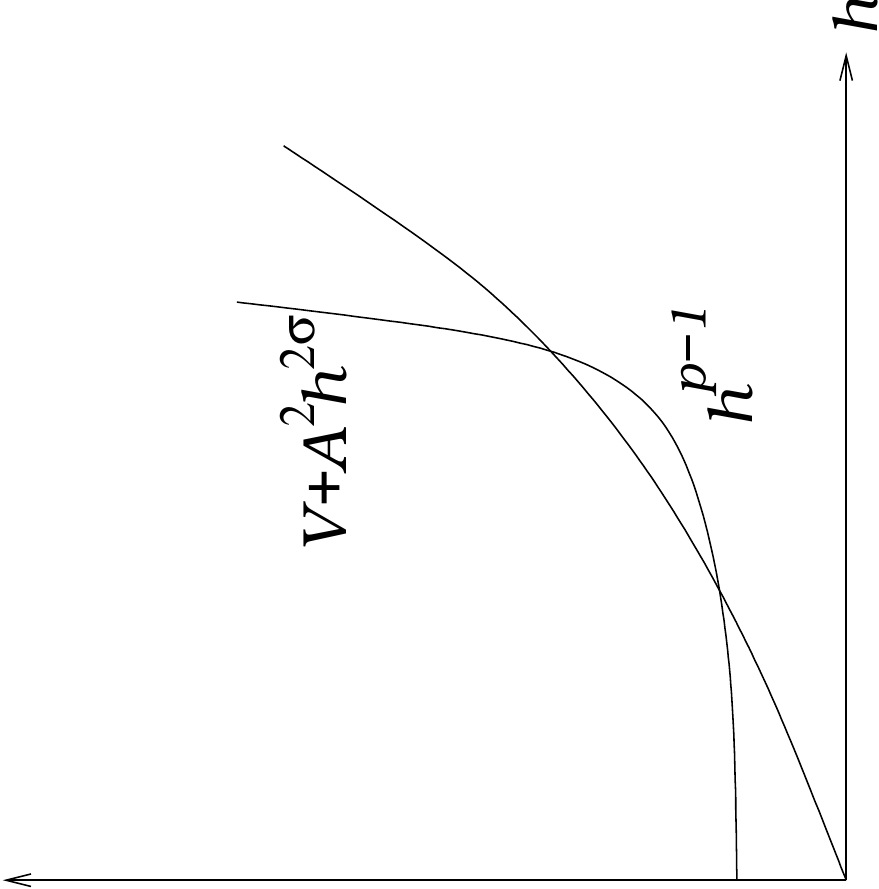}
\caption{the graphs of $V + \mathcal{A}^2 h^{2\s}$ and $h^{p-1}$
for $p = \frac{n+2}{n-2}$ with $\mathcal{A} \geq 1$ and for $p >
\frac{n+2}{n-2}$ with $\mathcal{A}$ small} \label{fig2}
\end{center}
\end{figure}
\end{rem}

\begin{rem}\label{r:heur}
In the above expansions, considering the terms of order $\e$, as
already noticed, we considered the metric $g$ to be flat near the
curve $\g$, and we tacitely assumed the potential $V$ to depend
only on the variable $\ov{s}$. Indeed, expanding the
Laplace-Beltrami operator and the potential $V$ taking the
variables $y$ into account, we obtain an extra term of order $\e$
which does not affect our computations since it turns out to be
odd in $y$, so it vanishes once multiplied by
$U\left(\frac{k(\ov{s}) y}{\e}\right)$ and integrated over
$\R^{n-1}$. For more details, we refer to Section \ref{s:as},
where precise estimates are worked out (in a system of coordinates
scaled in $\e$).
\end{rem}

\subsection{The Euler equation}\label{ss:euler}

Using the heuristic considerations of the previous subsection, we
now compute the energy of an approximate solution $\psi$
concentrated near a closed curve $\g$, and then find the $\g$'s for
which this energy is stationary. We let $\psi_{\g,\mathcal{A}}$
denote the function constructed in Subsection \ref{ss:cand}. In
order for the function $\psi_{\g,\mathcal{A}}$ to be globally well
defined, we need to impose one more restriction, namely that
$\psi_{\g,\mathcal{A}}$ is periodic in $\ov{s}$ with period $L$.
This is equivalent to require that $\int_0^L f'(\ov{s}) d\ov{s}$ is
an integer multiple of $2 \pi \e$, since we have the phase factor
$e^{- i \frac{f(\ov{s})}{\e}}$ in the expression of
$\psi_{\g,\mathcal{A}}$. From \eqref{eq:f'C}, then we find that also
$\int_0^L h(\ov{s})^\s d\ov{s}$ is an integer multiple of $2 \pi
\e$.

Multiplying \eqref{eq:pe} by $\psi_{\g,\mathcal{A}}$ and
integrating by parts, from the fact that $\psi_{\g,\mathcal{A}}$
is an approximate solution we find
\begin{eqnarray*}
  E_\e(\psi_{\g,\mathcal{A}}) & = & \frac 12 \int_M
  \left( \e^2 |\n_g \psi_{\g,\mathcal{A}}|^2 + V(x) |\psi_{\g,\mathcal{A}}|^2
  \right) dV_g - \frac{1}{p+1} \int_M |\psi_{\g,\mathcal{A}}|^{p+1} \\
  & \simeq & \left( \frac 12
  - \frac{1}{p+1} \right) \int_M |\psi_{\g,\mathcal{A}}|^{p+1} dV_g.
\end{eqnarray*}
Since $\psi_{\g,\mathcal{A}}$ is highly concentrated near $\g$,
using the coordinates $(\ov{s},y)$ introduced in \ref{ss:coord} we
have that
$$
  \int_M |\psi_{\g,\mathcal{A}}|^{p+1} dV_g \simeq \int_0^L d\ov{s} \int_{\R^{n-1}}
  h(\ov{s})^{p+1} \left| U \left( \frac{k(\ov{s}) y}{\e} \right) \right|^{p+1}
  dy.
$$
Using a change of variables, the last two formulas and
\eqref{eq:hk} we find that
\begin{equation}\label{eq:appen}
    E_\e(\psi_{\g,\mathcal{A}}) \simeq \ov{C}_0\e^{n-1} \int_\g h(\ov{s})^{\th} d\ov{s},
\end{equation}
where
$$
  \ov{C}_0 = \left( \frac 12 - \frac{1}{p+1}
    \right) \int_{\R^{n-1}} |U(y)|^{p+1} dy,
$$
and where we have set
\begin{equation}\label{eq:th}
    \th = p+1 - \frac 12 (p-1) (n-1)=p-\s-1.
\end{equation}

Consider now a one-parameter family of closed curves $\g_t : [0,L]
\to M$, where $t$ belongs to a neighborhood of $0$ in $\R$ and
where $\g_0 \equiv \g$. We compute next the approximate value of
the derivative in $t$ of the corresponding energy defined by
\eqref{eq:appen}.

As in Subsection \ref{ss:1st2ndvarlen} we let $\mathcal{V}_t$
denotes the vector field $\mathcal{V}_t(\ov{s}) = \frac{\pa
\g_t}{\pa t}(\ov{s})$ and we assume that $\mathcal{V} :=
\mathcal{V}_0$ is normal to $\g$. For any $t$ near zero, we let
$k_t(\ov{s}), h_t(\ov{s}), f_t(\ov{s})$ be defined by
\eqref{eq:hk} replacing $\g$ by $\g_t$ and $V(\ov{s})$ by
$V_t(\ov{s}) := V(\g_t(\ov{s}))$. Since we require periodicity of
each curve $\g_t$ in the variable $\ov{s}$, we also allow the
constant $\mathcal{A}$ given in \eqref{eq:f'C} to depend on $t$.
Denoting this by $\mathcal{A}_t$, by the above considerations we
choose $\mathcal{A}_t$ so that the following condition holds for
every value of $t$
\begin{equation}\label{eq:constr}
  \int_0^L \mathcal{A}_t h_t(\ov{s})^\s d\ov{s} = \int_0^L f_t'(\ov{s}) d\ov{s} = const.
\end{equation}
Below, we let $\mathcal{A}'_t = \frac{d}{dt} \mathcal{A}_t$ and we
will consider $h_t(\ov{s})$ as a function of $\mathcal{A}_t$ while
$V_t(\ov{s})$ as implicitly defined in \eqref{eq:hV}. From
\eqref{eq:varl} and  $\frac{\pa V_t(\ov{s})}{\pa t}|_{t=0} = \langle
\n^N V(\ov{s}), \mathcal{V}(\ov{s}) \rangle$, differentiating
\eqref{eq:constr} with respect to $t$ at $t = 0$ we get
\begin{equation*}
\int_0^L \mathcal{A} \s h^{\s-1} \frac{\partial h}{\partial V}
\langle \n^N V, \mathcal{V} \rangle d\ov{s} - \mathcal{A} \int_0^L
h^\s \langle \mathcal{V}, {\bf H} \rangle d\ov{s} + \mathcal{A}
\mathcal{A}' \s \int_0^L h^{\s-1} \frac{\partial h}{\partial
\mathcal{A}} d\ov{s} + \mathcal{A}' \int_0^L h^\s d\ov{s} = 0,
\end{equation*}
where we have set $\mathcal{A}' = \mathcal{A}'_0$ and where $\n^N
V$ stands for the component of $\n V$ normal to $\g$. From this
formula we obtain the following expression of $\mathcal{A}'$
\begin{equation}\label{eq:C'}
  \mathcal{A}' = -\mathcal{A} \frac{\int_0^L \left( \s h^{\s-1} \frac{\partial h}{\partial V}
  \langle \n^N V, \mathcal{V} \rangle - h^\s \langle \mathcal{V}, {\bf
H} \rangle  \right) d\ov{s}}{\int_0^L \left( \mathcal{A} \s
h^{\s-1} \frac{\partial h}{\partial \mathcal{A}} + h^\s \right)
d\ov{s}}.
\end{equation}
Similarly, computing the derivative of the (approximate) energy
with respect to $t$ we find
$$
  \frac{d E_\e(u_{\g_t,\mathcal{A}_t})}{dt}|_{t=0} = \int_0^L \left( \th h^{\th-1}
  \frac{\partial h}{\partial V} \langle \n^N V, \mathcal{V} \rangle -
  h^\th \langle \mathcal{V}, {\bf H} \rangle + \th \mathcal{A}' h^{\th-1}
  \frac{\partial h}{\partial \mathcal{A}} \right) d\ov{s}.
$$
Using \eqref{eq:C'} we deduce that the variation is given by
\begin{eqnarray*}
  \frac{d E_\e(u_{\g_t,\mathcal{A}_t})}{dt}|_{t=0} & = & \int_0^L \frac{\partial
  h}{\partial V} \langle \n^N V, \mathcal{V} \rangle
  \left[ \th h^{\th-1} - \frac{\mathcal{A} \s h^{\s-1} \int_0^L \th h^{\th-1}
  \frac{\partial h}{\partial \mathcal{A}} d\ov{s}}{\int_0^L \left( \mathcal{A} \s h^{\s-1}
  \frac{\partial h}{\partial \mathcal{A}} + h^\s \right) d\ov{s}} \right] d \ov{s} \\
  & - &
   \int_0^L \langle \mathcal{V}, {\bf H} \rangle \left[ h^{\th} - \frac{\mathcal{A} h^{\s}
  \int_0^L \th h^{\th-1} \frac{\partial h}{\partial \mathcal{A}} d\ov{s}}{\int_0^L
  \left( \mathcal{A} \s h^{\s-1} \frac{\partial h}{\partial \mathcal{A}} + h^\s \right) d\ov{s}}
  \right] d \ov{s}.
\end{eqnarray*}
Differentiating \eqref{eq:hV} with respect to $\mathcal{A}$ and
$V$ we get
\begin{equation}\label{eq:here}
    \frac{\partial h}{\partial \mathcal{A}} = \frac{2 \mathcal{A} h^{2\s}}{(p-1) h^{p-2} -
  2 \s \mathcal{A}^2 h^{2\s-1}} = 2 \mathcal{A} h^{2\s} \frac{\partial h}{\partial V},
\end{equation}
so it follows that
  \begin{equation}\label{eq:dhCV}
\mathcal{A} \s h^{\s-1} \frac{\partial h}{\partial \mathcal{A}} +
h^\s =
  \frac{(p-1) h^{p-1}}{(p-1) h^{\th} - 2 \s \mathcal{A}^2 h^{\s}}.
\end{equation}
Similarly, since $\th = p - \s - 1$, see \eqref{eq:s} and
\eqref{eq:th}, we deduce that
$$
  h^{\th-1} \frac{\partial h}{\partial \mathcal{A}} = \frac{2 \mathcal{A} h^{p+\s-2}}{(p-1)
  h^{p-2} -  2 \s \mathcal{A}^2 h^{2\s-1}}.
$$
Therefore we find also
\begin{equation}\label{eq:id}
    \frac{\th h^{\th-1} \frac{\partial h}{\partial \mathcal{A}}}{\mathcal{A} \s h^{\s-1}
    \frac{\partial h}{\partial \mathcal{A}} + h^\s} = \frac{2 \mathcal{A} \th}{p-1}.
\end{equation}
Hence from the last formulas the variation of the energy becomes
\begin{equation}\label{eq:qq}
   \frac{d E_\e(\psi_{\g_t,\mathcal{A}_t})}{dt}|_{t=0}=
    \int_0^L \langle \n^N V, \mathcal{V} \rangle \frac{\partial h}{\partial
 V} \left[ \th h^{\th - 1} - \frac{2 \mathcal{A}^2 \s \th}{p-1} h^{\s-1}
 \right] d\ov{s} - \int_0^L \langle \mathcal{V}, {\bf H} \rangle \left[ h^\th -
 \frac{2 \mathcal{A}^2 \th}{p-1} h^\s \right] d\ov{s}.
\end{equation}
Also, from the second equality in \eqref{eq:here}, dividing by
$h^\s$, multiplying by $\frac{\th}{p-1}$ and using the identity $p
- \s - 2 = \th - 1$ we obtain
$$
  h^{-\s} \frac{\th}{p-1} + 2 \s \mathcal{A}^2 h^{\s-1} \frac{\partial h}{\partial
  V} \frac{\th}{p-1} = \th h^{p-\s-2} \frac{\partial h}{\partial
  V} = \th h^{\th - 1} \frac{\partial h}{\partial V}.
$$
Using \eqref{eq:qq} and the last formula, we get the following
simplified expression
$$
  \frac{d E_\e(\psi_{\g_t,\mathcal{A}_t})}{dt}|_{t=0}=
  \int_0^L \frac{\th}{p-1} h^{-\s} \left[ \langle \n^N V, \mathcal{V}
  \rangle - \langle \mathcal{V}, {\bf H} \rangle \left( \frac{p-1}{\th}
h^{p-1} - 2 \mathcal{A}^2
  h^{2\s} \right) \right] d\ov{s}.
$$
Therefore the stationarity condition for the energy (under the
constraint \eqref{eq:constr}) becomes $\langle \n^N V, \mathcal{V}
  \rangle = \langle \mathcal{V}, {\bf H} \rangle \left(
  \frac{p-1}{\th} h^{p-1} - 2 \mathcal{A}^2 h^{2\s} \right)$ for
  every normal section $\mathcal{V}$, namely
\begin{equation}\label{eq:euler}
    \n^N V = {\bf H} \left( \frac{p-1}{\th} h^{p-1} - 2
    \mathcal{A}^2 h^{2\s} \right).
\end{equation}
We will see that this formula will be crucial later on to find
approximate solutions.

\begin{remark}\label{r:monC}
By \eqref{eq:dhCV}, we have that
$$
  \frac{\pa}{\pa \mathcal{A}} \left( \mathcal{A} h^\s \right) =
  \frac{(p-1) h^{p-1}}{(p-1) h^{\th} - 2 \s \mathcal{A}^2 h^{\s}}.
$$
If $\mathcal{A}$ is sufficiently small (depending on $V$ and $p$),
then we have $\frac{\pa}{\pa \mathcal{A}} \left( \mathcal{A} h^\s
\right) > 0$. This will be used in the second part \cite{mmm2} when,
for a fixed $\e$, we will adjust the value of the constant
$\mathcal{A}$ for obtaining periodicity of the function $f$.
\end{remark}

\subsection{Second variation and non-degeneracy
condition}\label{ss:2ndvarEuler}

We evaluate next the second variation of the Euler functional. As
in Subsection \ref{ss:1st2ndvarlen} we consider a two-parameter
family of closed curves $\g_{t_1,t_2}$, where $t_1, t_2$ are two
real numbers belonging to a small neighborhood of zero in $\R$,
and where $\g_{0,0} = \g$. As before, we require the constraint
\eqref{eq:constr} along the whole two-dimensional family of
curves, and we assume the functions $f, h, k$ and the constant
$\mathcal{A}$ to depend on $t_1$ and $t_2$, and we will use the
notation $\mathcal{A}_{t_1,t_2}$, etc.. Keeping this in mind, we
define the two vector fields $\mathcal{V}_{t_1,t_2} = \frac{\pa
\g_{t_1,t_2}}{\pa t_1}$, $\mathcal{W}_{t_1,t_2} = \frac{\pa
\g_{t_1,t_2}}{\pa t_2}$, and we can assume that $\mathcal{V} :=
\mathcal{V}_{0,0}, \mathcal{W} := \mathcal{W}_{0,0}$ are normal to
the initial curve $\g$. With some computations, which are worked
out in Section \ref{s:apx}, one finds that, at $(t_1,t_2)=(0,0)$
\begin{eqnarray}\label{eq:2ndvarfin} \nonumber
  \frac{\pa^2 E_\e(u_{\psi_{t_1,t_2},\mathcal{A}_{t_1,t_2}})}{\pa t_1 \pa t_2}
  & = & \int_0^L \left[ h^\th - \frac{2 \mathcal{A}^2 \th}{p-1} h^\s \right] \bigg[ \sum_j
  \dot{\mathcal{V}}^j \dot{\mathcal{W}}^j - \sum_{j,m} R_{1j1m}
  \mathcal{V}^j \mathcal{W}^m  \bigg] d\ov{s}
   \\
  & + & \frac{\th}{p-1} \int_0^L \bigg\{ ((\n^N)^2 V) [\mathcal{V},\mathcal{W}]
  - \langle \n^N V, \mathcal{V} \rangle
  \langle {\bf H}, \mathcal{W} \rangle - \langle \n^N V, \mathcal{W} \rangle
  \langle {\bf H}, \mathcal{V} \rangle \bigg\} h^{-\s} d\ov{s} \\
    & - & \frac{\s \th}{p-1} \int_0^L h^{-\s-1}\frac{\pa h}{\pa V} \langle \n^N V,
    \mathcal{V} \rangle \langle \n^N V, \mathcal{W} \rangle d\ov{s}
    \nonumber
    \\ & + & \mathcal{A}'_1 \mathcal{A}'_2 \frac{2 \th}{p-1}
    \int_0^L \left( \mathcal{A} \s h^{\s-1}
  \frac{\partial h}{\partial \mathcal{A}} + h^\s \right) d\ov{s}. \nonumber
    \nonumber
\end{eqnarray}
Here, $\mathcal{A}'_l$ stands for $\frac{\pa
\mathcal{A}_{t_1,t_2}}{\pa t_l}|_{(t_1,t_2)=(0,0)}$,
$\dot{\mathcal{V}}^j = \frac{d \mathcal{V}^j}{d \ov{s}}$, and
$\dot{\mathcal{W}}^j = \frac{d \mathcal{W}^j}{d \ov{s}}$, where the
$\mathcal{V}^j$'s, $\mathcal{W}^j$'s are the components of
$\mathcal{V}$ and $\mathcal{W}$ with respect to the basis $(E_j)_j$
introduced in Subsection \ref{ss:coord}.

Integrating by parts and using \eqref{eq:C'}, from the last formula
one derives that the non-degeneracy condition is equivalent to the
invertibility of the linear operator $\mathfrak{J} : \chi(N \g) \to
\chi(N \g)$ (from the family of smooth sections of the normal bundle
to $\g$ into itself) whose components are defined by
\begin{eqnarray}\label{eq:2ndvarfin2} \nonumber
  (\mathfrak{J} \mathcal{V})^m & = & - \left( h^\th - \frac{2 \mathcal{A}^2 \th}{p-1} h^\s \right)
  \left[\ddot{\mathcal{V}}^m + \sum_j R_{1j1m} \mathcal{V}^j \right]-
  \th \left( h^{\th-1} - \frac{2 \mathcal{A}^2 \s}{p-1} h^{\s-1} \right) h'
\dot{\mathcal{V}}^m
  \\ & + & \frac{\th}{p-1} h^{-\s} \left\{ ((\n^N)^2 V)(\mathcal{V},E_m)
  - H^m \langle \n^N V, \mathcal{V} \rangle - \langle {\bf H},
  \mathcal{V} \rangle \langle \n^N V, E_m \rangle \right\} \nonumber \\
  &-& \frac{\s
  \th}{p-1} h^{-(\s+1)} \frac{\partial h}{\partial V} \langle \n^N V, \mathcal{V}
  \rangle \langle \n^N V, E_m \rangle
   - \frac{2\th}{p-1} \mathcal{A} \mathcal{A}'_1
  \left( \s h^{\s-1} \frac{\partial h}{\partial V} \langle \n^N V, E_m \rangle
  - h^\s H^m \right), \nonumber
\end{eqnarray}
where $h' = \frac{d h(\ov{s})}{d \ov{s}}$. Recalling  formula
\eqref{eq:ctc}, using \eqref{eq:euler} and some other elementary
computations we obtain
\begin{eqnarray*}
  (\mathfrak{J} \mathcal{V})^m = - \left( h^\th - \frac{2
  \mathcal{A}^2 \th}{p-1} h^\s \right) \ddot{\mathcal{V}}^m - \th
  \left( h^{\th-1} - \frac{2 \mathcal{A}^2 \s}{p-1} h^{\s-1} \right) h'
  \dot{\mathcal{V}}^m
  + \frac{\th}{p-1} h^{-\s} ((\n^N)^2 V)(\mathcal{V},E_m) \\ \nonumber
  + \frac 12 \left( h^\th - \frac{2 \mathcal{A}^2 \th}{p-1} h^\s \right) \left(
  \sum_j \left(  \partial^2_{jm} g_{11} \right) {\mathcal{V}}^j \right)
  - \frac{2\th}{p-1} \mathcal{A} \mathcal{A}'_1 \left( \s h^{\s-1}
  \frac{\partial h}{\partial V} \langle \n^N V, E_m \rangle -
    h^\s H^m \right) \\
  +  H^m \langle {\bf H}, \mathcal{V} \rangle
  \left( \frac{2 \mathcal{A}^2 \th}{p-1} h^\s - h^\th \right) \left[
  \frac{(p-1) \left( 3 + \frac{\s}{\th} \right) h^{\th} - 8 \s \mathcal{A}^2
  h^{\s}}{(p-1) h^{\th} - 2 \mathcal{A}^2 \s h^{\s}} \right].
\end{eqnarray*}
For future convenience, we expand the last product explicitly,
finding
$$
  H^m \langle {\bf H}, \mathcal{V} \rangle
  \left[ \frac{- (p-1) \left( 3 + \frac{\s}{\th} \right)
  h^{2\th} - \frac{16 \s \th \mathcal{A}^4}{p-1} h^{2\s} + 2 \mathcal{A}^2 (5\s +
  3\th) h^{\th+\s}}{(p-1) h^{\th} - 2 \mathcal{A}^2 \s h^{\s}} \right].
$$
We also notice that
$$
  \s h^{\s-1} \langle \n^N V, E_m \rangle \frac{\pa h}{\pa V}
  - h^\s H^m = \frac{(p-1) (\th-\s)
  h^{p-1}}{\th \left[ (p-1) h^{\th} - 2 \s \mathcal{A}^2 h^{\s}
  \right]} H^m.
$$
In conclusion the non-degeneracy condition is equivalent to the
invertibility of the operator $\mathfrak{J} : \chi(N \g) \to \chi(N
\g)$ given in components by
\begin{eqnarray}\label{eq:2ndvarfin4} \nonumber
  (\mathfrak{J} \mathcal{V})^
  m & = & - \left( h^\th - \frac{2 \mathcal{A}^2
  \th}{p-1} h^\s \right) \ddot{\mathcal{V}}^m - \th
  \left( h^{\th-1} - \frac{2 \mathcal{A}^2 \s}{p-1} h^{\s-1}
  \right) h' \dot{\mathcal{V}}^m
  + \frac{\th}{p-1} h^{-\s} ((\n^N)^2 V)[\mathcal{V},E_m]
  \\  & + & \frac 12 \left( h^\th - \frac{2
  \mathcal{A}^2 \th}{p-1} h^\s \right) \left(
  \sum_j \left(  \partial^2_{jm} g_{11} \right) {\mathcal{V}}^j \right)
  - 2 \mathcal{A} \mathcal{A}'_1 \frac{(\th-\s)
  h^{p-1}}{\left[ (p-1) h^{\th} - 2 \s \mathcal{A}^2 h^{\s}
  \right]} H^m \\
  & + & H^m \langle {\bf H}, \mathcal{V} \rangle
  \left[ \frac{- (p-1) \left( 3 + \frac{\s}{\th} \right)
  h^{2\th} - \frac{16 \s \th \mathcal{A}^4}{p-1} h^{2\s} + 2 \mathcal{A}^2 (5\s +
  3\th) h^{\th+\s}}{(p-1) h^{\th} - 2 \mathcal{A}^2 \s h^{\s}}
  \right]. \nonumber
\end{eqnarray}

\no We summarize the results of this section in the following
Proposition.

\begin{pro}\label{p:eulnd} Consider the functional on curves $\int_\g h^\th(\ov{s}) \;
d \ov{s}$ restricted to the set $\G$ in \eqref{eq:Gamma}. Then the
stationarity condition is \eqref{eq:euler} and the non-degeneracy of
a critical point is equivalent to the invertibility of the operator
$\mathfrak{J}$ in \eqref{eq:2ndvarfin4}.
\end{pro}

\section{Approximate solutions}\label{s:as}

Using a change of variables, equation \eqref{eq:pe} is equivalent
to the following
\begin{equation}\label{eq:new}
  - \D_{g_\e} \psi + V(\e x) \psi = |\psi|^{p-1} \psi \qquad \hbox{ in
  } M_\e,
\end{equation}
where $M_\e$ denotes the manifold $M$ endowed with the scaled
metric $g_\e = \frac{1}{\e^2} g$. With an abuse of notation we
will often denote it through the scaling $M_\e = \frac 1 \e M$,
and if $x \in M_\e$ we write $\e x$ to indicate the corresponding
point on $M$.

In this section we find a family of approximate solutions to the
scaled equation \eqref{eq:new}. We consider a simple closed curve
$\g$ which is stationary within the class $\G$, namely satisfying
\eqref{eq:euler}. First, we introduce some convenient coordinates
near the scaled curve $\g_\e = \frac 1\e \g$, expanding the
Laplace-Beltrami operator with respect to the scaled metric in
powers of $\e$. Then, using these expansions, we construct the
approximate solutions solving formally \eqref{eq:new} up to order
$\e$. Since in the second part \cite{mmm2} we will need to work
out rigorous estimates, in order not to repeat later the
expansions we will treat some terms carefully and not only at a
formal level.

\subsection{Choice of coordinates in $M_\e$ and expansion of the metric coefficients}

Using the coordinates $(\ov{s},y)$ of Section \ref{s:not} defined
near $\g$, for some smooth normal section $\Phi(\ov{s})$ in $N \g$,
we define the following new coordinates $(s,z)$ (here and below we
use the notation $\ov{s} = \e s$) near $\frac 1\e \g$
\begin{equation}\label{eq:defzz}
    z = y - \Phi(\e s); \qquad \qquad z \in \R^{n-1}.
\end{equation}
 In this choice
we are motivated by the fact that in general we allow the
approximate solutions to be {\em tilted} normally to $\g_\e$, where
the tilting $\Phi$ depends (slowly) on the variable $s$: this allows
some extra flexibility in the construction, as in \cite{dkw},
\cite{mm} and \cite{mal} . As we will see, the choice of $\Phi$ is
irrelevant for solving \eqref{eq:new} up to order $\e$; on the other
hand,  the non-degeneracy assumption  will be necessary to guarantee
solvability of the equation up to higher orders.

We denote by $\tilde{g}_{AB}$ the metric coefficients in the new
coordinates $(s,z)$. Since $y = z + \Phi(\e s)$, it follows
$$
  \tilde{g}_{CD} = \sum_{AB} g_{AB} \left( \frac{\partial y_A}{\partial
  z_C} \right) \left( \frac{\partial y_B}{\partial z_D} \right).
$$
Explicitly, we then find
\begin{equation*}
    \tilde{g}_{11} = g_{11}|_{z+\Phi} + 2 \e \sum_j \Phi'_j
    g_{1j}|_{z+\Phi} + \e^2 \sum_{j,m} \Phi'_j(\e s) \Phi'_m(\e s) g_{jm}|_{z+\Phi};
\end{equation*}
\begin{equation*}
    \tilde{g}_{1j} = g_{1j}|_{z+\Phi} + \e \sum_m \Phi'_m(\e s)
    g_{jm}|_{z+\Phi}; \qquad \qquad \qquad \tilde{g}_{jm} = g_{jm}|_{z+\Phi}.
\end{equation*}
At this point, it is convenient to introduce some notation. For a
positive  integer $q$, we denote by $R_q(z)$, $R_q(z,\Phi)$ and
$R_q(z,\Phi,\Phi')$ error terms which satisfies respectively the
following bounds, for some positive constants $C$ and $d$
$$|R_q(z)|\leq C \e^q (1+|z|^d),$$
\begin{equation*}
    \left\{
      \begin{array}{ll}
        |R_q(z,\Phi)| \leq C \e^q (1+|z|^d); &  \\
        |R_q(z,\Phi) - R_q(z,\tPhi)| \leq C
    \e^q (1+|z|^d) [\,|\Phi-\tPhi|], &
      \end{array}
    \right.
\end{equation*}
and
\begin{equation*}
    \left\{
      \begin{array}{ll}
        |R_q(z,\Phi,\Phi')| \leq C \e^q (1+|z|^d); &  \\
        |R_q(z,\Phi,\Phi') - R_q(z,\tPhi,\tPhi')| \leq C
    \e^q (1+|z|^d) [\,|\Phi-\tPhi| + |\Phi'-\tPhi'|]. &
      \end{array}
    \right.
\end{equation*}
We also introduce error terms involving also second derivatives of
$\Phi$, $R_q(z,\Phi,\Phi',\Phi'')$ which satisfy
\begin{equation*}
    |R_q(z,\Phi,\Phi',\Phi'')| \leq C \e^q (1+|z|^d) + C \e^{q+1}
    (1+|z|^d) |\Phi''|;
\end{equation*}
\begin{eqnarray*}
 \nonumber |R_q(z,\Phi,\Phi',\Phi'') -
  R_q(z,\tPhi,\tPhi',\tPhi'')| & \leq & C
    \e^q (1+|z|^d) [|\Phi-\tPhi| + |\Phi'-\tPhi'|] \left(1 + \e (|\Phi''| +
    |\tilde{\Phi}''|) \right) \\ & + & C \e^{q+1} (1+|z|^d) |\Phi'' - \tPhi''|.
\end{eqnarray*}
Using the expansion of the metric coefficients $g_{AB}$ in Lemma
\ref{l:me} and this notation, we then obtain
\begin{eqnarray}\label{eq:tg112} \nonumber
    \tilde{g}_{11} & = & 1 - 2 \e \sum_{m=2}^n H^m (z_m +
    \Phi_m) + \frac 12 \e^2 \sum_{m,l=2}^n \partial^2_{ml} g_{11}|_{\g}
    (z_m + \Phi_m) (z_l + \Phi_l) \\ & + &
\e^2 |\Phi'|^2 +  R_3(z,\Phi,\Phi');
\end{eqnarray}
\begin{eqnarray}\label{eq:tg1j2} \nonumber
    \tilde{g}_{1j} & = &
    \e \Phi'_j + \frac 12 \e^2 \sum_{m,l=2}^n \partial^2_{ml}
    g_{1j}|_{\g} (z_m + \Phi_m) (z_l + \Phi_l) + R_3(z,\Phi,\Phi');
\end{eqnarray}
\begin{eqnarray}\label{eq:tgjk2} \nonumber
    \tilde{g}_{kj} & = & \d_{kj} + \frac 12 \e^2 \sum_{m,l=2}^n
    \partial^2_{ml} g_{kj}|_{\g} (z_m + \Phi_m) (z_l + \Phi_l)+ R_3(z,\Phi,\Phi').
\end{eqnarray}
Next we compute the inverse metric coefficients. Recall that, given
a formal expansion of a matrix as $M = 1 + \e A + \e^2 B$, we have
$$
  M^{-1} = 1 - \e A + \e^2 A^2 - \e^2 B.
$$
In our specific case the matrix $A$ is the following
\begin{equation}\label{eq:mA}
    A = \left(
\begin{array}{cc}
  - 2 \sum\limits_{m=2}^n H^m (z_m + \Phi_m) &
  \Phi'_j
  \\
  \Phi'_j
  & 0 \\
\end{array}
\right),
\end{equation}
and the elements of its square are given by
$$
  (A^2)_{11} =
  4 \left( \sum_{m=2}^n H^m (z_m + \Phi_m) \right)^2 +
  \sum_j \left( \Phi'_j
  \right)^2;
$$
$$
  (A^2)_{1j} = -2  \left( \sum_{m=2}^n H^m (z_m + \Phi_m)
  \right) \left(
  \Phi'_j
  \right);\qquad \hbox{and} \qquad
  (A^2)_{lj} = \left( \Phi'_l
  \right) \left( \Phi'_j
  \right).
$$
Therefore, using the above formula, the inverse coefficients are
\begin{eqnarray}\label{eq:tg112-1} \nonumber
    \tilde{g}^{11} & = & 1 + 2 \e \sum_{m=2}^n H^m (z_m +
    \Phi_m) - \frac 12 \e^2 \sum_{m,l=2}^n \partial^2_{ml} g_{11}|_{\g}
    (z_m + \Phi_m) (z_l + \Phi_l) \\ & + &
    4 \e^2 \left( \sum_{m=2}^n
    H^m (z_m + \Phi_m) \right)^2 +R_3(z,\Phi,\Phi'),
   \nonumber
\end{eqnarray}
We also get
\begin{eqnarray*}
    \tilde{g}^{1j} & = &
    - \e \Phi'_j - \frac 12 \e^2 \sum_{m,l=2}^n \partial^2_{ml}
    g_{1j}|_{\g} (z_m + \Phi_m) (z_l + \Phi_l) - 2 \e^2 \left(
    \sum_{m=2}^n H^m (z_m + \Phi_m) \right)
    \Phi'_j
\\
    &+& R_3(z,\Phi,\Phi').
\end{eqnarray*}
Moreover
\begin{eqnarray*}
   \pa_j(\tilde{g}^{1j}) = - \e^2 \sum_{l=2}^n\pa^2_{lj} g_{1j}|_\g (z_l
+ \Phi_l) - 2 \e^2 H^j \Phi'_j + R_3(z,\Phi,\Phi').
\end{eqnarray*}
Similarly, with some simple calculations one also finds
\begin{eqnarray}\label{eq:tg112-12} \nonumber
    \pa_1 (\tilde{g}^{11}) & = & \nonumber 2 \e^2 \sum_{m=2}^n (H^m)' (z_m +
   \Phi_m) + 2 \e^2 \sum_{m=2}^n H^m \Phi'_m + R_3(z,\Phi,\Phi',\Phi'').
\end{eqnarray}
Differentiating now $\tilde{g}^{1j}$ with respect to the first
variable we obtain
$$
  \pa_1(\tilde{g}^{1j}) = - \e^2 \Phi''_j - 2 \e^3 \left(
    \sum_{m=2}^n H^m (z_m + \Phi_m) \right) \Phi''_j +
 R_3(z,\Phi,\Phi',\Phi'').
$$
Analogously, we get
\begin{equation}\label{eq:tgjk2-1} \nonumber
    \tilde{g}^{kj} = \d_{kj} - \frac 12 \e^2 \sum_{m,l=2}^n
    \partial^2_{ml} g_{kj}|_{\g} (z_m + \Phi_m) (z_l + \Phi_l) +
    \e^2 \Phi'_k \Phi'_j + R_3(z,\Phi,\Phi');
\end{equation}
\begin{equation}\label{eq:tgjk2-122} \nonumber
    \pa_k(\tilde{g}^{kj}) =  - \e^2 \sum_{l=2}^n
    \partial^2_{kl} g_{kj}|_{\g} (z_l + \Phi_l) +
R_3(z,\Phi,\Phi').
\end{equation}

Finally, using the formal expansion $\tilde{g}_{CD} = \d_{CD} + \e
A_{CD} + \e^2 B_{CD} + o(\e^2)$, analyzing carefully the error
terms we obtain
$$
  \sqrt{\det \tilde{g}} = 1 + \frac 12 \e tr(A) + \e^2 \left( \frac 18
  (tr(A))^2 - \frac 14 tr(A^2) \right) + \frac 12 \e^2 tr(B)+O(\e^3).
$$
{From} the above expressions in \eqref{eq:tg112}, \eqref{eq:tg1j2}
we deduce that
\begin{eqnarray*}
  \sqrt{\det \tilde{g}} & = & 1 - \e \sum_m H^m (z_m +
  \Phi_m) \\ & + & \e^2 \left[ \frac 14 \sum_{m,l} \partial^2_{ml}
  g_{11} (z_m + \Phi_m) (z_l + \Phi_l) - \frac 12 \left( \sum_{m=2}^n
  H^m (z_m + \Phi_m) \right)^2 \right]\\
  & + & R_3(z,\Phi,\Phi'); \nonumber
\end{eqnarray*}
\begin{eqnarray*}
  \partial_m \sqrt{\det \tilde{g}} & = & - \e H^m +
  \e^2 \left[ \frac 12 \sum_l \partial^2_{ml} g_{11} (z_l + \Phi_l)
  - H^m \left( \sum_l H^l (z_l + \Phi_l)
  \right) \right]\\
  & + &R_3(z,\Phi,\Phi'),
\end{eqnarray*}
 moreover
\begin{eqnarray}\label{eq:d1sqrtdg} \nonumber
  \partial_1 \sqrt{\det \tilde{g}} & = & - \e^2 \sum_m (H^m)' (z_m +
\Phi_m) - \e^2 \sum_m H^m \Phi'_m + R_3(z,\Phi,\Phi',\Phi'').
\end{eqnarray}

\

The Laplacian of a smooth function $u$ in coordinates $(s,z)$ has
the following expression
$$
  - \D_{\tilde{g}} u = - \sum_{A,B} \tilde{g}^{AB} \pa^2_{AB}u - \sum_{A,B}\partial_B (\tilde{g}^{AB})
  \partial_A u - \frac{1}{\sqrt{\det \tilde{g}}} \sum_{A,B}\tilde{g}^{AB} \left( \partial_B
  \sqrt{\det \tilde{g}} \right) \partial_A u.
$$
We are going to expand next each of these terms. First, we consider
the determinant of  $\tilde{g}$. Recall that for a matrix of the
form $1 + \e A + \e^2 B$ the square root of the determinant admits
the formal expansion
\begin{equation}\label{dvg}
\sqrt{\hbox{det}g} = 1+ \frac{\varepsilon}{2} tr A +
\varepsilon^2\left(\frac{1}{8}(tr\; A)^2-\frac{1}{4}tr(A^2) +
\frac 12 tr B \right)+o(\varepsilon^2).
\end{equation}

\begin{lem}\label{l:explapltildeg}
Let $u$ be a smooth function. Then in the above coordinates
$(s,z)$ we have that
\begin{eqnarray*}
  & \D_{\tilde{g}} u  & = \pa^2_{ss}u+\D_{z}u
-\e \sum_jH^j\pa_ju -2\e\sum_j\Phi'_j\pa^2_{sj}u+2\e\langle {\bf
H},z+\Phi\rangle
\pa^2_{ss}u\\
&-&\e^2\langle {\bf H},z+\Phi\rangle\sum_{m,j}H^j\pa_j
u-\frac12\e^2\pa^2_{ml}g_{11}(z_m+\Phi_m)(z_l+\Phi_l)\pa^2_{ss}u+4\e^2\langle
{\bf H},z+\Phi\rangle^2\pa^2_{ss}u\\
&-&\e^2\pa^2_{ml}g_{1j}(z_m+\Phi_m)(z_l+\Phi_l)\pa^2_{sj}u-
4\e^2\langle {\bf
H},z+\Phi\rangle\sum_j\Phi'_j\pa^2_{sj}u+\e^2\sum_{t,j}\Phi'_t\Phi'_j\pa^2_{tj}u\\
&-&\frac12\e^2\sum_{m,l}\pa^2_{ml}g_{tj}(z_m+\Phi_m)(z_l+\Phi_l)\pa^2_{tj}u+
\e^2 \langle {\bf H}',z+\Phi\rangle \pa_s u - \e^2 \sum_{l,j}
\pa^2_{lj} g_{1j} (z_l + \Phi_l) \pa_s u \\
&-&\e^2\sum_j\Phi_j''\pa_j
u-\e^2\sum_{t,j,l}\pa^2_{tl}g_{tj}(z_l+\Phi_l)\pa_ju-2\e^3\langle
{\bf H},z+\Phi\rangle\sum_{j}\Phi_j'' \pa_j u
\\
&+&R_3(z,\Phi,\Phi')\pa^2_{ss}u+R_3(z,\Phi,\Phi')\pa^2_{sj}u+R_3(z,\Phi,\Phi')\pa^2_{lj}u+R_3(z,\Phi,\Phi',\Phi'')
\left( \pa_su+\pa_j u \right).
\end{eqnarray*}
Moreover, given two smooth normal sections $\Phi$ and
$\tilde{\Phi}$ and defining the corresponding  coordinates
$$(s,y-\Phi(\e s)) \qquad
\hbox{ and} \qquad (s,y-\tilde{\Phi}(\e s))
$$
and set $u_\Phi(s,y):=u(s,y-\Phi(\e s))$,
$u_{\tilde{\Phi}}(s,y):=u(s,y-\tilde{\Phi}(\e s))$. We then have
\begin{eqnarray*}
  \D_{\tilde{g}} u_{\Phi}-\D_{\tilde{g}} u_{\tilde \Phi}
  & = &-2\e\sum_j(\Phi'_j-\tilde{\Phi}'_j)\pa^2_{sj}u
+2\e\langle {\bf
H},\Phi-\tPhi\rangle\pa^2_{ss}u+\e^2\sum_{t,j}(\Phi'_t\Phi'_j-\tPhi'_t\tPhi'_j)\pa^2_{tj}u\\
&-&\frac12\e^2\sum_{m,l}\pa^2_{ml}g_{tj}\bigg[2z_m(\Phi_l-\tPhi_l)
+\Phi_l(\Phi_m-\tPhi_m)+\tPhi_l(\Phi_m-\tPhi_m)\bigg] \pa^2_{tj}u
\\
&-&\e^2\sum_{m,l}\pa^2_{ml}g_{1j}\bigg[2z_m(\Phi_l-\tPhi_l)
+\Phi_l(\Phi_m-\tPhi_m)+\tPhi_l(\Phi_m-\tPhi_m)\bigg] \pa^2_{sj}u
\\
&-&\frac12\e^2\sum_{m,l}\pa^2_{ml}g_{11}\bigg[2z_m(\Phi_l-\tPhi_l)
+\Phi_l(\Phi_m-\tPhi_m)+\tPhi_l(\Phi_m-\tPhi_m)\bigg] \pa^2_{ss}u\\
&-&2\e^2\sum_l
H^j\bigg[z_l(\Phi'_l-\tPhi'_l)+\Phi_l(\Phi'_l-\tPhi'_l)+\tPhi'_l(\Phi_l-\tPhi_l)
\bigg]\pa^2_{sj}u
\\
&+& 4\e^2\sum_{m,l}H^mH^l\bigg[2z_m(\Phi_l-\tPhi_l)
+\Phi_l(\Phi_m-\tPhi_m)+\tPhi_l(\Phi_m-\tPhi_m)\bigg]
\pa^2_{ss}u\\
&-&\e^2\sum_j(\Phi_j''-\tPhi''_j)\pa_ju-\e^2\sum_{t,j,l}\pa^2_{tl}
g_{tj}(\Phi_l-\tPhi_l)\pa_ju -\e^2\langle{\bf
H},\Phi-\tilde{\Phi}\rangle\sum_{j}H^j\pa_j
u\\
&+& \e^2 \langle{\bf H}',\Phi-\tilde{\Phi}\rangle \pa_s u - \e^2
\sum_{l,j}
\pa^2_{lj} g_{1j} (\Phi_l - \tilde{\Phi}_l) \pa_s u \\
&-&2\e^3\sum_{mj}H^m\bigg[(z_m+\Phi_m)(\Phi_j''-\tPhi_j'')+\tPhi_j''(\Phi_m-\tPhi_m)
\bigg]\pa_j u \\
&+&O(1+|z|^d)\bigg[\e^4(|\Phi-\tPhi|+|\Phi'-\tPhi'|)|\pa^2_{ss}u|+
\e^3(|\Phi-\tPhi|+|\Phi'-\tPhi'|)(|\pa^2_{sj}u|+|\pa^2_{lj}u)|\bigg]
\\
&+&O(1+|z|^d)\bigg[\e^3(|\Phi-\tPhi|+|\Phi'-\tPhi'|)+
\e^4(|\Phi''||\Phi-\tPhi|+|\Phi''-\tPhi''|)\bigg]
(|\pa_{s}u|+|\pa_{j}u|).
\end{eqnarray*}
\end{lem}
\begin{pf}
The proof is based on the Taylor expansion of the metric
coefficients given above. We recall that the Laplace-Beltrami
operator is given by
\[ \Delta_{
\tilde{g}}=\sum_{A,B}\frac{1}{\sqrt{\det
\tilde{g}}}\,\partial_A(\,\sqrt{\det
\tilde{g}}\,(g_\e)^{AB}\,\partial_B\,)\,,\] where indices $A$ and
$B$ run between $1$ and $n$. We can also write~
\[\Delta_{\tilde{g}}=\sum_{A,B}\bigg({\tilde{g}}^{AB}\,\partial^2_{AB}+\left(\partial_A\,{
\tilde{g}}^{AB}\right)\,\partial_B+\frac{1}{\sqrt{\det \tilde{g}}}
\tilde{g}^{AB} \left( \partial_B
  \sqrt{\det \tilde{g}} \right) \partial_A \bigg).
\]
Using the expansion of the metric coefficients determined above and
\eqref{dvg}, one can easily prove that
\begin{eqnarray*}
\sum_{AB} \tilde{g}^{AB}
\pa^2_{AB}u&=&\D_{z}u+\pa^2_{ss}u-2\e\sum_j\Phi'_j\pa^2_{sj}u+2\e\langle{\bf
H},z+\Phi\rangle\pa^2_{ss}u
+\e^2\sum_{l,j}\Phi'_l\Phi'_j\pa^2_{lj}u\\
&+&4\e^2\langle{\bf H},z+\Phi\rangle^2\pa^2_{ss}u-\frac12\e^2
\sum_{m,l}\pa^2_{ml}g_{kj}(z_m+\Phi_m)(z_l+\Phi_l)
\pa^2_{kj}u\\
&-&\frac12\e^2\sum_{m,l}\pa^2_{ml} g_{1j}(z_m+\Phi_m)(z_l+\Phi_l)
\pa^2_{sj}u -  4\e^2\langle{\bf H},z+\Phi\rangle\Phi'_j\pa^2_{sj}u\\
&-&\frac12\e^2
\sum_{m,l}\pa^2_{ml}g_{11}(z_m+\Phi_m)(z_l+\Phi_l)\pa^2_{ss}u
+R_3(z,\Phi,\Phi')\left(\pa^2_{ss}u+\pa^2_{sj}u+\pa^2_{lj}u\right)
\end{eqnarray*}
\begin{eqnarray*}
 \sum_{A,B} \pa_A\tilde{g}^{AB} \pa_Bu& =& -\e^2\sum_j\Phi_j''\pa_j
u-\e^2\sum_{i,j,l}\pa^2_{kl}g_{kj}(z_l+\Phi_l)\pa_ju-2\e^3\langle{\bf
H},z+\Phi\rangle\sum_{j}\Phi_j''
\pa_j u \\
&+& 2  \e^2 \langle{\bf H}',z+\Phi\rangle \pa_s u - \e^2
\sum_{l,j} \pa^2_{lj} g_{1j} (z_l + \Phi_l) \pa_s u
+R_3(z,\Phi,\Phi',\Phi'')(\pa_su+\pa_j u)
\end{eqnarray*}

\begin{eqnarray*}
\sum_{A,B}\frac{1}{\sqrt{\det \tilde{g}}} \tilde{g}^{AB} \left(
\partial_B
  \sqrt{\det \tilde{g}} \right) \partial_A u&=&-\e\sum_j H^j\pa_j
  u-\e^2\langle{\bf H},z+\Phi\rangle\sum_jH^j\pa_j u-\e^2\langle{\bf H}',z+\Phi\rangle\pa_su\\
  &+&\frac12\sum_j\pa^2_{lj}g_{11}(z_l+\Phi_l)\pa_ju
  +R_3(z,\Phi,\Phi')(\pa_s u+\pa_j u)
\end{eqnarray*}
The result then follows by collecting these three terms.
\end{pf}

\subsection{Expansion at first order in $\e$}\label{s:exp1st}

In this subsection we solve formally equation \eqref{eq:new} up to
order $\e$, discarding the terms which turn out to be of order
$\e^2$ and higher. \

 For the approximate solution as in \eqref{eq:prof}, we make a
more precise {\em ansatz} of the following form
  \begin{equation}\label{eq:u0}
  \psi_{1,\e}(s,z) = e^{- i \frac{\widetilde{f}_0(\e s)}{\e}} \left\{ h(\e s) U \left( k( \e
  s) z \right) + \e \left[ w_r + i w_i \right] \right\}, \qquad \qquad s
  \in [0, 2 \pi], y \in \R^{n-1},
\end{equation}
where $\widetilde{f}_0(\e s)=f(\e s)+\e f_1(\e s)$.  By direct
computation, the first and second derivatives of $\psi_{1,\e}$
satisfy
  \begin{eqnarray*}
\partial_s \psi_{1,\e} & = & e^{- i \frac{\widetilde{f}_0(\e s)}{\e}} \left[
  - i \widetilde{f}_0'(\e s) h(\e s) U (k (\e
  s) z) + \e h'(\e s) U (k(\e s)
  z) + \e h(\e s) k'(\e s) \n U
  (k (\e s) z) \cdot z \right] \\ & + & e^{- i \frac{\widetilde{f}_0(\e s)}{\e}}
  \left[ - i \e f' w_r + \e f' w_i \right] + O(\e^2);
\end{eqnarray*}
  \begin{equation*}
\partial_{i} \psi_{1,\e} = e^{- i \frac{\widetilde{f}_0(\e s)}{\e}} \left[ h(\e s)
k (\e s )\partial_{i} U \left( k( \e s) z \right) + \e
\partial_{i} w_r + i \e \partial_{i} w_i \right];
\end{equation*}
\begin{eqnarray*}
    \partial^2_{ss} \psi_{1,\e} & = & - (\widetilde{f}_0')^2 h
    U (k z) e^{- i \frac{\widetilde{f}_0(\e s)}{\e}}
    - i\e e^{- i \frac{\widetilde{f}_0(\e s)}{\e}} \left[
    f'' h U (k z) + 2 f' h'
    U (k z) + 2 f' h k' \n U
    (k z) \cdot z \right] \nonumber \\
    & - & \e f'^2e^{- i \frac{\widetilde{f}_0(\e s)}{\e}} [w_r + i w_i] + O(\e^2);
\end{eqnarray*}
  \begin{equation*}
\partial^2_{lj} \psi_{1,\e} = e^{- i \frac{\widetilde{f}_0(\e s)}{\e}} \left[ h(\e
s) k^2 (\e s )
  \partial^2_{lj} U \left( k( \e s) z \right) + \e
  \partial^2_{lj} w_r + i \e \partial^2_{lj} w_i \right];
\end{equation*}
\begin{eqnarray*}
  \partial^2_{sj} \psi_{1,\e}  & = & e^{-
  i \frac{\widetilde{f}_0(\e s)}{\e}} \bigg[ - i \widetilde{f}_0'(\e s) h(\e s)
  + \e h'(\e s)  \bigg]k (\e s) {\partial_j}
  U(k z) \nonumber \\ & + & \e e^{- i \frac{\widetilde{f}_0(\e s)}{\e}} h (\e s)
  k'(\e s) \left[ k \sum_l \partial^2_{lj}
  U (k z) z_l + \partial_j U (k z) \right] \\ & - &
   i \e f'(\e s) \partial_j w_r (\e s, z) e^{- i \frac{\widetilde{f}_0(\e s)}{\e}}
   + \e f'(\e s) \partial_j w_i (\e s, z) e^{- i \frac{\widetilde{f}_0(\e s)}{\e}} +
   O(\e^2). \nonumber
\end{eqnarray*}
Similarly, the potential $V$ satisfies
$$
  V(\e x) = V(\e s) + \e \langle \n^N V, z + \Phi \rangle + O(\e^2).
$$
Expanding \eqref{eq:new} in powers of $\e$, we obtain
$$
 e^{i \frac{\widetilde{f}_0(\e s)}{\e}} \bigg( - \D_g \psi_{1,\e}
 + V(\e x) \psi_{1,\e} - |\psi_{1,\e}|^{p-1}
\psi_{1,\e}\bigg) = \e \mathcal{R}_r
  + i \e \mathcal{R}_i + O(\e^2),
$$
with
\begin{equation}\label{eq:Rr}
    \mathcal{R}_r = \mathcal{L}_r w_r +2f'f_1'hU+ 2 f'^2 h U(kz) \langle
    {\bf H}, z + \Phi \rangle +
    h k \langle {\bf H}, \n U(kz) \rangle +  \langle \n^N V, z + \Phi \rangle h U(kz);
\end{equation}
\begin{equation}\label{eq:Ri}
    \mathcal{R}_i = \mathcal{L}_i w_i + \left[ f'' h U(kz) + 2 f' h' U(kz) + 2
    f' h k' \n U(kz) \cdot z \right] - 2 \sum_j [\Phi'_j f' h k \pa_j
    U(kz)],
\end{equation}
and where we have defined the two operators $\mathcal{L}_r$ and
$\mathcal{L}_i$ as
$$
  \mathcal{L}_r w = - \D_z w + (V + f'^2) w - p h^{p-1}
    U(kz)^{p-1} w;
$$
$$
  \mathcal{L}_i w = - \D_z w + (V + f'^2) w - h^{p-1}
    U(kz)^{p-1} w.
$$
It is well-known, see for example \cite{o}, that the kernel of
$\mathcal{L}_r$ is generated by the functions $\partial_2 U (k
\cdot), \dots, \partial_n U (k \cdot)$, while that of
$\mathcal{L}_i$ is one-dimensional and generated by $U(k \cdot)$.

We choose the functions $w_r$ and $w_i$ in such a way that
$\mathcal{R}_r$ and $\mathcal{R}_i$ vanish. Since $\mathcal{L}_r$
is Fredholm, the solvability condition for $w_r$ is that the
right-hand side of this equation is orthogonal in $L^2(\R^{n-1})$
to $\partial_2 U (k \cdot), \dots, \partial_n U (k \cdot)$.
Therefore, to get solvability, we should multiply the right-hand
side by each of these functions and get $0$. The same holds true
for $w_i$, but replacing the functions $\partial_{z_j} U (k
\cdot)$ by $U(k \cdot)$. We discuss the solvability in $w_i$
first. Writing this equation as $\mathcal{L}_i w_i =
\mathfrak{f}$, we can multiply it by $U (k \cdot)$ and use the
self-adjointness of $\mathcal{L}_i$ to get
$$
  0 = \int_{\R^{n-1}} w_i \mathcal{L}_i U(k \cdot) = \int_{\R^{n-1}}
  U(k \cdot) \mathcal{L}_i w_i = \int_{\R^{n-1}} \mathfrak{f} U(k \cdot).
$$
Following the computations of Subsection \ref{ss:cand}, this
condition yields
$$
  f'' h k^{-(n-1)} + 2 f' h' k^{-(n-1)} = (n-1) h k' k^{-n} f',
$$
which implies
$$
  f' = \mathcal{A} \frac{k^{n-1}}{h^2} = \mathcal{A} h^\s.
$$
This equation is nothing but \eqref{eq:f'C}, and hence the
solvability is guaranteed. Since $\mathcal{L}_i$ clearly preserves
the parity in $z$, we can decompose $w_i$ in its even and odd
parts as
$$
  w_i = w_{i,e} + w_{i,o},
$$
with $w_{i,e}$ and $w_{i,o}$ solving respectively the equations
$$
  \mathcal{L}_i w_{i,e} = - \left[ f'' h U(kz) + 2 f' h' U(kz) + 2
    f' h k' \n U(kz) \cdot z \right]; \qquad \mathcal{L}_i w_{i,o} = 2 \sum_j
    \left[ \Phi'_j f' h k \partial_{j} U (kz) \right],
$$
where the right-hand sides are respectively the even and odd parts
of the datum in \eqref{eq:Ri}. We notice that, since the kernel of
$\mathcal{L}_i$ consists of even functions, only the even part of
the equation plays a role in the solvability, since the product
with the odd part vanishes by oddness.

Indeed, \eqref{eq:Rr} and \eqref{eq:Ri} can be solved explicitly,
and the solutions are given by
$$
  w_{i,e} = \frac{p-1}{4} f' h' |z|^2 U(kz); \qquad \qquad
  w_{i,o} = - \sum_j \Phi'_j f' h z_j U(kz).
$$
In fact, as one can easily check, we have the following relations
$$
  \mathcal{L}_i (z_j U(kz)) = - 2 k \partial_j U(kz); \qquad \qquad
  \mathcal{L}_i (|z|^2 U(kz)) = - 2 (n-1) U(kz) - 4 k \n U(kz) \cdot z,
$$
which imply the above claim (here we also used \eqref{eq:hk} and
some manipulations).

\

Turning to $w_r$, if we multiply by $\partial_j U$, we integrate by
parts and use some scaling, we find that the following conditions
holds true, for $j = 2, \dots, n$
$$
  2 H^j \left( (f')^2 \int_{\R^{n-1}} U(z)^2 dz -
  \frac{k^2}{n-1} \int_{\R^{n-1}} |\n U(z)|^2 dz \right) +
  \langle \n^N V, E_j \rangle \int_{\R^{n-1}} U(z)^2 dz = 0.
$$
Using \eqref{eq:hk}, we get equivalently
$$
  2 H^j \left( \mathcal{A}^2 h^{2 \s} \int_{\R^{n-1}} U(z)^2 dz -
  \frac{h^{p-1}}{n-1} \int_{\R^{n-1}} |\n U(z)|^2 dz
  \right) + \langle \n^N V, E_j \rangle \int_{\R^{n-1}} U(z)^2 dz = 0,
  \quad j = 2, \dots, n.
$$
{From} a Pohozaev-type identity (playing with \eqref{eq:ovv} and
integrating by parts) one finds
\begin{equation}\label{eq:poho1}
    \int_{\R^{n-1}} |\n U(z)|^2 dz = \frac{(n-1)(p-1)}{(3-n)(p+1)
  + 2 (n-1)} \int_{\R^{n-1}} U(z)^2 dz = \frac{(n-1)(p-1)}{2 \th}
  \int_{\R^{n-1}} U(z)^2 dz.
\end{equation}
Using this formula the solvability condition then becomes
$$
  H^j \left( (p-1) \frac{h^{p-1}}{\th} - 2 \mathcal{A}^2 h^{2 \s}
  \right) = \langle \n^N V, E_j \rangle, \qquad j = 2,
  \dots, n,
$$
which is nothing but the stationary condition \eqref{eq:euler}.
Therefore, since we are indeed assuming this condition, also the
solvability for $w_r$ is guaranteed. As for $w_i$, we can
decompose $w_r$ in its even and odd parts as
$$
  w_r = w_{r,e} + w_{r,o},
$$
where $w_{r,e}$ and $w_{r,o}$ solve respectively
$$
  \mathcal{L}_r w_{r,e} = -2f'f_1'hU- 2 (f')^2 h U(kz) \langle {\bf H},
  \Phi \rangle - \langle \n^N V, \Phi \rangle h U(kz);
$$
$$
  \mathcal{L}_r w_{r,o} = - 2 (f')^2 h U(kz) \langle {\bf H}, z \rangle
  - h k \sum_j H^j \partial_j U(kz) - \langle \n^N V, z \rangle h
  U(kz).
$$
Using the Euler equation, one gets
$$
\mathcal{L}_r w_{r,o} =-h\sum_j H^j\left(k\pa_j U+ h^{p-1}
\frac{p-1}{\th}z_jU \right).
$$
It is also convenient to have the explicit expression of $w_r$. We
notice first that
$$
  \mathcal{L}_r \left( - \frac{1}{(p-1) h^{p-1}} U(kz) -
  \frac{1}{2k} \n U(kz) \cdot z \right) = U(kz).
$$
Hence it follows
$$
  w_{r,e} = \left[h \langle \n^N V + 2 (f')^2 {\bf H}, \Phi \rangle +2f'f_1'h\right] \left(
  \frac{1}{(p-1) h^{p-1}} U(kz) + \frac{1}{2k} \n U(kz) \cdot z \right).
$$
Using \eqref{eq:euler} we finally find
$$
  w_{r,e} = \left[\frac{p-1}{\th} h^p \langle {\bf H}, \Phi \rangle+2f'f_1'h\right] \left(
  \frac{1}{(p-1) h^{p-1}} U(kz) + \frac{1}{2k} \n U(kz) \cdot z \right).
$$
By the above computations we obtain the following result.

\begin{lem}\label{l:solve} Suppose $h(\ov{s})$ and $f(\ov{s})$ satisfy
\eqref{eq:ovhovk3} and \eqref{eq:f'Cintr} for some $\mathcal{A} >
0$: assume also that the curve $\g$ verifies \eqref{eq:eulerintr}.
Then there exist two smooth functions $w_r(\ov{s},z), w_i(\ov{s},z)$
for which the terms $\mathcal{R}_r$ and $\mathcal{R}_i$ in
\eqref{eq:Rr}-\eqref{eq:Ri} vanish identically. Therefore, the
function $\psi_{1,\e}$ in \eqref{eq:u0} satisfies \eqref{eq:new} up
to an error $O(\e^2)$.
\end{lem}

\subsection{Expansions at second order in $\e$}\label{ss:fe}

Next we compute the terms of order $\e^2$ in the above expression.
Adding a correction $\e^2 [v_r + i v_i]$ to the function in
\eqref{eq:u0} we define an approximate solution of the form
  \begin{equation}\label{eq:u02}
  \psi_{2,\e}(s,z) = e^{- i \frac{\widetilde{f}_0(\e s)}{\e}}
  \left\{ h(\e s) U \left( k( \e
  s) z \right) + \e \left[ w_r + i w_i \right] +
  \e^2 \left[ v_r + i v_i \right]\right\}; \qquad \qquad s
  \in [0, 2 \pi], y \in \R^{n-1},
\end{equation}
where $\tilde{f}_0=f(\e s)+\e f_1(\e s)$. The first and second
derivatives of $\psi_{2,\e}$ are given by
  \begin{eqnarray*}
e^{ i \frac{\tilde{f}_0(\e s)}{\e}}\partial_s \psi_{2,\e} & = &
\left[
  - i \tilde{f}_0'(\e s) h(\e s) U (k (\e
  s) z) + \e h'(\e s) U (k(\e s)
  z) + \e h(\e s) k'(\e s) \n U
  (k (\e s) z) \cdot z \right] \\ & + &
  \left[ - i \e \tilde{f}_0' w_r + \e \tilde{f}_0' w_i \right]
  +\left[ - i \e^2 \tilde{f}_0' v_r + \e^2 \tilde{f}_0' v_i \right]+
   \e^2(\pa_sw_r+iw_i)+O(\e^3);
\end{eqnarray*}
  \begin{equation*}
e^{ i \frac{\tilde{f}_0(\e s)}{\e}}\partial_{j} \psi_{2,\e} =
\left[ h(\e s) k (\e s )\partial_{j} U \left( k( \e s) z \right) +
\e\partial_{j} w_r + i \e \partial_{j} w_i +\e^2 \partial_{j} v_r
+ i \e^2 \partial_{j} v_i\right];
\end{equation*}
\begin{eqnarray*}
   e^{i \frac{\tilde{f}_0(\e s)}{\e}} \frac{\partial^2
   \psi_{2,\e}}{\partial s^2} & = & - (\tilde{f}_0')^2 h
    U (k z)
    - i\e  \left[\tilde{f}_0'' h U (k z) + 2 \tilde{f}_0' h'
    U (k z) + 2 \tilde{f}_0' h k' \n U
    (k z) \cdot z \right] -\e \tilde{f}_0'^2 [w_r + i w_i]\nonumber
    \\[3mm]
    & - & \e^2 \tilde{f}_0'^2 [v_r + i v_i]+\e^2\left[ 2\tilde{f}_0'\pa_s
    w_i+h''U(kz)+2h'k'\n U\cdot z+hk''\n U\cdot z \right. \\ [3mm] & + & \left.
    hk'^2\n^2U(kz)[z,z] +\tilde{f}_0''
    w_i \right]
    - i \e^2\left[ 2\tilde{f}_0'\pa_s w_r+\tilde{f}_0''w_r
    \right]+O(\e^3); \nonumber
\end{eqnarray*}
  \begin{equation*}
e^{ i \frac{\tilde{f}_0(\e s)}{\e}}\partial^2_{lj} \psi_{1,\e} =
\left[ h(\e s) k^2 (\e s )
  \partial^2_{lj} U \left( k( \e s) z \right) + \e
  \partial^2_{lj} w_r + i \e \partial^2_{lj} w_i + \e^2
  \partial^2_{lj} v_r + i \e^2 \partial^2_{lj} v_i\right];
\end{equation*}
\begin{eqnarray*}
 e^{ i \frac{\tilde{f}_0(\e s)}{\e}} \partial^2_{sj}
 \psi_{2,\e}  & = &  \bigg[ - i \tilde{f}_0'(\e s) h(\e s)k (\e s)
  + \e h'(\e s)k (\e s)+\e  h (\e s)k'(\e s  \bigg] {\partial_j}U(k z) \nonumber \\
  & + & \e  h (\e s)k(\e s)k'(\e s) \sum_l
  \partial^2_{lj}U (k z) z_l + \partial_j U (k z) -i\e\tilde{f}_0'\pa_j w_r \\
   & + & \e \tilde{f}_0'(\e s) \partial_j
   w_i-i\e^2\tilde{f}_0'\pa_j v_r+\e^2 \tilde{f}_0'(\e s) \partial_j v_i
   + \e^2  \pa^2_{sj} w_r+i \e^2  \pa^2_{sj} w_r. \nonumber
\end{eqnarray*}
We also have the formal expansion
\begin{eqnarray*}
    e^{i \frac{\tilde{f}_0(\e s)}{\e}} |\psi_{2,\e}|^{p-1} \psi_{2,\e} & = &
    h^p|U|^{p-1} U + p\e h^{p-1}|U|^{p-1} w_r + i \e h^{p-1}|U|^{p-1} w_i + \frac
    12 p (p-1)\e^2 h^{p-2}|U|^{p-3} U w_r^2 \\ & + &
    \frac 12 (p-1)\e^2 h^{p-2}|U|^{p-3} U w_i^2 + i (p-1)\e^2 h^{p-2} |U|^{p-3}
    U w_r w_i \\ & + & p\e^2 h^{p-1}|U|^{p-1} v_r + i \e^2 h^{p-1}|U|^{p-1} v_i +
    O(\e^3).
\end{eqnarray*}
Similarly, expanding $V$ up to order $\e^2$, we have
$$
V(\e x)=V(\e s)+\e \langle \n^N V,z+\Phi \rangle +\frac
12\e^2(\n^N)^2V[z+\Phi,z+\Phi]+R_3(z,\Phi).
$$

Using the expansions of Subsection \ref{ss:cand}, we obtain
\begin{eqnarray*}
 e^{i \frac{\tilde{f}_0(\e s)}{\e}} \bigg( - \D_g \psi_{2,\e}
 + V(\e x) \psi_{2,\e} - |\psi_{2,\e}|^{p-1}
\psi_{2,\e}\bigg)&=& \e^2 (\tilde{\mathcal{R}}_r + i
\tilde{\mathcal{R}}_i)
\\ & = & \e^2(\tilde{R}_{r,e}+\tilde{R}_{r,o})
+\e^2i(\tilde{R}_{i,e}+\tilde{R}_{i,o})\\
&+&\e^2(\tilde{R}_{r,e,f_1}+\tilde{R}_{r,o,f_1})
+\e^2i(\tilde{R}_{i,e,f_1}+\tilde{R}_{i,o,f_1})
\\ & + & \e^2 \mathcal{L}_r v_r + \e^2 i \mathcal{L}_i v_i +
O(\e^3),
\end{eqnarray*}
where
\begin{eqnarray}\label{eq:Rre} \nonumber
    \tilde{R}_{r,e} & = & - \frac 12 (f')^2 h U(kz)  \sum_{l,m} \partial^2_{lm}
    g_{11} (z_m z_l + \Phi_m \Phi_l)  +2(f')^2 \langle {\bf H},w_{r,e}\Phi  + w_{r,o}z\rangle  \\
    & + & 4(f')^2 h U(kz)\left[ \nonumber
    \langle {\bf H},z\rangle^2 + \langle {\bf H},\Phi
    \rangle^2 \right] \\ \nonumber  & - & \left[ h'' U(kz) + 2 h' k'
    \n U(kz) \cdot z + h k'' \n U(kz) \cdot z + h (k')^2 \n^2 U (kz) [z,z]
    \right] \\ & + & 2 f' \partial_s w_{i,e} + f'' w_{i,e} + 2
    \Phi'_j f' \partial_{j} w_{i,o} + \left[ \frac 12 \sum_{l,m} \partial^2_{lm}
    g_{tj} (z_m z_l + \Phi_m \Phi_l) - \Phi'_t \Phi'_j \right] h k^2
    \partial^2_{tj} U(kz) \\ \nonumber  & + &h k  \sum_{l,m,j} \partial^2_{lm} g_{mj} z_l
     \partial_{j} U(kz) + \sum_m H^m \partial_{m} w_{r,o}
    + h k \langle {\bf H},z\rangle \sum_m H^m  \partial_{m} U(kz) \\ \nonumber  & + &
     k h\sum_m \left[\langle {\bf H},z\rangle H^m - \frac 12
    \sum_l \partial^2_{ml} g_{11} z_l \right] \partial_{m} U(kz) \\ \nonumber  & - &
    \frac 12 p (p-1) h^{p-2} U(kz)^{p-2} (w_{r,e}^2 + w_{r,o}^2) - \frac 12 (p-1) h^{p-2}
    U(kz)^{p-2} (w_{i,e}^2 + w_{i,o}^2) \\ \nonumber  & + & \langle \n^N V,w_{r,o}z +w_{r,e}\Phi\rangle
      + \frac 12 \sum_{m,j} \partial^2_{mj} V (z_m z_j + \Phi_m \Phi_j) h U(kz);
\end{eqnarray}
\begin{eqnarray}\label{eq:Rro} \nonumber
    \tilde{R}_{r,o} & = & - (f')^2 h U(kz) \sum_{l,m} \partial^2_{lm}
    g_{11} z_m \Phi_l + 8(f')^2 h U(kz) \langle
    {\bf H},z\rangle \langle {\bf H},\Phi\rangle \\ \nonumber  & + &
    2(f')^2   \langle {\bf H},w_{r,e}z+w_{r,o}\Phi \rangle
     - 2 f' \partial_s w_{i,o} - f'' w_{i,o} + 2
    h' k\sum_j\Phi'_j  \partial_{j} U(kz) \\ \nonumber  & + & 2 h' k\sum_j\Phi'_j \left[
    k \sum_l \partial^2_{lj} U(kz) z_l + \partial_j U(kz)
    \right] + 2 f'\sum_j\Phi'_j\partial_{j} w_{i,e} +
    h k\langle {\bf H},\Phi \rangle \sum_m H^m \partial_{m} U(kz) \\
    & + & \left[
\frac 12 \sum_{l,m} \partial^2_{lm} g_{tj} (z_m \Phi_l + z_l
\Phi_m) \right] h k^2 \partial^2_{tj} U(kz) +
hk\sum\Phi_j''\partial_{j} U(kz) + \left( \sum_{j,l,m}
\partial^2_{lm} g_{mj} \Phi_l
    \right) h k \partial_{j} U(kz) \\ \nonumber  & +&  \sum_m H^m \partial_{m} w_{r,e} +
    h k \sum_m \left[ \langle {\bf H},\Phi \rangle H^m - \frac 12
    \sum_l \partial^2_{ml} g_{11} \Phi_l \right] \partial_{m}
    U(kz) \\ & - & p (p-1) h^{p-2} U(kz)^{p-2} w_{r,e} w_{r,o} -
    (p-1) h^{p-2} U(kz)^{p-2} w_{i,e} w_{i,o} \nonumber
    \\ & + &\langle \n^N V,w_{r,e}z +w_{r,o}\Phi\rangle \nonumber
    +
    \sum_{j,l} \partial^2_{jl} V z_j \Phi_l h U(kz);
\end{eqnarray}
\begin{eqnarray}\label{eq:Rie} \nonumber
    \tilde{R}_{i,e} & = & 2\left[ f'' h U(kz) + 2 f' h' U(kz) + 2 f' h k' \n U(kz)
    \cdot z \right]\langle {\bf H},\Phi \rangle  \\ \nonumber
    & + &2 (f')^2 \langle {\bf H},  w_{i,e}\Phi +w_{i,o} z  \rangle
     + 2 f' \partial_s w_{r,e} + f'' w_{r,e} - 2f' \sum_j\Phi'_j\partial_{j} w_{r,o} \\
     & -& 2 f' h k \sum_j\partial_{j} U(kz)
    \left[  2\langle {\bf H},z \rangle  \Phi'_j
    + \frac 12 \sum_{l,m} \partial^2_{lm} g_{1j} (z_m \Phi_l + z_l \Phi_m)
    \right] \\ & - & f' h U(kz) \left( \sum_m ( \partial^2_{1m} g_{11} \Phi_m
    -2\langle {\bf H},\Phi' \rangle   ) \right) -f' h \left[2\langle {\bf H},\Phi' \rangle
    +\sum_{j,l}\partial^2_{lj} g_{1j} \Phi_l \right]  U(kz) \nonumber  \\ & + &
    \frac12 f' h\left(\sum_l \partial_{1l} g_{11} \Phi_l \right)  U(kz)
    +\sum_jH^j \partial_{j} w_{i,o} \nonumber  \\ & - &
    (p-1) h^{p-2} U(kz)^{p-2} (w_{r,e} w_{i,e}
    + w_{r,o} w_{i,o}) +\langle \n^N V,w_{i,o} z +w_{i,e}\Phi\rangle; \nonumber
\end{eqnarray}
\begin{eqnarray}\label{eq:Rio} \nonumber
    \tilde{R}_{i,o} & = & 2 \left[ f'' h U(kz) + 2 f' h' U(kz) + 2 f' h k'
    \n U(kz) \cdot z \right] \langle {\bf H},z \rangle
     + \sum_i H^j\partial_{j} w_{i,e} \\ \nonumber
    & + & 2(f')^2\langle {\bf H},w_{i,e}z+ w_{i,o}\Phi \rangle
     + 2 f' \partial_s w_{r,o} + f'' w_{r,o}-2f' \sum_j\Phi'_j
    \partial_{j} w_{r,e} \\ & -&   2 f' h k \sum_j\partial_{j} U(kz)
    \left[ 2\langle {\bf H},\Phi \rangle  \Phi'_j
    + \frac 12 \sum_{l,m} \partial^2_{lm} g_{1j} (z_m z_l + \Phi_l \Phi_m)
    \right] \\ & - & f' h U(kz) \left( \sum_m \partial^2_{1m} g_{11} z_m
    \right) -f'h \left( \sum_{j,l} \partial^2_{lj} g_{1j} z_l \right)  U(kz) + \frac 12f' h
\left( \sum_l \partial^2_{1l} g_{11} z_l
    \right)  U(kz) \nonumber  \\ & - & (p-1) h^{p-2} U(kz)^{p-2} (w_{r,e} w_{i,o}
    + w_{r,o} w_{i,e}) + \langle \n^N V,w_{i,e}z + w_{i,o} \Phi\rangle. \nonumber
\end{eqnarray}
We used the notation $\tilde{R}_{r,e,f_1}$, $\tilde{R}_{r,o,f_1}$,
$\tilde{R}_{i,e,f_1}$ and $\tilde{R}_{i,o,f_1}$ for the terms
involving $f_1$, namely
\begin{eqnarray} \label{eq:Rref1} \nonumber
  \tilde{R}_{r,e,f_1} & = & (f_1')^2 hU + 2f'f_1'w_{r,e}+4\langle {\bf{H}},
  \Phi \rangle f'f_1'hU - 2 p (p-1) h^{p-2}|U|^{p-2}  h^2f'^2f_1'^2\tilde{U}^2 \\
  & + & 2 f' f'_1 h \langle \n^N V, \Phi \rangle \tilde{U} + 4 \langle {\bf H},
  \Phi \rangle (f')^3 f'_1 h \tilde{U} - 2 p \frac{(p-1)^2}{\th} h^{2p-1}
  f' f'_1 \langle {\bf H}, \Phi \rangle U^{p-2} \tilde{U}^2;
\end{eqnarray}
\begin{eqnarray} \nonumber
  \tilde{R}_{r,o,f_1} & = & 2f'f_1'w_{r,o}+4\langle {\bf{H}},z
  \rangle f'f_1'hU
  - 2 p(p-1) f' f'_1 h^{p-1} U^{p-2} \tilde{U} w_{r,o} \\
   & + & 2 f' f'_1 h \langle \n^N V, z \rangle \tilde{U}  + 2 H^j
  f' f'_1 h k \pa_j \tilde{U} + 4 \langle {\bf H}, z \rangle
  (f')^2 f' f'_1 h \tilde{U};
\end{eqnarray}
\begin{eqnarray} \nonumber
  \tilde{R}_{i,e,f_1} & = & 2h'f_1'U+2hf_1'k'\n
U\cdot z+2f'f_1'w_{i,e}+f_1''hU+4f'\pa_s(hf'f_1'\tilde{U}) \\
  & + & 2f''hf'f_1'\tilde{U}-2(p-1) h^{p-1} |U|^{p-2} f'f_1'
  \tilde{U} w_{i,e};
\end{eqnarray}
\begin{equation}\label{eq:Riof2}
\tilde{R}_{i,o,f_1}= 2f'f_1' w_{i,o}-2(p-1) h^{p-1}
|U|^{p-2} f'f_1' \tilde{U} w_{i,o} - 4 (f')^2 h k f'_1
\Phi'_j \partial_j \tilde{U} - 2 h k f'_1 \Phi'_j \pa_j U,
\end{equation}
where we wrote for brevity
$\tilde{U}=\frac{1}{(p-1)h^{p-1}}U(kz)+\frac{1}{2k}\n U(kz)\cdot z$.
Again, we collect the results of this section in one proposition.

\begin{pro}\label{p:errors} Suppose $\Phi, f_1$ are smooth functions
on $[0,L]$, let $z$ be the normal coordinates given in
\eqref{eq:defzz} and $w_i, w_r$ be be as in Lemma \ref{l:solve}.
Then, if $\psi_{2,\e}$ is as in \eqref{eq:u02}, in the coordinates
$(s,z)$ we have
\begin{eqnarray}\label{eq:exp2nd} \nonumber
 & & - \D_{\tilde{g}_\e}
\psi_{2,\e} + V(\e x) \psi_{2,\e} - |\psi_{2,\e}|^{p-1} \psi_{2,\e}  = \e^2 (\mathcal{L}_r
v_r + i \mathcal{L}_i v_i) \\ & + &
\e^2 (\tilde{R}_{r,e} + \tilde{R}_{r,o} + i \tilde{R}_{i,e} + i \tilde{R}_{i,o} +
\tilde{R}_{r,e,f_1} + \tilde{R}_{r,o,f_1} + i \tilde{R}_{i,e,f_1} + i \tilde{R}_{i,o,f_1}) +
o(\e^2),
\end{eqnarray}
where the above error terms are given respectively in
\eqref{eq:Rre}-\eqref{eq:Riof2}.
\end{pro}

\section{Final expansions and proof of Theorem \ref{t:approx}}\label{s:proj}

In this Section we prove existence of approximate solutions to
\eqref{eq:pe} up to any power of $\e$. To do this, we first evaluate
the projections of the error terms on the kernels of the operators
$\mathcal{L}_i$ and $\mathcal{L}_r$, and adjust $\Phi$ and $f_1$ so
that these projections vanish, in order to obtain solvability of
\eqref{eq:new} up to order $\e^2$. Then, with an iterative
procedure, we turn to the general case.

\subsection{Projection onto the kernel of
$\mathcal{L}_i$}\label{ss:prkerLi}

If one wants to find $v_i$ so that the imaginary terms in
\eqref{eq:exp2nd} vanish, by Fredholm's alternative the imaginary
$\tilde{R}$'s must be orthogonal for every $s$ to the kernel of
$\mathcal{L}_i$, which is given by  $i U(k(\e s) \cdot)$. To compute
this projection, by parity reasons, we need to multiply
$\tilde{R}_{i,e}$ and $\tilde{R}_{i,e,f_1}$ by $U(k(\e s) \cdot)$
and to integrate over $\R^{n-1}$. We evaluate the two terms
separately.

\

\noindent {\bf Contribution of $\tilde{R}_{i,e}$}. After some manipulation
we obtain
$$
  \tilde{R}_{i,e} = C_1 + C_2 + C_3 + C_4 + C_5,
$$
where
\begin{eqnarray*}
    C_1 & = & 2 \left[ f'' h U(kz) + 2 f' h' U(kz) + 2 f' h k' \n U(kz)
    \cdot z \right] \langle{\bf H},\Phi\rangle
    \\ & - & f' h U(kz) \left( \sum_m ( \partial^2_{1m} g_{11} \Phi_m
    -2\langle{\bf H},\Phi'\rangle  ) \right) + \frac 12f' h \left(
    \sum_l \partial_{1l} g_{11} \Phi_l \right)  U(kz);
\end{eqnarray*}
$$
  C_2 = - 2 f' h k\sum_j \partial_{j} U(kz) \left( \sum_{l,m}
  \partial^2_{lm} g_{1j} z_m \Phi_l \right) -f' h
    \sum_{j,l} \partial^2_{lj} g_{1j} \Phi_l U(kz);
$$
$$
  C_3 = -4 \sum_j f' h k \partial_{j} U(kz)
  \langle{\bf H},z\rangle  \Phi'_j -2f' h \sum_jH^j \Phi'_j  U(kz);
$$
\begin{eqnarray*}
    C_4 & = & - (f')^2   \langle{\bf H},w_{i,e}\Phi+w_{i,o}z\rangle
    + \sum_l H^l \partial_{l}w_{i,o} \\ & + & \langle \n^N V,w_{i,e}\Phi+w_{i,o}z\rangle
    - (p-1) h^{p-2} U(kz)^{p-2} (w_{r,e} w_{i,e} + w_{r,o}w_{i,o});
\end{eqnarray*}
$$
  C_5 = 2 f' \partial_s w_{r,e} + f'' w_{r,e} - 2 \sum_j\Phi'_j f'
    \partial_{j} w_{r,o}.
$$
We now evaluate $\int C_1 U(kz)$. Arguing as for the derivation of
\eqref{eq:f'C} we find the following identity
$$
  \int_{\R^{n-1}} U(kz) \left[ f'' h U(kz) + 2 f' h' U(kz) + 2 f' h k' \n U(kz)
    \cdot z \right] = 0,
$$
which implies
$$
  \int_{\R^{n-1}} C_1 U(kz) = f' h \left[ 2\langle{\bf H},\Phi'\rangle
  - \frac 12 \left( \sum_l \partial^2_{1l} g_{11} \Phi_l \right)
  \right] \int_{\R^{n-1}} U(kz)^2.
$$
On the other hand, using integration by parts, we obtain easily
$$
  \int_{\R^{n-1}} C_2 U(kz) = \int_{\R^{n-1}} C_3 U(kz) = 0.
$$
Turning to $C_4$, we recall that
$$
  \mathcal{L}_r w_r = -2(f')^2 h U(kz) \langle {\bf H},\Phi+z\rangle - h\langle \n^N V,\Phi+z\rangle
  U(kz) + \frac 12 \sum_j \partial_j g_{11} h k \partial_j U(kz).
$$
Therefore, after some computations we deduce
$$
  \int_{\R^{n-1}} C_4 U(kz) = - \frac 1h \int_{\R^{n-1}} (\mathcal{L}_r w_r) w_i -2 k\sum_j
  H^j \int_{\R^{n-1}} \partial_j U(kz) w_i - (p-1) h^{p-2}
  \int_{\R^{n-1}} U(kz)^{p-1} w_r w_i.
$$
Recalling the expressions of $\mathcal{L}_r$ and $\mathcal{L}_i$
and integrating by parts we then get
$$
  \int_{\R^{n-1}} C_4 U(kz) = - \frac 1h \int_{\R^{n-1}} (\mathcal{L}_i w_i) w_r -2k \sum_j
  H^j\int_{\R^{n-1}} \partial_j U(kz) w_i.
$$
From the definition of $w_i$ and using some cancelation we then find
\begin{eqnarray*}
    \int_{\R^{n-1}} (C_4 + C_5) U(kz) & = & \frac 1h \partial_s \left[ 2 h f' \int_{\R^{n-1}}
  w_{1} U(kz) \right] -2k \sum_j H^j\int_{\R^{n-1}} \partial_j
  U(kz) w_i \\ & = & \frac 1h \partial_s \left[ 2 h f' \int_{\R^{n-1}}
  w_{r,e} U(kz) \right] -2k \sum_j H^j \int_{\R^{n-1}} \partial_j
  U(kz) w_i.
\end{eqnarray*}
To evaluate the last integral we need to use the explicit expression
of $w_{r,e}$: recall that we have
$$
  w_{r,e} = h \left[    \langle \n^N V,\Phi \rangle
  +2(f')^2  \langle {\bf H},\Phi \rangle  \right] \left( \frac{1}{(p-1)
  h^{p-1}} U(kz) + \frac{1}{2k} \n U(kz) \cdot z \right).
$$
Therefore, adding all the terms and using some scaling we obtain
\begin{eqnarray*}
    \int_{\R^{n-1}} \tilde{R}_{i,e} U(kz) & = & \frac{f'h}{k^{n-1}}
    \left[ +2  \langle {\bf H},\Phi' \rangle  - \frac 12 \left(
    \sum_l \partial_{1l} g_{11} \Phi_l \right) \right]
    \int_{\R^{n-1}} U^2 -\frac{f'h}{k^{n-1}}
      \langle{\bf H},\Phi \rangle  \int_{\R^{n-1}} U^2
      \\ & + & \frac 1h \partial_s \left\{ 2 h^2 f'
    \left[ \langle \n^N V,\Phi \rangle
  +2 (f')^2  \langle {\bf H},\Phi \rangle \right] \frac{1}{k^{n+1}} \left(
    \frac{1}{p-1} - \frac{n-1}{4} \right)\right\} \int_{\R^{n-1}} U^2.
\end{eqnarray*}
From the Euler equation \eqref{eq:euler} it follows that
$$
 \langle \n^N V,E_m \rangle
  +2 (f')^2 H^m = \frac{p-1}{\th} h^{p-1} H^m.
$$
Therefore after some manipulation we find
$$
  \int_{\R^{n-1}} \tilde{R}_{i,e} U(kz) = -\frac{2\mathcal{A}}{h} \left( \frac{p-1}{2\th} - 1
  \right) \int_{\R^{n-1}} U^2 \partial_s \langle {\bf H},\Phi \rangle.
$$

\

\noindent {\bf Contribution of $\tilde{R}_{i,e,f_1}$.} Multiplying
$\tilde{R}_{i,e,f_1}$ by $hU(kz)$ and integrating, recalling the
expression of $w_{i}$ determined in subsection \ref{s:exp1st}
($w_{i,o}=-\sum_j \Phi_j'f'h z_jU(kz)$ and
$w_{i,e}=\frac{p-1}{4}f'h'|z|^2U(kz)$), we obtain

\begin{eqnarray*}
h\int_{\R^{n-1}}\tilde{R}_{i,e,f_1}U&=&2hh'f_1'
\int_{\R^{n-1}}U^2+2h^2f_1'k'\int_{\R^{n-1}}\n
U\cdot zU+\frac{p-1}{2}f'^2f_1'hh'\int_{\R^{n-1}}|z|^2U^2\\
&+&f_1''h^2\int_{\R^{n-1}}U^2+4h^2f'^2f_1''
\int_{\R^{n-1}}U\tilde{U}+4hf'f_1'\int_{\R^{n-1}}U\pa_s(hf'\tilde{U})
\\
&+&2h^2f''f'f_1'\int_{\R^{n-1}}\tilde{U}U-\frac{(p-1)^2}{2}
h^{p}h' f'^2f_1' \int_{\R^{n-1}}U^p|z|^2\tilde{U}
\end{eqnarray*}
We will need the following observations
$$
\int_{\R^{n-1}}\n U(kz)\cdot
zU(kz)dz=\frac{1}{k^n}\sum_j\int_{\R^n}\pa_jU(z)z_jU(z)dz=-\frac{n-1}{2k^n}\int_{\R^n}U^2(z)dz
=-\frac{n-1}{2k^n}\int_{\R^n}U^2(z)dz;
$$
which imply
\begin{equation}\label{eq:utu}
\int_{\R^{n-1}}U\tilde{U}dz=\int_{\R^{n-1}}\frac{1}{(p-1)h^{p-1}}U^2(kz)+\frac{1}{2k}\n
U(kz)\cdot z
U(kz)dz=-\frac{\s}{2(p-1)k^{n+1}}\int_{\R^{n-1}}U^2(z)dz;
\end{equation}
\begin{eqnarray*}
&&\int_{\R^{n-1}} |z|^2\n
 U(kz)\cdot zU^p(kz)dz=\frac{1}{k^{n+2}}\sum_{l,j}\int_{\R^n}z_l^2\pa_jU(z)z_jU^p(z)dz\\
 &=&-\frac{1}{k^{n+2}}\sum_{j}\int_{\R^n}(n+1)z_j^2
 U^{p+1}(z)dz-\frac{1}{k^{n+2}}\sum_{l,j}\int_{\R^n}pz_l^2\pa_jU(z)z_jU^p(z)dz\\
 &=&-\frac{n+1}{k}\int_{\R^n}|z|^2U^{p+1}(kz)dz-p\int_{\R^n}|z|^2\n
 U(kz)\cdot zU^p(kz)dz.
\end{eqnarray*}
From these we get
\begin{equation*}
\int_{\R^{n-1}} |z|^2\n U(kz)\cdot
zU^p(kz)dz=-\frac{n+1}{(p+1)k}\int_{\R^n}|z|^2
U^{p+1}(kz)dz=-\frac{n+1}{(p+1)k^{n+1}}\int_{\R^n}|z|^2U^{p+1}(z)dz.
\end{equation*}
Moreover one can easily see that
\begin{eqnarray*}
  2hh'f_1'\int_{\R^{n-1}}U^2dz+2h^2f_1'k'\int_{\R^{n-1}}\n
U\cdot zU dz+f_1''h^2\int_{\R^{n-1}}U^2dz & = & \pa_{\ov
s}\left(h^2f_1'
\int_{\R^n}U^{2}(kz)   \right)dz \\
  & = & \pa_{\ov s}\left(\frac{h^2f_1'}{k^{n-1}}\int_{\R^n}U^{2}(z)dz
\right),
\end{eqnarray*}
and
\begin{eqnarray*}
 4hf'f_1'\int_{\R^{n-1}}U\pa_s(hf'\tilde{U})dz&=&4\pa_{\ov
 s}\left(h^2f'^2f'_1\int_{\R^{n-1}}U\tilde{U}\right)dz+4f'h\pa_{\ov
 s}[hf'f_1']\int_{\R^{n-1}}U\tilde{U}dz \\ & + &
 4f'^2h^2f_1'\int_{\R^{n-1}}\tilde{U}\pa_{\ov s}Udz\\
 &=&4\pa_{\ov
 s}\left(h^2f'^2f'_1\int_{\R^{n-1}}U\tilde{U}\right)dz+4f'h\pa_{\ov
 s}[hf'f_1']\int_{\R^{n-1}}U\tilde{U}dz\\
 &+&4f'^2h^2f_1'k'\frac{1}{(p-1)h^{p-1}}\int_{\R^{n-1}}U(kz)\n
 U(kz)\cdot z dz\\ & + & 2f'^2h^2f_1'\frac{k'}{k}\int_{\R^{n-1}}|z|^2|\n
 U(kz)|^2dz.
\end{eqnarray*}
On the other hand, by \eqref{eq:unuup} (see Section
\ref{ss:apppoh}) we have
\begin{eqnarray*}
\frac{p-1}{2}f'^2f_1'hh'\int_{\R^{n-1}}|z|^2U^2 dz& = &
\frac{p-1}{p+1}f'^2f_1'hh' \int_{\R^{n-1}} |z|^2 U^{p+1}(kz)dz \\
& - &  \frac{(n-5)(p-1)}{2(n+1)}f'^2f_1'hh' \int_{\R^{n-1}} |z|^2
|\n U|^2 dz.
\end{eqnarray*}
Hence with some easy manipulations and by  \eqref{eq:nuuup} one
finds
\begin{eqnarray*}
&&\frac{p-1}{2}f'^2f_1'hh'\int_{\R^{n-1}}|z|^2U^2dz-\frac{(p-1)^2}{2}
h^{p}h' f'^2f_1'
\int_{\R^{n-1}}U^p|z|^2\tilde{U}dz+2f'^2h^2f_1'\frac{k'}{k}\int_{\R^{n-1}}|z|^2|\n
U(kz)|^2dz\\
&&=-\frac{(p-1)(n-1)}{4}f'^2f_1'hh'\left[
\frac{6}{n+1}\int_{\R^{n-1}}|z|^2|\n
U(kz)|^2dz-\frac{p-1}{p+1}\int_{\R^{n-1}} |z|^2 U^{p+1}(kz)dz \right]\\
&&=-\frac{(p-1)(n-1)^2}{4k^{n+1}}f'^2f_1'hh'\int_{\R^{n-1}}
U^{2}(z)dz.
\end{eqnarray*}
Collecting the above expressions, and letting $C_0=\int_{\R^{n-1}}
U^{2}(z)dz$ one proves that
\begin{eqnarray*}
\frac{h}{C_0}\int_{\R^{n-1}}{\tilde R}_{i,e,f_1}Udz&=& \pa_{\ov
s}\left(\frac{h^2f_1'}{k^{n-1}} \right)-2\pa_{\ov s}\left(f'^2h^2
f_1'
\frac{\s}{(p-1)k^{n+1}}\right)-\frac{(p-1)(n-1)^2}{4k^{n+1}}f'^2hh'f_1'\\
&+&2f'^2f_1'hh'\frac{\s}{(p-1)k^{n+1}}+(n-1)f'^2f_1'
\frac{hh'}{k^{n+1}}+f'f''f_1'h^2\frac{\s}{(p-1)k^{n+1}}\\
&=&\pa_{\ov s}\left(\frac{h^2f_1'}{k^{n-1}}\left[1- 2 f'^2
\frac{\s}{(p-1)k^{2}} \right]
\right)-\frac{\s}{(p-1)k^{n+1}}hf_1'\left( \s f'^2h'-f'f''
\right);
\end{eqnarray*}
using  \eqref{eq:f'C}, we have that $\s f'^2h'-f'f''=0$, which
implies
\begin{eqnarray*}
  \frac{h}{C_0}\int_{\R^{n-1}}{\tilde R}_{i,e,f_1}U(kz)& = &
  \pa_{\ov s}\left(\frac{h^2f_1'}{k^{n-1}}\left[1- 2 f'^2
\frac{\s}{(p-1)k^{2}} \right]  \right) \\
  & = & \pa_{\ov s}\left( \frac{h^2f_1'}{(p-1)k^{n+1}}
  \left[(p-1)h^{p-1}- 2 \s \mathcal{A}^2h^{2\s} \right] \right).
\end{eqnarray*}
Recalling the discussion at the beginning of this subsection, we
want to find $f_1$ such that
$$
\int_{\R^{n-1}}{\tilde R}_{i,e,f_1}U(kz)+\int_{\R^{n-1}} \tilde{R}_{i,e}
U(kz)
=0.
$$
This is equivalent to
\begin{equation}\label{eq:Tf1}
    T f_1 := \pa_{\ov s}\left( \frac{h^2f_1'}{(p-1)k^{n+1}}\left[(p-1)h^{p-1}-
2 \s \mathcal{A}^2h^{2\s} \right] \right) = 2\mathcal{A} \left(
\frac{p-1}{2\th} - 1 \right)  \partial_{\ov s} \langle {\bf H},\Phi \rangle.
\end{equation}
Hence we get
$$
f'_1=\frac{2\mathcal{A}(p-1)k^{n+1}}{(p-1)h^{p+1}- 2 \s
\mathcal{A}^2h^{2\s+2} } \left( \frac{p-1}{2\th} - 1
  \right)   \langle {\bf H},\Phi \rangle
  + c \frac{(p-1) k^{n+1}}{(p-1) h^{p+1} - 2 \s \mathcal{A}^2
  h^{2\s+2}},
$$
where $c$ is a constant to be chosen so that $\int_0^{L} f'_1 ds =
0$. Noticing that, by \eqref{eq:dhCV}, we have
$$
  \frac{(p-1) k^{n+1}}{(p-1) h^{p+1} - 2 \s \mathcal{A}^2
  h^{2\s+2}} = A \s h^{\s-1} \frac{\pa h}{\pa \mathcal{A}} + h^\s,
$$
the required condition becomes
$$
c = - \left( \frac{p-1}{2\th} - 1
  \right) 2\mathcal{A}(p-1) \frac{\int_0^{L}\frac{k^{n+1}}{(p-1)h^{p+1}- 2
\s \mathcal{A}^2h^{2\s+2} }   \langle {\bf H},\Phi \rangle
ds}{\int_0^L A \s
  h^{\s-1} \frac{\pa h}{\pa \mathcal{A}} + h^\s ds}.
$$
As one can easily check from \eqref{eq:C'} and \eqref{eq:euler}, $c$
coincides with $\mathcal{A}'$ and therefore we have in conclusion
\begin{equation}\label{eq:f1'final}
    f'_1=\frac{2\mathcal{A}(p-1)k^{n+1}}{(p-1)h^{p+1}- 2 \s
\mathcal{A}^2h^{2\s+2} } \left( \frac{p-1}{2\th} - 1
  \right)   \langle {\bf H},\Phi \rangle
  + \mathcal{A}' \frac{(p-1) k^{n+1}}{(p-1) h^{p+1} - 2 \s \mathcal{A}^2
  h^{2\s+2}}.
\end{equation}

\subsection{Projection onto the kernel of $\mathcal{L}_r$}\label{ss:hp}

Similarly to the previous subsection, we need to annihilate the
projection of the $\tilde{R}$'s in \eqref{eq:exp2nd} onto the kernel
of $\mathcal{L}_r$. This  corresponds to multiplying the error terms
by $\pa_m U(k \cdot)$, $m = 1, \dots, n-1$, integrating over
$\R^{n-1}$ and taking the real part. As before, we are left to
consider only two terms: $\tilde{R}_{r,o}$ and
$\tilde{R}_{r,o,f_1}$.

We begin by multiplying $\tilde{R}_{r,o}$ by $\partial_{m} U$ and
integrating. Recalling the expression of $w_{i}$ determined in
Subsection \ref{s:exp1st} ($w_{i,o}=-\sum_j \Phi_j'f'h z_jU(kz)$ and
$w_{i,e}=\frac{p-1}{4}f'h'|z|^2U(kz)$), we obtain
$$
  \int_{\R^{n-1}} (\tilde{R}_{r,o} + {\tilde R}_{r,o,f_1}) \pa_m U(kz) = A_1 + A_2 + A_3 +
  A_4 + A_5 + A_6,
$$
where
$$
  A_1 = 2 (f')^2 \Phi''_m h \int_{\R^{n-1}} z_m U(kz) \pa_m
  U(kz) + \Phi''_m h k \int_{\R^{n-1}} (\pa_m U(kz))^2;
$$
\begin{eqnarray*}
    A_2 & = & \Phi'_m \bigg[ 2 f' f'' h \int_{\R^{n-1}}z_m U(kz) \pa_m U(kz)
    + 2 (f')^2 h' \int_{\R^{n-1}} z_m U(kz) \pa_m U(kz)
    \\ & + & f' f'' h \int_{\R^{n-1}}
    z_m U(kz) \pa_m U(kz) + 2 h' k \int (\pa_m
    U(kz))^2 +  2 h k' \int_{\R^{n-1}} (\pa_m
    U(kz))^2 \\ & + &  2 (f')^2 \frac{k'}{k} h \int_{\R^{n-1}}
    z_m U(kz) \pa_m U(kz)
      \bigg];
\end{eqnarray*}
\begin{eqnarray*}
    A_3 & = & \Phi'_m \left[ 2 (f')^2 h k'
   \int_{\R^{n-1}} z_m \pa_m U(kz) \n U(kz) \cdot z + 2 h k' k \sum_l
   \int_{\R^{n-1}} z_l \pa^2_{ml} U(kz) \pa_m U(kz) \right. \\ &
    & \left. + (f')^2 k' h \int_{\R^{n-1}} |z|^2 (\pa_m U(kz))^2 +
    \frac 12 (p-1) (f')^2 h^p \frac{k'}{k} \int_{\R^{n-1}} |z|^2 (U(kz))^p z_m
    \pa_m U(kz) \right].
\end{eqnarray*}
Here we used the fact that
$\frac{k'}{k}=\frac{p-1}{2}\frac{h'}{h}$. We have next
\begin{eqnarray*}
    A_4 & = & 2 (f')^2 H^m \int_{\R^{n-1}} z_m
    \pa_m U(kz) w_{r,e} +2(f')^2 \langle {\bf H},\Phi \rangle
    \int_{\R^{n-1}} w_{r,o} \pa_m U(kz) \\
    & + & H^m \int_{\R^{n-1}} \pa_m w_{r,e} \pa_m U(kz) -
    p (p-1) h^{p-2} \int_{\R^{n-1}} U(kz)^{p-2} w_{r,e} w_{r,o}
    \pa_m U(kz) \\ & + &
    \int_{\R^{n-1}}  \langle \n^N V,w_{r,e}z+w_{r,o}\Phi \rangle
    \pa_mU(kz);
\end{eqnarray*}
\begin{eqnarray*}
    A_5 & = & - (f')^2 h \sum_l \Phi_l \partial^2_{lm} g_{11} \int_{\R^{n-1}} z_m U(kz)
    \pa_m U(kz) + 8 (f')^2 h \langle {\bf H},\Phi \rangle H^m\int_{\R^{n-1}} z_m U(kz) \pa_m
    U(kz) \\ & + &h k \langle {\bf H},\Phi \rangle H^m
    \int_{\R^{n-1}} (\pa_m U(kz))^2 + h k^2 \sum_{l,s,t,j} \partial^2_{ls}
    g_{tj} \Phi_l \int_{\R^{n-1}} \partial^2_{tj} U(kz) z_s
    \pa_m U(kz) \\ & + & h k \sum_{l,t} \partial^2_{lt}
    g_{tm} \Phi_l \int_{\R^{n-1}} (\pa_m U(kz))^2 + h k \langle {\bf H},\Phi \rangle
    H^m \int_{\R^{n-1}} (\pa_m U(kz))^2 \\ & - & \frac 12 h k
    \left( \sum_l \partial^2_{ml} g_{11} \Phi_l \right) \int_{\R^{n-1}} (\pa_m
    U(kz))^2 + \sum_{j} \partial^2_{mj} V\Phi_j h \int_{\R^{n-1}} z_m U(kz) \pa_m
    U(kz),
\end{eqnarray*}
and
\begin{eqnarray*}
  A_6 & = & 2f'f_1' \int_{\R^{n-1}} w_{r,o} \pa_m U(k z) + 4 f'f_1'h H^m
  \int_{\R^{n-1}} U z_m \pa_m U(k z)
  \\ & - & 2 p(p-1) f' f'_1 h^{p-1} \int_{\R^{n-1}} U^{p-2} \tilde{U} w_{r,o}
  \pa_m U(k z) + 2 f' f'_1 h (\n^N V)^m \int_{\R^{n-1}}
  z_m \tilde{U} \partial_m U(k z) \\
   &  + & 2 H^m f' f'_1 h \int_{\R^{n-1}} \pa_m \tilde{U} \pa_m U(k z) + 4
   H^m (f')^2 f' f'_1 \int_{\R^{n-1}} z_m \pa_m U(k z)
h \tilde{U}.
\end{eqnarray*}
Integrating by parts and using the above relations we obtain
\begin{eqnarray*}
    A_1 & = & \Phi''_m \frac{h}{k^n} \left[ k^2 \frac{p-1}{2 \th} -
    (f')^2 \right] \int_{\R^{n-1}} U^2(z) dz = \Phi''_m \frac{h}{k^n}
    \frac{p-1}{2 \th} h^\s \left[ h^\th - \frac{2 \mathcal{A}^2 \th}{p-1}
    h^\s \right] \int_{\R^{n-1}} U^2(z) dz.
\end{eqnarray*}
Similarly, for $A_2$ after some integration by parts and some
rescaling one finds
$$
  A_2 = \frac{\Phi_m'}{k^{n-1}} \left\{ \frac{2 (hk)'}{n-1} \int_{\R^{n-1}}
  |\n U(z)|^2 dz - \frac hk \left[ \left( \frac 32 f' f'' + (f')^2 \frac{k'}{k}
  \right) + \frac{h'}{k} (f')^2 \right] \int_{\R^{n-1}} U^2(z) dz \right\}.
$$
In particular, recalling the identity \eqref{eq:poho1} one then
gets
$$
  A_2 = \frac{\Phi_m'}{k^{n-1}} \left\{ \frac{p-1}{\th} (hk)' - \left[ \frac hk
  \left( \frac 32 f' f'' + (f')^2 \frac{k'}{k}
  \right) + \frac{h'}{k} (f')^2 \right] \right\} \int_{\R^{n-1}} U^2(z) dz.
$$
Now we turn to the term $A_3$. First of all we can write
$$
  \sum_l \int_{\R^{n-1}} z_l \pa^2_{ml} U(kz) \pa_m U(kz)
  = - \frac{n-1}{k} \int_{\R^{n-1}} (\pa_m U(kz))^2 - \sum_l \int_{\R^{n-1}}
  z_l \pa^2_{ml} U(kz) \pa_m U(kz).
$$
Hence, using a simple scaling and again formula \eqref{eq:poho1}
we obtain that
$$
  \sum_l \int_{\R^{n-1}} z_l \pa^2_{ml} U(kz) \pa_m U(kz)
  = - \frac{(n-1)(p-1)}{4 \th k^n} \int_{\R^{n-1}} U^2(z) dz.
$$
Regarding the other terms in $A_3$, we write them in radial
coordinates, obtaining
$$
  3 \frac{(f')^2 h k'}{(n-1)} \o_{n-2} \int_0^\infty r^n (U'(kr))^2 dr
  + \frac 12 \frac{(p-1) (f')^2 h^p}{(n-1)} \frac{k'}{k} \o_{n-2}
  \int_0^\infty r^{n+1} U(kr)^p U'(kr) dr,
$$
where $\o_{n-2}$ is the volume of $S^{n-2}$. Using a change of
variables and integrating by parts (recalling the relation $k^2 =
h^{p-1}$) we then find
$$
  \frac{\o_{n-2} (f')^2}{n-1} \frac{k' h}{k^{n+1}} \left[ 3 \int_0^\infty
  r^n (U'(r))^2 dr - \frac 12 \frac{(n+1)(p-1)}{p+1} \int_0^\infty
  r^n U^{p+1}(r) dr \right].
$$
Using the formula in the appendix we find that this expression
becomes
$$
  \frac{n+1}{2} \frac{(f')^2 k' h}{k^{n+1}} \int_{\R^{n-1}} U^2(z)
  dz.
$$
Therefore we also obtain that
$$
  A_3 = \Phi'_m \left[ \frac{n+1}{2} \frac{(f')^2 k' h}{k^{n+1}}
  \int_{\R^{n-1}} U^2(z) dz - \frac{(n-1)(p-1)}{2 \th k^{n-1}}
  h k' \int_{\R^{n-1}} U^2(z) dz \right].
$$
Finally, after some tedious but straightforward computation one also
deduces
$$
  A_2 + A_3 = \Phi'_m \left[ \frac{p-1}{2} \frac{h'}{k^{n-2}}
  - \mathcal{A}^2 \s \frac{h^{2\s} h'}{k^n} \right] \int_{\R^{n-1}} U^2 = \Phi'_m \frac{p-1}{2}
  \frac{h'}{k^n} \left[ h^{p-1} - \frac{2 \mathcal{A}^2 \s}{p-1} h^{2\s}
  \right] \int_{\R^{n-1}} U^2.
$$
It remains now to evaluate the contribution of $A_4$ and of $A_5$.
Regarding $A_4$ we recall that $w_r$ satisfies
$$
  \mathcal{L}_r w_r = F,
$$
where
$$
  F = (f')^2 h U(kz) \left( \sum_m \partial_m g_{11}
    (z_m + \Phi_m) \right) + \frac 12 \sum_j \partial_j g_{11}
    h k \partial_j U(kz) - \sum_m \frac{\partial V}{\partial z_m}
    (z_m + \Phi_m) h U(kz).
$$
Differentiating this relation with respect to $z_m$ we obtain
$$
  \mathcal{L}_r \pa_m w_r = \pa_m F + p(p-1) h^{p-1} k
  U(kz)^{p-2} \pa_m U(kz) w_r.
$$
Multiplying by $w_r$, integrating and using the self-adjointness of
$\mathcal{L}_r$ we find
$$
  \frac 12 p(p-1) h^{p-2} \int_{\R^{n-1}} U(kz)^{p-2} \pa_m U(kz)
  w_r^2 - \frac{1}{2hk} \int_{\R^{n-1}} F \left( F \pa_m w_r -
  w_r \pa_m F \right).
$$
We have clearly
$$
  0 =\int_{\R^{n-1}} \pa_m (F w_r) =\int_{\R^{n-1}} \left( w_r \pa_m F
  + F \pa_m w_r \right).
$$
Therefore it follows that
$$
  - \frac 12 p (p-1) h^{p-2} \int_{\R^{n-1}} U(kz)^{p-2} \pa_m U(kz)
  w_r^2 = \frac{1}{hk} \int_{\R^{n-1}} w_r \pa_m F.
$$
Hence $A_4$ can be written as
\begin{eqnarray*}
    A_4 & = &  -\int_{\R^{n-1}} w_r \left[ k\sum_lH^l
  \pa^2_{ml} U(kz) + (f')^2  \pa_m U(kz)
    \left( \sum_l H^l (z_l + \Phi_l) \right) \right] \\
    & + & \frac{1}{k} \int_{\R^{n-1}} w_r \bigg[ -(f')^2  U(kz) H^m
    -  k \pa_m U(kz) \langle \n^N V,z+\Phi \rangle   - h U(kz) \langle \n^N V,E_m \rangle
    \bigg] \\ & + & \int_{\R^{n-1}} \pa_m U(kz) \left[ w_r
    \langle \n^N V,z+\Phi \rangle +2 (f')^2
    \langle {\bf H},z+\Phi \rangle w_r \right] \\ & + &h^2 k^2
     \int_{\R^{n-1}} \pa_m U(kz)\sum_l H^l  \partial_{l} w_r.
\end{eqnarray*}
After some cancelations and some integration by parts we find the
following formula
$$
  A_4 = -2k \sum_l H^l \int_{\R^{n-1}} w_r \pa^2_{ml} U(kz)
  + \frac{p-1}{2 k \th} k^2 \int_{\R^{n-1}} w_r U(kz).
$$
Using the symmetries of the integrals we find
\begin{eqnarray*}
    A_4 & = &-2 k H^m \left[ \frac{1}{n-1}\int_{\R^{n-1}} w_r \D U(kz)
    + \frac{p-1}{2\th} \int_{\R^{n-1}} w_r U(kz) \right] \\ & = &-2 k
    H^m\int_{\R^{n-1}} w_{r,e} \left[ U(kz) \left( \frac{1}{n-1} +
\frac{p-1}{2\th} \right)
    - \frac{1}{n-1} U^p(kz) \right].
\end{eqnarray*}
From the explicit expression of $w_r$ and some integration by parts
we obtain
\begin{eqnarray*}
    A_4 & = &2 k H^m \frac{(p-1) h}{ \th k^{n-1}}
    \langle {\bf H},\Phi \rangle \\ & \times &
    \left[ \frac{1}{n-1} \left( \frac{1}{p-1} - \frac{n-1}{2(p+1)}
    \right) \int_{\R^{n-1}} U^{p+1} - \left( \frac{1}{p-1} - \frac{n-1}{4}
    \right) \left( \frac{1}{n-1} + \frac{p-1}{2\th} \right) \int_{\R^{n-1}} U^{2}
    \right].
\end{eqnarray*}
Using \eqref{eq:poho1} we find that $\int_{\R^{n-1}} U^{p+1} =
\left( 1 + \frac{(n-1)(p-1)}{2 \th} \right) \int_{\R^{n-1}} U^2$,
and therefore it follows that
\begin{eqnarray*}
    A_4 & = & k H^m \frac{(p-1)^2 h}{2 \th^2 k^{n-1}}
    \langle {\bf H},\Phi \rangle  \int_{\R^{n-1}} U^2.
\end{eqnarray*}

Next we turn to $A_5$. First of all we consider the two terms
$$
  B_1 =  h k^2\sum_{l,s,t,j} \partial^2_{ls}
    g_{tj} \Phi_l \int_{\R^{n-1}} \partial^2_{tj} U(kz) z_s
    \pa_m U(kz) + h k \sum_{l,t} \partial^2_{lt}
    g_{tm} \Phi_l \int_{\R^{n-1}} (\pa_m U(kz))^2.
$$
Looking at the first one, by symmetry reasons the summands do not
vanish if, either $s = m$ and $t = j$, if $t = m$ and $s = j$ or
if $j = m$ and $s = t$. In the first case, we see appearing the
second derivative of $g_{tt}$, which vanishes by our choice of the
geodesic coordinates. Therefore we are left with the terms
$$
  B_2 = h k^2\sum_{l,j} \partial^2_{lj}
    g_{mj}  \Phi_l \int_{\R^{n-1}} \partial^2_{mj} U(kz) z_j
    \pa_m U(kz) +h k^2 \sum_{l,t} \partial^2_{lt}
    g_{tm}\Phi_l\int_{\R^{n-1}} \partial^2_{tm} U(kz) z_t
    \pa_m U(kz).
$$
Integrating by parts one easily finds
$$
  B_2 = - \frac 12h k \sum_{l,j} \partial^2_{lj}
    g_{mj}  \Phi_l \int_{\R^{n-1}} (\pa_m U(kz))^2
    - \frac 12h k \sum_{l,t} \partial^2_{lt}
    g_{tm}  \Phi_l \int_{\R^{n-1}} (\pa_m U(kz))^2,
$$
therefore it follows that
$$
  B_1 = 0.
$$
We turn now to the other terms. Integrating by parts and using
\eqref{eq:poho1} one deduces
\begin{eqnarray*}
    - (f')^2 h \sum_l \Phi_l \partial^2_{lm} g_{11} \int_{\R^{n-1}} z_m U(kz)
    \pa_m U(kz) - \frac 12 h k \left( \sum_l \partial^2_{ml}
    g_{11} \Phi_l \right)\int_{\R^{n-1}} (\pa_m U(kz))^2 \\ = \frac 12
    \left( \sum_l \Phi_l \partial^2_{lm} g_{11} \right) \frac{h}{k^n}
    \left[ \mathcal{A}^2 h^{2\s} - h^{p-1} \frac{p-1}{2\th} \right].
\end{eqnarray*}
Hence after some integration one finds
\begin{eqnarray*}
    A_5 & = & \frac 12 \left( \sum_l \Phi_l \partial^2_{lm} g_{11}
    \right) \frac{h}{k^n} \left[ \mathcal{A}^2 h^{2\s} - h^{p-1} \frac{p-1}{2\th}
    \right] \int_{\R^{n-1}} U^2 - \frac{h}{2k^n}\sum_{j} \partial^2_{mj} V \Phi_j
    \int_{\R^{n-1}} U^2 \\ & + & 4\frac{h}{k^n}\langle{\bf H},\Phi\rangle H^m
     \left[ h^{p-1} \frac{p-1}{4\th} - \mathcal{A}^2 h^{2\s}
    \right] \int_{\R^{n-1}} U^2.
\end{eqnarray*}

Finally we turn to $A_6$: we first evaluate the terms involving
$w_{r,o}$, whose explicit expression is not known, but which can
be handled via some integration by parts. Differentiating the
equation $\mathcal{L}_r \tilde{U} = - U$ with respect to $z_m$ we
find that
$$
  \mathcal{L}_r (\partial_m \tilde{U}) = - \pa_m U + p(p-1)
  h^{p-1} U^{p-2} \tilde{U} \pa_m U.
$$
Therefore, integrating by parts and recalling the definition of
$w_{r,o}$ we obtain
\begin{eqnarray*}
  & & 2f'f_1' \int_{\R^{n-1}} w_{r,o} \pa_m U(k z) - 2 p(p-1) f'
  f'_1 h^{p-1} \int_{\R^{n-1}} U^{p-2} \tilde{U} w_{r,o} \pa_m U(k
  z) \\ & = & - 2 f' f'_1 \int_{\R^{n-1}} w_{r,o} \mathcal{L}_r
  (\pa_m \tilde{U}) = - 2 f' f'_1 \int_{\R^{n-1}} \pa_m \tilde{U}
  \mathcal{L}_r w_{r,0} \\ & = & 2 f' f'_1 h H^m \int_{\R^{n-1}}
  \pa_1
  \tilde{U} \left( k \pa_1 U + \frac{p-1}{\th} h^{p-1} z_1 U
  \right).
\end{eqnarray*}
Using \eqref{eq:f'C} and \eqref{eq:euler}, we can also combine the
last term in the second row of $A_6$ with the last one in the
third row to obtain
\begin{eqnarray*}
   & & 2 f' f'_1 h (\n^N V)^m \int_{\R^{n-1}}
  z_1 \tilde{U} \partial_1 U(k z) + 4
   H^m (f')^2 f' f'_1 h \int_{\R^{n-1}} z_1 \pa_1 U(k z) \tilde{U} \\
  & = & 2 f' f'_1 h \left[ H^m \left( \frac{p-1}{\th} h^{p-1}
  - 2 \mathcal{A}^2 h^{2 \s} \right) \int_{\R^{n-1}}
  z_1 \tilde{U} \partial_1 U(k z) + 2 \mathcal{A}^2 h^{2\s}
  H^m \int_{\R^{n-1}} z_1 \pa_1 U(k z) \tilde{U} \right] \\ & = &
  2 f' f'_1 h \frac{p-1}{\th} h^{p-1} H^m \int_{\R^{n-1}}
  z_1 \tilde{U} \partial_1 U(k z).
\end{eqnarray*}
With some manipulation, the sum of the last two formulas gives
$$
 2 f' f'_1 h H^m \left( \int_{\R^{n-1}} \pa_1 \tilde{U} \pa_1 U
 - \frac{p-1}{\th} \frac{h^{p-1}}{k} \int_{\R^{n-1}} U \tilde{U}
 \right).
$$
Collecting all the terms in $A_6$ then we get
$$
 2 f' f'_1 h H^m \left( 2 k \int_{\R^{n-1}} \pa_1 \tilde{U} \pa_1 U
 - \frac{p-1}{\th} \frac{h^{p-1}}{k} \int_{\R^{n-1}} U \tilde{U}
 + 2 \int_{\R^{n-1}} U z_1 \pa_1 U(k z) \right).
$$
With some integration by parts, one finds that $$\int_{\R^{n-1}}
\pa_1 \tilde{U} \pa_1 U = \frac{1}{n-1} \int_{\R^{n-1}} \tilde{U} (-
\D U) = \frac{1}{n-1} \int_{\R^{n-1}} \tilde{U} (U^p - U),$$ and
$\int_{\R^{n-1}} U z_1 \pa_1 U(k z) = - \frac{1}{2k} \int_{\R^{n-1}}
U^2$, and therefore the last quantity becomes
$$
  2 f' f'_1 h H^m \left[ - \left( \frac{2}{n-1} + \frac{p-1}{\th} \right)
  \int_{\R^{n-1}} U \tilde{U} + \frac{2}{n-1} \int_{\R^{n-1}} \tilde{U}
  U^p - 1 \right].
$$
Using \eqref{eq:poho1} one also finds
\begin{eqnarray*}
  \int_{\R^{n-1}} \tilde{U}(kz) U^p(kz)dz &=& \frac{1}{k^{2}}
  \left( \frac{1}{p-1} - \frac{n-1}{2(p+1)}
  \right) \int_{\R^{n-1}} U^{p+1}(kz)dz \\
  &=& \frac{1}{k^{2}}\frac{\th}{(p+1)(p-1)} \frac{p+1}{\th}
  \int_{\R^{n-1}} U^2(kz) = \frac{1}{k^{2}}\frac{1}{p-1} \int_{\R^{n-1}} U^2(kz).
\end{eqnarray*}
{From} this formula, \eqref{eq:utu} and some manipulation we then
get
\begin{eqnarray*}
   & & 2 f' f'_1 h H^m \left( \frac{2 \s (p+1)}{2 (n-1) (p-1) \th} +
  \frac{2}{(p-1)(n-1)} - 1 \right) \int_{\R^{n-1}} U^2(kz) \\
  & = & 2 f' f'_1 \frac{h}{k^n} H^m \left( \frac{p-1}{2 \th} - 1 \right)
  \int_{\R^{n-1}} U^2(z).
\end{eqnarray*}
Now, recalling the expression of $f'_1$ given by
\eqref{eq:f1'final}, the total projection onto the kernel of
$\mathcal{L}_r$ is
\begin{eqnarray*}
    &  & C_0 \left\{\frac 12\frac{h}{k^n}  \left( \mathcal{A}^2 h^{2\s}
    - h^{p-1} \frac{p-1}{2\th} \right)  \sum_l \Phi_l
    \partial^2_{lm} g_{11}   -
    \frac{h}{2k^n}\sum_{j} \partial^2_{mj} V\Phi_j \right. \\ & + & 4\frac{h}{k^n}\langle{\bf H},\Phi\rangle H^m
     \left[ h^{p-1} \frac{p-1}{4\th} - \mathcal{A}^2 h^{2\s}
    \right]  +   \frac{(p-1)^2 h}{2 \th^2
    k^{n-2}}H^m \langle{\bf H},\Phi\rangle   \\
    & + & \Phi''_m \frac{h}{k^n} \frac{p-1}{2 \th} h^\s \left[ h^\th -
    \frac{2 \mathcal{A}^2 \th}{p-1} h^\s \right]  +
    \Phi'_m \frac{p-1}{2} \frac{h'}{k^n} \left[ h^{p-1} - \frac{2
    \mathcal{A}^2 \s}{p-1} h^{2\s} \right] \\
    & + & \left. \left( \frac{p-1}{2\th} - 1\right)^2\frac{h}{k^n} \frac{(2 \mathcal{A})^2
    (p-1) h^{\s+p-1}}{(p-1) h^{\th} - 2 \s \mathcal{A}^2 h^{\s}}
    \langle {\bf H}, \Phi \rangle H^m  -
    \left( \frac{\s-\th}{2\th} \right)\frac{h}{k^n} \frac{2 \mathcal{A}
    \mathcal{A}' (p-1) h^{\s+p-1}}{(p-1) h^{\th} - 2 \s
    \mathcal{A}^2 h^{\s}} H^m \right\},
\end{eqnarray*}
where $C_0=\int_{\R^{n-1}}U^2$. Now, recalling the definitions of
$\s$ and $\th$, choosing $f_1$ as in \eqref{eq:f1'final}, with some
calculation we find
\begin{eqnarray}\label{eq:2ndvarfin42} \nonumber
    \int_{\R^{n-1}} (\tilde{R}_{r,o} + \tilde{R}_{r,o,f_1}) \pa_m U(k(\e s) z) dz =
    \qquad \qquad \qquad \qquad \qquad \qquad  \\ \nonumber
   -\frac{p-1}{2\th}\frac{1}{h k} C_0\left\{- \left( h^\th -
 \frac{2 \mathcal{A}^2 \th}{p-1} h^\s \right)
 \Phi''_m - \th
  \left( h^{\th-1} - \frac{2 \mathcal{A}^2 \s}{p-1} h^{\s-1} \right) h' \Phi'_m
  + \frac{\th}{p-1} h^{-\s} ((\n^N)^2 V)\Phi_m\right.
  \\  + \frac 12 \left( h^\th - \frac{2 \mathcal{A}^2 \th}{p-1} h^\s \right) \left(
  \sum_j \left(  \partial^2_{jm} g_{11} \right) \Phi_j \right)
  - 2 \mathcal{A} \mathcal{A}'_1 \frac{(\th-\s)
  h^{p-1}}{\left[ (p-1) h^{\th} - 2 \s \mathcal{A}^2 h^{\s}
  \right]} H^m \\
 \left. + H^m \langle {\bf H}, \Phi \rangle
  \left[ \frac{- (p-1) \left( 3 + \frac{\s}{\th} \right)
  h^{2\th} - \frac{16 \s \th \mathcal{A}^4}{p-1} h^{2\s} + 2 \mathcal{A}^2 (5\s +
  3\th) h^{\th+\s}}{(p-1) h^{\th} - 2 \mathcal{A}^2 \s h^{\s}} \right]\right\}. \nonumber
\end{eqnarray}
We notice that  the operator between brackets coincides precisely
with the one in \eqref{eq:2ndvarfin4}, corresponding to the second
variation of the reduced functional which we determined in
Subsection \ref{ss:2ndvarEuler}.

\begin{rem}\label{r:phase}

According to the considerations in Subsection \ref{ss:euler}, to
every normal variation of $\g$ it corresponds some variation in
the phase due to both the variation of position and the variation
of the constant $\mathcal{A}$. Recall that the phase of the
approximate solution is the following
$$
  F_\e = \frac 1\e f(\e s) = \frac 1\e \int_0^{\e s} f' dl.
$$
Differentiating with respect to a variation $\nu$ (see
\eqref{eq:f'Cintr}) we obtain
\begin{eqnarray*}
    \frac{\partial}{\partial \nu} F_\e & = & \frac 1\e \int_0^{\e s} \mathcal{A}_\nu
    h^\s + \frac 1\e \int_0^{\e s} \mathcal{A} \s h^{\s-1} \left( \frac{\partial h}{\partial
    \mathcal{A}} \mathcal{A}_\nu + \frac{\partial h}{\partial V} \frac{\partial V}{\partial \nu}
    \right) + \frac 12 \frac 1\e \int_0^{\e s} \partial_\nu g_{11} \mathcal{A} h^\s.
\end{eqnarray*}
Recalling formula \eqref{eq:id} we find
\begin{eqnarray*}
    \frac{\partial}{\partial \nu} F_\e & = & \frac 1\e \frac{p-1}{2\mathcal{A}} \mathcal{A}_\nu
    \int_0^{\e s} h^{\th-1} \frac{\partial h}{\partial \mathcal{A}} + \frac 1\e
    \int_0^{\e s} \mathcal{A} \s h^{\s-1} \frac{\partial h}{\partial V}
    \frac{\partial V}{\partial \nu} + \frac 12 \frac 1\e \int_0^{\e s}
    \partial_\nu g_{11} \mathcal{A} h^\s.
\end{eqnarray*}
Therefore, when we take a variation $\nu_2$ of $\g$,  this also
corresponds to a variation of the phase of $\frac{\partial}{\partial
\nu_2} F_\e$. Notice that multiplying the horizontal part by $h
\pa_m U$ corresponds to adding a variation of $- \frac{\nu_2}{k}$.

Hence integrating by parts we get
$$
  \mathcal{A} \left( \frac{p-1}{2\th} - 1 \right) \int \left( \sum_m \Phi_m \partial_m
  g_{11} \right) \left[ \frac{p-1}{2\mathcal{A}} \mathcal{A}_{\nu_2} h^{\th-1} \frac{\partial
  h}{\partial \mathcal{A}} - \mathcal{A} \s h^{\s-1} \frac{\partial h}{\partial V} \frac{\partial
  V}{\partial \nu_2} - \frac 12 \mathcal{A} h^\s \partial_{\nu_2} g_{11}
  \right].
$$
\end{rem}

\subsection{Proof of Theorem \ref{t:approx}}\label{ss:proof}

The proof of Theorem \ref{t:approx} can be deduced with an iterative
procedure, adding higher order corrections (at any arbitrary order)
to the above approximate solutions. This method has been used for
other singularly perturbed equations, and is described in detail for
example in \cite{mmp}, \cite{mm}, \cite{mal}, \cite{malm2}: hence,
we will limit ourselves to a formal proof, since rigorous estimates
can be derived as in the aforementioned papers.

For $m \in \N$ we consider an approximate solution of the form
$$
  \psi_{m,\e} = e^{- i \frac{\tilde{f}(\e s)}{\e}} \left[ h(\e s) U(k(\e s )z)
  + \sum_{j=1}^m \e^j (w_{r,j}(\e s, z) + i w_{i,j}(\e s, z)) \right]; \qquad
  z = y - \Phi(\e s).
$$
Here $y$ stands for normal Fermi coordinates as in Subsection
\ref{ss:coord}, while we have set
$$
  \tilde{f}(\e s) = f + \sum_{j=1}^m \e^j f_j(\e s); \qquad \qquad
  \Phi(\e s) = \sum_{j=0}^m \e^j \Phi_j(\e s),
$$
where $f_j$, $\Phi_j$ are smooth real functions and normal sections
respectively. We then write
$$
  - \D_{\tilde{g}_\e} \psi_{m,\e} + V(\e x) \psi_{m,\e} - |\psi_{m,\e}|^{p-1}
  \psi_{m,\e} = e^{- i \frac{\tilde{f}(\e s)}{\e}} \left( \sum_{j = 0}^m \e^j
   \mathcal{R}_j \right) + o(\e^m),
$$
where $\mathcal{R}_j$ are error terms depending on $M$, $V$, $p$,
$\Phi$ and $\tilde{f}$.

In Lemma \ref{l:solve} we showed how to choose $w_{r,1}$ and
$w_{i,1}$ so that $\mathcal{R}_1$ vanishes identically. In the
previous two subsections instead we proved that $\mathcal{R}_2$ can
also be canceled provided $f_1$ satisfies \eqref{eq:f1'final} and
$\Phi_0$ (we are using the above notation) satisfies $\mathfrak{J}
\Phi_0 = 0$. By the invertibility of $\mathfrak{J}$ it is indeed
sufficient to take both $f_1$ and $\Phi_0$ identically equal to
zero.

Turning to higher order terms, we will find that the coefficient of
$\e^3$ will be of the form (up to the phase factor)
$$
  \mathcal{L}_r w_{3,r} + i \mathcal{L}_i w_{3,i} + \mathcal{G}_3(\e s, z),
$$
where $\mathcal{G}_3(\e s, z)$ is an expression depending on $V$,
$w_{1,r}, w_{1,i}$, $w_{2,r}, w_{2,i}$, $f$, $f_2$ and $\Phi_1$. As
before, we will find that the above expression can be made vanish
provided $f_2$ and $\Phi_1$ satisfy are periodic solutions of a
system of the form (see \eqref{eq:Tf1} and \eqref{eq:2ndvarfin42})
$$
   \left\{
     \begin{array}{ll}
       T f_2 = 2\mathcal{A} \left( \frac{p-1}{2\th} - 1
  \right)  \partial_{\ov s} \langle {\bf H},\Phi_1
  \rangle + \mathcal{W}_{3,1} & \hbox{ in } [0,L]; \\
      \mathfrak{J} \Phi_1 =  \mathcal{W}_{3,2} & \hbox{ in } [0,L].
     \end{array}
   \right.
$$
Here $\mathcal{W}_{3,1}$ and $\mathcal{W}_{3,2}$ are smooth
functions of $\ov{s}$ independent of $f_2$ and $\Phi_1$. Again by
our assumptions on $\g$ the above system can be solved in $f_2,
\Phi_1$, and solvability up to order $\e^3$ can be guaranteed. For
higher powers of $\e$ one  can proceed similarly.

\begin{rem}\label{orders} We stated Theorem \ref{t:approx} for
general powers in $\e$ for expository reasons. Indeed, for our
purposes in \cite{mmm2} we will need approximate solutions up to
order $\e^3$ only. However, an accurate analysis of the error terms
will be needed.
\end{rem}

\begin{rem}\label{r:sa} If we multiply the operators $\mathfrak{J}$ and
$T$ (see \eqref{eq:2ndvarfin4} and \eqref{eq:Tf1}) by $h(\ov{s})
k(\ov{s})$ and $h(\ov{s})$ respectively, they become self-adjoint.
This fact will be used crucially in the second part \cite{mmm2}, see
in particular Subsection 2.3 there.
\end{rem}

\section{Appendix}\label{s:app}

In this appendix we collect some technical results: some integral
identities first, and then the derivation of the second variation
formula \eqref{eq:2ndvarfin4}.

\subsection{Some generalized Pohozaev identities}\label{ss:apppoh}

In this Subsection we derive some identities which are useful in
the above computations. They follow from a simple integration by
parts, but the proof is rather involved so we give some
explanation here. We recall that the function $U$ in
\eqref{eq:ovv} is radial, and so it satisfies the following
ordinary differential equation
\begin{equation}\label{eq:Ur}
  - U'' - \frac{n-2}{r} U' + U = U^p \qquad \hbox{ in } \R_+.
\end{equation}
Differentiating this relation with respect to $r$ we obtain
\begin{equation}\label{eq:Ur'}
  - U''' - \frac{n-2}{r} U'' + \frac{n-2}{r^2} U' + U' = p U^{p-1}
  U' \qquad \hbox{ in } \R_+.
\end{equation}
We multiply \eqref{eq:Ur} by $r^4 U''$ and integrate by parts to get
(recall that we are in $\R^{n-1}$)
\begin{eqnarray*}
    \int_0^\infty U' r^{n+2} U''' dr & + &
    (n+2) \int_0^\infty U' U'' r^{n+1} dr
    + (n-2) \int_0^\infty U' r^{n+1} U'' dr \\ & + & (n-2) (n+1) \int_0^\infty (U')^2 r^n dr
    \\ & - & \int_0^\infty (U')^2 r^{n+2} dr - (n+2) \int_0^\infty U U' r^{n+1} dr
    \\ & = & - p \int_0^\infty (U')^2 U^{p-1} r^{n+2} dr - (n+2) \int_0^\infty
    U^p U' r^{n+1} dr.
\end{eqnarray*}
Using \eqref{eq:Ur'} the above identity simplifies as follows
\begin{equation}\label{eq:pid1}
    (n-5) \int_0^\infty r^n (U')^2 dr + (n+1) \int_0^\infty r^n
    U^2 dr = 2 \frac{n+1}{p+1} \int_0^\infty r^n U^{p+1} dr.
\end{equation}
Similarly, if we multiply \eqref{eq:Ur} by $r^3 U'$ and integrate
again by parts we find
\begin{eqnarray*}
    \int_0^\infty U r^{n+1} U''' dr + (n+1) \int_0^\infty U U'' r^{n} dr
    + (n-2) \int_0^\infty U r^{n} U'' dr + n (n-2) \int_0^\infty U U' r^{n-1}
    dr \\ - \int_0^\infty U U'
    r^{n+1}dr  - (n+1) \int_0^\infty U^2 r^{n} dr = - p \int_0^\infty U
    U' U^{p-1} r^{n+1} dr - (n+1) \int_0^\infty U^{p+1} r^{n} dr.
\end{eqnarray*}
Using \eqref{eq:Ur'} the last identity simplifies as
\begin{equation}\label{eq:pid2}
    (n-1) \int_0^\infty r^{n-2} (U)^2 dr = \int_0^\infty r^n (U')^2 dr
    + \int_0^\infty r^n U^2 dr - \int_0^\infty r^n U^{p+1} dr.
\end{equation}
{From} \eqref{eq:pid1} it then follows that
\begin{equation}\label{eq:unuup}
  \int_0^\infty r^n U^2 dr = \frac{2}{p+1} \int_0^\infty r^n U^{p+1}
  dr - \frac{(n-5)}{n+1} \int_0^\infty r^n (U')^2 dr.
\end{equation}
If we plug this identity into \eqref{eq:pid2} we get
\begin{equation}\label{eq:nuuup}
  (n-1) \int_{\R^{n-1}} U^2 dz = \frac{6}{n+1} \int_{\R^{n-1}} |z|^2
  |\n U|^2 dz - \frac{p-1}{p+1} \int_{\R^{n-1}} |z|^2 U^{p+1} dz.
\end{equation}

\subsection{Second variation}\label{s:apx}
The aim of this Subsection is to prove  formula \eqref{eq:2ndvarfin}
for the second variation of the length functional stated in the
beginning of Subsection
\ref{ss:2ndvarEuler}. \\

\no  {\em Proof of \eqref{eq:2ndvarfin}.}   Differentiating
\eqref{eq:constr} twice with respect to $t_1, t_2$ (at $(t_1, t_2)
= (0,0)$), recalling our notations in Subsection
\ref{ss:2ndvarEuler}, taking into account \eqref{eq:varl} and
\eqref{eq:2ndvarl} we obtain
\begin{equation}\label{eq:S1S2}
    \Sigma_1 + \Sigma_2 = 0,
\end{equation}
where
\begin{eqnarray*}
    \Sigma_1 & = & \int_\g \mathcal{A}'_2 \s h^{\s-1} \frac{\partial h}{\partial
    V} \langle\n^N V,\mathcal{V}  \rangle d\ov{s} + \int_\g \mathcal{A} \s (\s-1)
    h^{\s-2} \left[ \frac{\partial h}{\partial V} \langle\n^N V,\mathcal{W}\rangle
    + \frac{\partial h}{\partial \mathcal{A}} \mathcal{A}'_2 \right] \frac{\partial h}{\partial
    V} \langle\n^N V,\mathcal{V}  \rangle d\ov{s} \\ & + & \int_\g \mathcal{A} \s
    h^{\s-1} \left[ \frac{\partial^2 h}{\partial V^2}\langle\n^N V,\mathcal{W}\rangle +
    \frac{\partial^2 h}{\partial V \partial \mathcal{A}}
    \mathcal{A}'_2 \right] \langle\n^N V,\mathcal{V}  \rangle d\ov{s} +
    \int_\g \mathcal{A} \s h^{\s-1} \frac{\partial h}{\partial V}
    \left((\n^N)^2V\right)[\mathcal{V},\mathcal{W}] d\ov{s} \\
    &-& \mathcal{A}'_2 \int_\g h^\s \langle{\bf H},\mathcal{V}\rangle
    d\ov{s}  -   \int_\g \mathcal{A} \s h^{\s-1} \frac{\partial h}{\partial
    V} \langle\n^N V,\mathcal{V}  \rangle \langle{\bf H},\mathcal{W}\rangle
    d\ov{s} \\ & - & \mathcal{A} \s \int_\g h^{\s-1} \left[ \frac{\partial h}{\partial V}
    \langle\n^N V,\mathcal{W}\rangle + \frac{\partial h}{\partial \mathcal{A}}
    \mathcal{A}'_2
    \right] \langle{\bf H},\mathcal{V}\rangle d\ov{s};
\end{eqnarray*}
\begin{eqnarray*}
    \Sigma_2 & = & \int_\g \mathcal{A} h^\s
    \left[ \sum_j \dot{\mathcal{V}}^j \dot{\mathcal{W}}^j - \sum_{j,m}
    R_{1j1m} \mathcal{V}^j \mathcal{W}^m \right]
    d\ov{s} + \mathcal{A}'_1 \mathcal{A}'_2
    \s \int_\g h^{\s-1} \frac{\partial h}{\partial \mathcal{A}} d\ov{s} \\ & + &
    \mathcal{A} \mathcal{A}''_{12} \s \int_\g h^{\s-1} \frac{\partial h}{\partial
    \mathcal{A}} d\ov{s} + \mathcal{A} \mathcal{A}'_1 \s (\s-1) \int_\g
    h^{\s-2} \left[ \frac{\partial h}{\partial V}
    \langle\n^N V,\mathcal{W}  \rangle  + \frac{\partial h}{\partial \mathcal{A}}
    \mathcal{A}'_2
    \right] \frac{\partial h}{\partial \mathcal{A}} d\ov{s} \\ & + & \mathcal{A} \mathcal{A}'_1
    \s \int_\g h^{\s-1} \left[ \frac{\partial^2 h}{\partial \mathcal{A} \partial V}
    \langle\n^N V,\mathcal{W}\rangle  + \frac{\partial^2 h}{\partial \mathcal{A}^2}
    \mathcal{A}'_2 \right] d\ov{s} - \mathcal{A} \mathcal{A}'_1 \s \int_\g
    h^{\s-1} \frac{\partial h}{\partial \mathcal{A}} \langle{\bf H},\mathcal{W}\rangle
    d\ov{s}  \\ & + & \mathcal{A}'_1 \s
    \int_\g h^{\s-1} \left[ \frac{\partial h}{\partial V}
    \langle\n^N V,\mathcal{W}\rangle + \frac{\partial h}{\partial \mathcal{A}} \mathcal{A}'_2
    \right] d\ov{s} - \mathcal{A}'_1 \int_\g h^\s \langle{\bf H},\mathcal{W}\rangle d\ov{s}
    + \mathcal{A}''_{12} \int_\g h^{\s} d\ov{s}.
\end{eqnarray*}
Then the second variation of the energy is given by the following
formula at $(t_1,t_2)=(0,0)$
\begin{eqnarray*}
    \frac{\pa^2 E_\e(u_{\psi_{t_1,t_2},\mathcal{A}_{t_1,t_2}})}{\pa t_1 \pa t_2} & = &
    \int_\g \th (\th-1) h^{\th-2} \left[ \frac{\partial h}{\partial V}
    \langle \n^NV,\mathcal{W}\rangle + \frac{\partial h}{\partial \mathcal{A}}
    \mathcal{A}'_2 \right] \frac{\partial h}{\partial
    V} \langle\n^NV,\mathcal{V}\rangle d\ov{s} \\ & + & \int_0^L \th
    h^{\th-1} \frac{\partial h}{\partial V} \left((\n^N)^2V\right)[\mathcal{V},\mathcal{W}] d\ov{s} \\
    & + & \int_0^L \th h^{\th-1} \left[
    \frac{\partial^2 h}{\partial V^2} \langle\n^NV,\mathcal{W}\rangle +
    \frac{\partial^2 h}{\partial V \partial \mathcal{A}} \mathcal{A}'_2 \right]
    \langle\n^NV,\mathcal{V}\rangle d\ov{s} \\
    &-& \int_0^L \th
    h^{\th-1} \left[ \frac{\partial h}{\partial V} \langle\n^NV,\mathcal{W}\rangle
    + \frac{\partial h}{\partial \mathcal{A}} \mathcal{A}'_2
    \right] \langle{\bf H},\mathcal{V}\rangle d\ov{s}
    \\ & + & \int_0^L h^\th \left[ \sum_j \dot{\mathcal{V}}^j \dot{\mathcal{W}}^j - \sum_{j,m}
    R_{1j1m}
    \mathcal{V}^j \mathcal{W}^m  \right]
    d\ov{s} \\ & + & \mathcal{A}'_1 \int_0^L \th (\th-1) h^{\th-2} \left[
    \frac{\partial h}{\partial V} \langle \n^NV,\mathcal{W}\rangle +
    \frac{\partial h}{\partial \mathcal{A}} \mathcal{A}'_2 \right] \frac{\partial h}{\partial
    \mathcal{A}} d\ov{s}\\
    & +& \int_0^L \th \mathcal{A}'_1 h^{\th-1} \left[ \frac{\partial^2 h}{\partial \mathcal{A}
    \partial V} \langle \n^NV,\mathcal{W}\rangle + \frac{\partial^2 h}{\partial
    \mathcal{A}^2} \mathcal{A}'_2 \right] d\ov{s} \\ & + & \int_0^L
    \th \mathcal{A}''_{12} h^{\th-1}
    \frac{\partial h}{\partial \mathcal{A}} d\ov{s} - \int_0^L \left[ \th h^{\th-1}
    \frac{\partial h}{\partial V} \langle\n^N V,\mathcal{V}\rangle  + \th
    \mathcal{A}'_1
    h^{\th-1} \frac{\partial h}{\partial \mathcal{A}} \right] \langle{\bf H},\mathcal{W}\rangle.
\end{eqnarray*}
Now some cancelation will occur. We plug in the value of
$\mathcal{A}''_{12}$ from \eqref{eq:S1S2} into the last equality to
obtain
\begin{eqnarray*}
    \frac{\pa^2 E_\e(u_{\psi_{t_1,t_2},\mathcal{A}_{t_1,t_2}})}{\pa t_1 \pa t_2}
    & = & \mathcal{E}_1 + \mathcal{E}_2 + \mathcal{E}_3 + \mathcal{E}_4 +\mathcal{E}_5 +
    \mathcal{E}_6 + \mathcal{E}_7 + \mathcal{E}_8 + \mathcal{E}_9,
\end{eqnarray*}
where the terms $(\mathcal{E}_i)_i$ are given by
\begin{eqnarray*}
    \mathcal{E}_1 =  \int_0^L h^\th
    \left[ \sum_j \dot{\mathcal{V}}^j \dot{\mathcal{W}}^j - \sum_{j,m}
    R_{1j1m} \mathcal{V}^j \mathcal{W}^m  \right]
    d\ov{s} -  \frac{2 \mathcal{A}^2 \th}{p-1} \int_0^L h^\s \left[ \sum_j \dot{\mathcal{V}}^j
    \dot{\mathcal{W}}^j - \sum_{j,m} R_{1j1m} \mathcal{V}^j \mathcal{W}^m \right] d\ov{s};
\end{eqnarray*}
\begin{eqnarray*}
    \mathcal{E}_2 & = & \int_0^L \th (\th-1) h^{\th-2} \frac{\partial h}{\partial
    V} \frac{\partial h}{\partial \mathcal{A}} \left[ \langle\n^NV,\mathcal{V}\rangle
    \mathcal{A}'_2 + \langle\n^NV,\mathcal{W}\rangle \mathcal{A}'_1
    \right] \\ & + & \int_0^L \th h^{\th-1} \frac{\partial^2 h}{\partial V \partial
    \mathcal{A}} \left[ \mathcal{A}'_2 \langle\n^NV,\mathcal{V}\rangle +
    \mathcal{A}'_1 \langle\n^NV,\mathcal{W}\rangle \right]
    \\ & - & \frac{2 \mathcal{A} \th}{p-1} \int_0^L \mathcal{A} \s (\s-1)
    h^{\s-2} \frac{\partial h}{\partial \mathcal{A}} \frac{\partial h}{\partial
    V} \left[ \mathcal{A}'_1 \langle\n^NV,\mathcal{W}\rangle +
    \mathcal{A}'_2 \langle\n^NV,\mathcal{V}\rangle \right] \\ & - & \frac{2 \mathcal{A} \th}{p-1}
    \int_0^L \mathcal{A} \s h^{\s-1} \frac{\partial^2 h}{\partial \mathcal{A} \partial
    V} \left[ \mathcal{A}'_1 \langle\n^NV,\mathcal{W}\rangle + \mathcal{A}'_2
    \langle\n^NV,\mathcal{V}\rangle \right] \\ & - & \frac{2 \mathcal{A} \th}{p-1}
    \int_0^L \s h^{\s-1} \frac{\partial h}{\partial V} \left[ \mathcal{A}'_1
    \langle\n^NV,\mathcal{W}\rangle + \mathcal{A}'_2 \langle\n^NV,\mathcal{V}\rangle
    \right];
\end{eqnarray*}
\begin{eqnarray*}
    \mathcal{E}_3 & = & - \int_0^L \th h^{\th-1} \frac{\partial h}{\partial
    \mathcal{A}} \left[ \mathcal{A}'_1 \langle {\bf H}, \mathcal{W} \rangle
    + \mathcal{A}'_2 \langle {\bf H}, \mathcal{V} \rangle \right]
    + \frac{2 \mathcal{A} \th}{p-1} \int_0^L h^\s
    \left[ \mathcal{A}'_1 \langle {\bf H}, \mathcal{W} \rangle
    + \mathcal{A}'_2 \langle {\bf H}, \mathcal{V} \rangle \right] \\
    & + & \frac{2 \mathcal{A} \th}{p-1} \mathcal{A} \int_\g
    \s h^{\s-1} \frac{\partial h}{\partial \mathcal{A}} \left[ \mathcal{A}'_1 \langle {\bf H},
    \mathcal{W} \rangle + \mathcal{A}'_2 \langle {\bf H}, \mathcal{V} \rangle \right];
\end{eqnarray*}
\begin{eqnarray*}
    \mathcal{E}_4 & = & \mathcal{A}'_1 \mathcal{A}'_2 \th (\th-1) \int_0^L h^{\th-2}
    \left( \frac{\partial h}{\partial \mathcal{A}} \right)^2 + \int_0^L
    \th \mathcal{A}'_1 \mathcal{A}'_2 h^{\th-1} \frac{\partial^2 h}{\partial
    \mathcal{A}^2} \\ & - & \frac{4 \mathcal{A} \th}{p-1} \mathcal{A}'_1 \mathcal{A}'_2 \s
    \int_0^L h^{\s-1} \frac{\partial h}{\partial \mathcal{A}} - \frac{2 \mathcal{A}
    \th}{p-1} \mathcal{A} \mathcal{A}'_1 \s (\s-1) \int_0^L h^{\s-2} \left(
    \frac{\partial h}{\partial \mathcal{A}} \right)^2 \mathcal{A}'_2 \\ & - &
    \frac{2 \mathcal{A} \th}{p-1} \mathcal{A} \mathcal{A}'_1 \s \int_0^L h^{\s-1} \frac{\partial^2
    h}{\partial \mathcal{A}^2} \mathcal{A}'_2;
\end{eqnarray*}
\begin{eqnarray*}
    \mathcal{E}_5 = \int_0^L \th h^{\th-1} \frac{\partial h}{\partial
    V} ((\n^N)^2 V)[\mathcal{V},\mathcal{W}] - \frac{2 \mathcal{A}
    \th}{p-1} \mathcal{A} \s \int_0^L h^{\s-1} \frac{\partial h}{\partial
    V} ((\n^N)^2 V)[\mathcal{V},\mathcal{W}];
\end{eqnarray*}
\begin{eqnarray*}
    \mathcal{E}_6 = \int_0^L \th h^{\th-1} \left[ \frac{\partial^2 h}{\partial V^2}
    \langle \n^N V, \mathcal{V} \rangle \langle \n^N V, \mathcal{W} \rangle
    \right] - \frac{2 \mathcal{A} \th}{p-1} \int_0^L \mathcal{A}
    \s h^{\s-1} \frac{\partial^2 h}{\partial
    V^2} \langle \n^N V, \mathcal{V} \rangle \langle \n^N V, \mathcal{W} \rangle;
\end{eqnarray*}
\begin{eqnarray*}
    \mathcal{E}_7 = - \int_0^L \th h^{\th-1} \frac{\partial h}{\partial
    V} \langle \n^N V, \mathcal{W} \rangle \langle {\bf H}, \mathcal{V}
    \rangle + \frac{2 \mathcal{A} \th}{p-1} \mathcal{A} \s \int_0^L h^{\s-1} \frac{\partial
    h}{\partial V} \langle \n^N V, \mathcal{W} \rangle \langle {\bf H},
    \mathcal{V} \rangle;
\end{eqnarray*}
\begin{eqnarray*}
    \mathcal{E}_8 = - \int_0^L \th h^{\th-1} \frac{\partial h}{\partial
    V} \langle \n^N V, \mathcal{V} \rangle \langle {\bf H}, \mathcal{W}
    \rangle + \frac{2 \mathcal{A} \th}{p-1} \mathcal{A} \s \int_0^L h^{\s-1} \frac{\partial
    h}{\partial V} \langle \n^N V, \mathcal{V} \rangle \langle {\bf H},
    \mathcal{W} \rangle;
\end{eqnarray*}
\begin{eqnarray*}
    \mathcal{E}_9 = \int_0^L \left[ \th (\th-1) h^{\th-2} - \frac{2 \mathcal{A}^2 \th}{p-1} \s (\s-1)
    h^{\s-2} \right] \left( \frac{\partial h}{\partial V} \right)^2 \langle
    \n^N V, \mathcal{V} \rangle \langle \n^N V, \mathcal{W} \rangle.
\end{eqnarray*}
Now we will simplify each of these terms. We get immediately
\begin{eqnarray*}
    \mathcal{E}_1 = \int_0^L \left[ h^\th - \frac{2 \mathcal{A}^2 \th}{p-1} h^\s \right]
    \left[ \sum_j \dot{V}^j \dot{W}^j - \sum_{j,m} R_{1j1m} V^j W^m  \right]
    d\ov{s}.
\end{eqnarray*}
Recall the identity \eqref{eq:id}
$$
  h^{\th-1} \frac{\partial h}{\partial \mathcal{A}} = \frac{2\mathcal{A}}{p-1} \left(
  h^\s + \mathcal{A} \s h^{\s-1} \frac{\partial h}{\partial \mathcal{A}} \right):
$$
differentiating  with respect to $V$ and $\mathcal{A}$ we get the
two formulas
$$
  (\th-1) h^{\th-2} \frac{\partial h}{\partial V} \frac{\partial h}{\partial
  \mathcal{A}} + h^{\th-1} \frac{\partial^2 h}{\partial
  \mathcal{A} \partial V} = \frac{2\mathcal{A}}{p-1}
  \left( \s h^{\s-1} \frac{\partial h}{\partial V} + \mathcal{A} \s (\s-1) h^{\s-2}
  \frac{\partial h}{\partial V} \frac{\partial h}{\partial \mathcal{A}} + \mathcal{A} \s h^{\s-1}
  \frac{\partial^2 h}{\partial \mathcal{A} \partial V} \right);
$$
\begin{eqnarray*}
    (\th-1) h^{\th-2} \left( \frac{\partial h}{\partial \mathcal{A}} \right)^2
  + h^{\th-1} \frac{\partial^2 h}{\partial \mathcal{A}^2} & = & \frac{2\mathcal{A}}{p-1}
  \left( 3 \s h^{\s-1} \frac{\partial h}{\partial \mathcal{A}} + \mathcal{A} \s (\s-1) h^{\s-2}
  \left( \frac{\partial h}{\partial \mathcal{A}} \right)^2 + \mathcal{A} \s h^{\s-1}
  \frac{\partial^2 h}{\partial \mathcal{A}^2} \right) \\ & + & \frac{2}{p-1}
  h^\s.
\end{eqnarray*}
{From} these  it follows immediately that
$$
  \mathcal{E}_2 = \mathcal{E}_3 = 0; \qquad \mathcal{E}_4 = \mathcal{A}'_1 \mathcal{A}'_2
  \frac{2 \th}{p-1} \left( \mathcal{A} \s h^{\s-1}
  \frac{\partial h}{\partial \mathcal{A}} + h^\s \right).
$$
By means of \eqref{eq:dhCV}, the above identity \eqref{eq:id} can
also be written as
$$
  h^{\th-1} 2 \mathcal{A} h^{2 \s} \frac{\partial h}{\partial V} =
  \frac{2\mathcal{A}}{p-1} \left( \mathcal{A} \s h^{\s-1}
  2 \mathcal{A} h^{2\s} \frac{\partial h}{\partial V}
  + h^\s \right),
$$
from which it follows that
  \begin{equation}\label{eq:ddhV}
h^{\th-1} \frac{\partial h}{\partial V} - \frac{2 \mathcal{A}^2
\s}{p-1}
  h^{\s-1} \frac{\partial h}{\partial V} = \frac{1}{p-1} h^{-\s}.
\end{equation}
This formula implies
$$
  \mathcal{E}_5 = \frac{\th}{p-1} \int_0^L ((\n^N)^2 V)[\mathcal{V},\mathcal{W}]
  h^{-\s};
$$
$$
  \mathcal{E}_7 = - \frac{\th}{(p-1)} \int_0^L
  \langle \n^N V, \mathcal{W} \rangle \langle {\bf H}, \mathcal{V} \rangle
  h^{-\s}; \qquad \qquad \mathcal{E}_8 = - \frac{\th}{(p-1)} \int_0^L \langle \n^N V,
  \mathcal{V} \rangle \langle {\bf H}, \mathcal{W} \rangle h^{-\s}.
$$
Therefore in conclusion we get at $(t_1,t_2)=(0,0)$
\begin{eqnarray*}
  \frac{\pa^2 E_\e(u_{\psi_{t_1,t_2},\mathcal{A}_{t_1,t_2}})}{\pa t_1 \pa t_2}  & = & \int_0^L
  \left[ h^\th - \frac{2 \mathcal{A}^2 \th}{p-1} h^\s \right] \left[ \sum_j
  \dot{\mathcal{V}}^j \dot{\mathcal{W}}^j - \sum_{j,m} R_{1j1m} \mathcal{V}^j
  \mathcal{W}^m  \right] d\ov{s} \\ &
  + & \int_0^L \left(\th h^{\th-1}- \frac{2 \mathcal{A}^2 \th}{p-1}\s h^{\s-1} \right)
  \frac{\partial^2 h}{\partial V^2}
    \langle \n^N V,\mathcal{V} \rangle\langle \n^N V,\mathcal{W} \rangle
      \\ & + & \frac{\th}{p-1} \int_0^L \bigg\{((\n^N)^2 V)[\mathcal{V},\mathcal{W}]
  -\langle \n^N V, \mathcal{W} \rangle \langle {\bf H}, \mathcal{V}\rangle
  -  \langle \n^N V,
  \mathcal{V} \rangle \langle {\bf H}, \mathcal{W} \rangle \bigg\}h^{-\s} \\ & + &
   \int_0^L \left[ \th (\th-1) h^{\th-2} - \frac{2 \mathcal{A}^2 \th}{p-1} \s (\s-1)
    h^{\s-2} \right] \left( \frac{\partial h}{\partial V} \right)^2 \langle
    \n^N V, \mathcal{V} \rangle \langle \n^N V, \mathcal{W} \rangle
    \\ & + & \mathcal{A}'_1 \mathcal{A}'_2 \frac{2 \th}{p-1}
    \int_0^L \left( \mathcal{A} \s h^{\s-1}
  \frac{\partial h}{\partial \mathcal{A}} + h^\s \right).
\end{eqnarray*}
Differentiating \eqref{eq:ddhV} with respect to $V$, we obtain
\begin{eqnarray*}
    \left[ (p-1) h^{\th-1} - 2 \s \mathcal{A}^2 h^{\s-1} \right] \frac{\partial^2 h}{\partial
    V^2} + \left[ (p-1) (\th-1) h^{\th-2} - 2 \s (\s-1) \mathcal{A}^2 h^{\s-2}
    \right] \left( \frac{\partial h}{\partial V} \right)^2  = - \s
    h^{-\s-1} \frac{\partial h}{\partial V},
\end{eqnarray*}
which yields
\begin{eqnarray*}
    \int_0^L \left(\th h^{\th-1}- \frac{2
    \mathcal{A}^2 \th}{p-1}\s h^{\s-1} \right) \frac{\partial^2 h}{\partial V^2}
    \langle \n^N V,\mathcal{V} \rangle\langle \n^N V,\mathcal{W} \rangle
      = - \frac{\th \s}{p-1} \int_0^L\frac{\partial
    h}{\partial V} \langle \n^N V,\mathcal{V}
    \rangle\langle \n^N V,\mathcal{W} \rangle h^{-\s-1} \\ + \int_0^L\left(\frac{2 \s (\s-1)
 \mathcal{A}^2\th}{p-1}h^{\s-2}- \th (\th-1) h^{\th-2}\right)  \left( \frac{\partial
    h}{\partial V} \right)^2 \langle \n^N V,\mathcal{V} \rangle\langle \n^N V,\mathcal{W} \rangle.
\end{eqnarray*}
Collecting the above computation together with some further
cancelation, one finds, at $(t_1,t_2)=(0,0)$
\begin{eqnarray*}
  \frac{\pa^2 E_\e(u_{\psi_{t_1,t_2},\mathcal{A}_{t_1,t_2}})}{\pa t_1 \pa t_2}
  & = & \int_0^L \left[ h^\th - \frac{2 \mathcal{A}^2 \th}{p-1} h^\s \right] \left[ \sum_j
  \dot{\mathcal{V}}^j \dot{\mathcal{W}}^j - \sum_{j,m} R_{1j1m}
  \mathcal{V}^j \mathcal{W}^m  \right] d\ov{s}
   \\
  & + & \frac{\th}{p-1} \int_0^L \left\{ ((\n^N)^2 V) [\mathcal{V},\mathcal{W}]
  - \langle \n^N V, \mathcal{V} \rangle
  \langle {\bf H}, \mathcal{W} \rangle - \langle \n^N V, \mathcal{W} \rangle
  \langle {\bf H}, \mathcal{V} \rangle \right\} h^{-\s} d\ov{s} \\
    & - & \frac{\s \th}{p-1} \int_0^L h^{-\s-1} \frac{\pa h}{\pa V} \langle \n^N V,
    \mathcal{V} \rangle \langle \n^N V, \mathcal{W} \rangle d\ov{s}
     \\ & + & \mathcal{A}'_1 \mathcal{A}'_2 \frac{2 \th}{p-1}
     \int_0^L \left( \mathcal{A} \s h^{\s-1}
  \frac{\partial h}{\partial \mathcal{A}} + h^\s \right) d\ov{s}.
\end{eqnarray*}
This proves formula \eqref{eq:2ndvarfin}. \hfill{$\blacksquare$}

\begin{center}
{\bf Acknowledgments}
\end{center}

\noindent This work started in the fall of 2004, when the second
author visited Departamento de Matem\'atica in Universidade Estadual
de Campinas: he is very grateful to the institution for the kind
hospitality. F. M and A. M are  supported by M.U.R.S.T within the
PRIN 2006 {\em Variational methods and nonlinear differential
equations}. F. M is grateful to Sissa for the kind hospitality.
 M. M. was supported by FAEP-UNICAMP, FAPESP and CNPq.

\

\end{document}